\documentclass[12pt]{article}
\usepackage{graphics}
\usepackage{amsmath,amsthm,amsfonts,amssymb,amscd} 
\usepackage[all]{xy}
\usepackage{pstricks}
\usepackage{tikz}
\usetikzlibrary{patterns}

\usepackage[french]{babel}

\newtheorem{theorem}{Th\'eor\`eme}[section]
\newtheorem{lemma}[theorem]{Lemme}
\newtheorem{corollary}[theorem]{Corollaire}
\newtheorem{proposition}[theorem]{Proposition}

\newtheorem{definis}[theorem]{D\'efinitions}
\newenvironment{definitions}{\begin{definis} \em}{\end{definis}}
\newtheorem{defini}[theorem]{D\'efinition}
\newenvironment{definition}{\begin{defini} \em}{\end{defini}}
\newtheorem{exa}[theorem]{Exemple}
\newenvironment{example}{\begin{exa} \em}{\end{exa}}
\newtheorem{rema}[theorem]{Remarque}
\newenvironment{remarque}{\begin{rema} \em}{\end{rema}}
\newtheorem{remas}[theorem]{Remarques}
\newenvironment{remarques}{\begin{remas} \em}{\end{remas}}
\newtheorem{construc}[theorem]{Construction}

\newenvironment{preuve}{{\noindent \small \sc Preuve.} }{\mbox{ }\hfill$\Box$  
                                                                \vspace{1.5ex} \par}

\def\rond{\mathaccent"7017}  

\def\omg{\omega}       
\def\sg{\sigma}
\def\Rho{{\rm P}}
\def\R{{\bf R}}

\def\vol{{\rm vol}}
\def\san{sous-analytique }

\def\libel#1{\label{#1}}

\def\hfleche#1#2{\smash{\mathop{\vbox{\hbox to 8.5mm{{#1}}}}\limits^{#2}}}
\def\vfleche#1#2{\bigg #1 \rlap{$\vcenter{\hbox{$\scriptstyle#2$}}$}}

\begin{document}

\title{Formes de Whitney et primitives relatives de formes diff\'erentielles sous-analytiques}

\author{Jean-Paul Brasselet\footnote{CNRS, IML Luminy} \ 
et Bernard Teissier\footnote{CNRS, IMJ Paris\newline  AMS classification: 32B20, 32B25, 58C35}}
\maketitle
{\sl (version 50)\rm} 
\bigskip
\hfill{\it D\'edi\'e \`a Heisuke Hironaka \rm} 
\noindent

\tableofcontents

\newpage
\vglue 5 truecm

\listoffigures
\newpage

\setcounter{section}{-1}  

\section{Introduction}
\setcounter{equation}{0}  
L'origine de ce travail est une question pos\'ee par Fran{\c c}ois Treves au second
auteur en 1982 : 

{\it Etant donn\'es un morphisme analytique $g\colon {\bf S}^n\to
{\bf R}$ o\`u ${\bf S}^n$ est la sph\`ere de dimension $n$, et une $r$-forme
diff\'erentielle $\omg$ de classe ${\cal C}^\infty$ sur ${\bf S}^n$ dont la
restriction \`a chaque fibre non singuli\`ere de $g$ est exacte, existe-t-il une
$(r-1)$-forme $H$ h\"old\'erienne sur  ${\bf S}^n$ telle que l'on ait  
$$dg \wedge (\omg-dH)=0 ,$$ 
la diff\'erentielle $dH$ \'etant prise au sens des distributions ?}
\smallskip

\noindent Dans ce travail nous d\'emontrons un analogue de ce r\'esultat pour une forme $\omg$
sous-analytique et continue, dans le cadre plus g\'en\'eral des morphismes sous
analytiques propres triangulables entre espaces non singuliers (voir le Corollaire \ref{THEODG} ci-dessous). Rappelons que Masahiro Shiota a d\'emontr\'e que tout morphisme sous-analytique d'un espace compact vers $\R$ est triangulable (voir \cite{Sh1}, \cite{Sh2}, Chap. II, \S 3) et que le second auteur a prouv\'e (voir \cite{Te1}) que tout morphisme sous-analytique propre devient triangulable localement sur la base apr\`es des changements de base qui sont des compos\'es finis d'\'eclatements locaux du but.\par\noindent
S. Chanillo et F. Treves ont prouv\'e en 1997 (voir \cite{C-T}, Lemme 2.2) un
r\'esultat analogue \`a l'\'enonc\'e ci-dessus: sous l'hypoth\`ese d'annulation \`a l'ordre infini de $\omg$ le
long des fibres singuli\`eres de $g$, ils obtiennent une primitive relative de classe
${\cal C}^\infty$ satisfaisant la m\^eme propri\'et\'e d'annulation.\par   Rappelons aussi
qu'une fonction $h$ sur une vari\'et\'e analytique
$U$ est dite h\"old\'erienne si tout point  de $U$  poss\`ede un voisinage $V$ tel
qu'il existe des constantes positives $\alpha$ et $C$ telles que pour $x,x'\in V$ on
ait $\vert h(x)-h(x')\vert \leq C\vert x-x'\vert ^\alpha $. Une forme
diff\'erentielle sur  $U$  est dite h\"old\'erienne si ses coefficients sont des
fonctions  h\"old\'eriennes sur $U$  et si sa
diff\'erentielle au sens des distributions peut \^etre repr\'esent\'ee par une forme \`a coefficients h\"old\'eriens. Rappelons enfin que pour une fonction, \^etre sous-analytique et continue implique la propri\'et\'e de H\"older  gr\^ace aux in\'egalit\'es de \L ojasiewicz \'etendues au cas sous-analytique
par Hironaka ({\it cf} \cite{Hi}, \cite{Ha} et \S 1).

L'id\'ee principale est de transformer ce probl\`eme d'analyse en un
probl\`eme g\'eom\'etrique en repr\'esentant les formes diff\'erentielles au moyen
des {\it formes de Whitney}. Rappelons que la forme diff\'erentielle que Whitney
associe ({\it cf} \cite[Chap. IV, \S 27]{Whi}) \`a un simplexe d'une vari\'et\'e
diff\'erentielle triangul\'ee est essentiellement une r\'egularisation de la forme
volume de ce simplexe exprim\'ee en  coordonn\'ees barycentriques; c'est cette
derni\`ere que nous appellerons forme de Whitney. L'hypoth\`ese de
triangulabilit\'e du morphisme permet de transporter le probl\`eme sur un morphisme
simplicial $f\colon \Delta \to T$ au moyen d'hom\'eo\-mor\-phis\-mes sous-analytiques.

Le cadre sous-analytique convient bien pour cette raison et parce que les formes
volumes des simplexes en coordonn\'ees barycentriques sont lin\'eaires par morceaux (dans l'\'etoile de chaque simplexe) et continues et donc sous-analytiques et continues. Elles sont d'ailleurs non seulement h\"old\'eriennes mais m\^eme lipschitziennes comme
formes diff\'e\-ren\-tiel\-les d\'efinies sur l'espace entier. Il n'y a donc pas \`a les r\'egulariser, et l'on a une
correspondance biunivoque entre simplexes et formes de Whitney, ce qui joue un
r\^ole important dans la suite. 

La premi\`ere difficult\'e est que la d\'ecomposition simpliciale de $\Delta$
n'induit pas sur les fibres de $f$ une d\'ecomposition simpliciale. Mais nous avons
remarqu\'e qu'elle induisait une d\'ecomposition en prismes (produits de
simplexes), et que ces d\'ecompositions \'etaient naturellement isomorphes pour
toutes les fibres de points de l'int\'erieur d'un simplexe de la base $T$. L'utilit\'e de cette g\'eom\'etrie prismale a \'et\'e remarqu\'ee ind\'ependamment et pour des raisons analogues (int\'egration dans les fibres), par J. Dupont et ses collaborateurs (voir \cite{Du1}, en particulier \S 3, et \cite{Du2}). Pour
rendre compte des d\'eg\'enerescences de la structure prismale des fibres qui se
produisent par sp\'ecialisation sur les faces d'un simplexe de $T$ nous avons
introduit les {\it faisceaux prismaux}. Le morphisme simplicial d\'efinit
canoniquement un faisceau prismal mais en fait il est bien plus commode de
travailler avec un autre faisceau prismal, qui est un ``\'eclat\'e" de celui-ci et
est trivialis\'e au dessus de chaque simplexe ferm\'e de $T$. Sur ce faisceau
trivialis\'e il est facile de d\'efinir les formes de Whitney; la forme de Whitney
d'un prisme est essentiellement le produit ext\'erieur des formes de Whitney des
simplexes facteurs du prisme. Le calcul permet de d\'efinir naturellement les
formes de Whitney {\it relatives}. Apr\`es avoir  repr\'esent\'e (modulo les
formes exactes) la forme $\omg$ par une combinaison lin\'eaire de formes de Whitney
relatives, on cherche $H$ sous la m\^eme forme, et le probl\`eme de la recherche d'une primitive relative se ram\`ene
finalement \`a la r\'esolution dans chaque simplexe d'\'equations aux d\'eriv\'ees partielles d'une forme tr\`es particuli\`ere qui permet de d\'emontrer assez facilement l'existence de solutions sous-analytiques.\par
La premi\`ere partie du texte (\S1 et \S2) pr\'ecise le point de vue que nous utilisons sur la combinatoire des morphismes simpliciaux. La seconde partie (\S 3) d\'etaille le comportement des formes de Whitney vis \`a vis de cette combinatoire. On y trouve en particulier une remarque (\ref{bord}) que nous pensons int\'eressante sur la version "formes de Whitney" de l'homomorphisme bord des cha\^\i nes simpliciales. Nous avons d\'etaill\'e les calculs dans l'espoir de faciliter l'application de cette g\'eom\'etrie prismale \`a d'autres probl\`emes. Dans la troisi\`eme partie (\S 4 et \S 5) nous exprimons explicitement une forme sous-analytique comme combinaison \`a coefficients sous-analytiques de formes de Whitney. Dans la derni\`ere partie, les r\'esultats des parties pr\'ec\'edentes permettent de poser de mani\`ere pr\'ecise les \'equations \`a r\'esoudre pour trouver une primitive relative; cette partie contient les r\'esultats d'int\'egration locale et de globalisation qui ach\`event la preuve du th\'eor\`eme \ref{THEOPR} et de son corollaire \ref{THEODG}, 
qui est notre r\'esultat principal.\par\medskip
\textit{Nous remercions les coll\`egues qui nous ont encourag\'es pendant la pr\'eparation de ce texte, et principalement L\^e D.T. et Claude Weber, ainsi que le D\'epartement de Math\'e\-ma\-ti\-ques de l'Universit\'e de Gen\`eve et la School of Mathematics de l'ICTP de Trieste qui nous ont offert la possibilit\'e de venir y travailler \`a ce projet.}

\section{G\'eom\'etrie des morphismes triangul\'es}  
\setcounter{equation}{0}

\subsection{Rappels}\label{debut}

\begin{definition} \libel{ANALYT} Soit $X \subset {\bf R}^n
\times {\bf R}^m $ un sous-ensemble sous-analytique tel que la restriction $g \colon X
\rightarrow {\bf R}^n$ de la premi\`ere projection soit un morphisme propre. On
rappelle qu'une triangulation de cette situation est la donn\'ee~: 
\begin{enumerate}
\item d'hom\'eomorphismes sous-analytiques $t$ de ${\bf R}^n \times {\bf R}^m $
dans lui-m\^eme et $t_0$ de ${\bf R}^n$ dans lui-m\^eme, tels que le
diagramme suivant soit commutatif 
$$
\begin{matrix}
{\bf R}^n \times {\bf R}^m  & \hfleche\rightarrowfill{t}& {\bf R}^n
\times {\bf R}^m&\supset & X \cr 
\vfleche\downarrow{pr_1}&&\vfleche\downarrow{pr_1}&&\vfleche\downarrow{g}\cr 
{\bf R}^n&\hfleche\rightarrowfill{t _0} &{\bf R}^n&=&{\bf R}^n\cr 
\end{matrix}$$ 
\item de d\'ecompositions de ${\bf R}^n \times {\bf R}^m$ et de ${\bf R}^n$ en
simplexes lin\'eaires de telle fa\c con que les hom\'eomorphismes $t$ et $t_0$
soient analytiques \`a l'int\'erieur de chaque simplexe et que $pr_1$ soit une
application simpliciale. Ces donn\'ees  doivent \^etre compatibles avec $X$~: si
l'image par $t$ d'un simplexe  rencontre $X$, elle est contenue dans $X$. 
\end{enumerate}
\end{definition}

Dans la suite, \'etant donn\'e un morphisme $g$ triangul\'e comme ci-dessus, nous
noterons $\Delta$ le sous-complexe simplicial de ${\bf R}^n \times {\bf R}^m$ form\'e
des simplexes dont l'image par $t$ rencontre $X$. Le sous-complexe simplicial de
${\bf R}^n$ form\'e des simplexes dont l'image par $t_0$ rencontre $g(X)$ sera not\'e $T$, enfin $f\colon
\Delta \to T$ d\'esignera le morphisme induit par $\hbox{\rm pr}_1$.

\begin{remarques} 
1) D'apr\`es l'in\'egalit\'e de \L ojasiewicz, tout hom\'eomorphisme
sous-ana\-ly\-ti\-que $t$ est h\"old\'erien sur tout compact: pour tout sous-ensemble
compact $K$ de $\bf R \rm ^n \times \bf R \rm ^m $, il existe des nombres r\'eels
positifs $C$ et $\alpha$ tels que l'on ait :  $\Vert t (y) -t (x) \Vert \leq C \Vert
y -x \Vert^{\alpha}$. \par  
\noindent 2) L'image par $pr_1$ d'un simplexe est un simplexe et l'image 
r\'eciproque par $pr_1$ d'un simplexe est une r\'eunion de simplexes. 
\end{remarques}

\subsection{Morphismes et faisceaux prismaux }

\begin{definition} \libel{PRISME} Un {\sl prisme} (resp. un prisme
ouvert, resp. un prisme ferm\'e) est un produit de simplexes (resp. simplexes ouverts, resp.
ferm\'es) lin\'eairement plong\'es dans un espace euclidien et muni de la topologie induite ; l'ensemble vide est donc un prisme 
auquel, par convention, on attribue la dimension $-\infty$.
Une {\sl face} d'un prisme est vide ou est un produit de faces de ses facteurs. Nous
appellerons {\sl coordonn\'ees barycentriques} d'un point d'un prisme la famille des
coordonn\'ees barycentriques de ses projections sur les facteurs du prisme. 
\end{definition}

Un {\sl ensemble prismal}
de $\bf R \rm^N$ est une r\'eunion localement finie de prismes ferm\'es telle que
l'intersection de deux quelconques d'entre eux soit une face de chacun. Un sous-complexe
simplicial d'une triangulation de $\bf R \rm^N$ est un ensemble prismal. \par 
Sauf mention du contraire, les simplexes et les prismes consid\'er\'es dans la suite
sont ferm\'es.

Un {\sl ensemble prismal de dimension pure}
$d$  est une r\'eunion de prismes de dimension $d$.

Une {\sl orientation} d'un prisme
$\pi=\sigma_0 \times \sigma_1 \times \cdots \times \sigma_p$ est la donn\'ee d'un
ordre sur l'ensemble $\lbrace 0,\ldots ,p\rbrace$ et d'une orientation sur chaque
simplexe $\sg_i$. Une orientation d'un prisme d\'etermine une orientation de chacune de ses faces. 
Une {\sl orientation locale} d'un ensemble prismal est la
donn\'ee d'une orientation de chacun de ses prismes.

Si un ensemble prismal de dimension pure $d$ est une pseudovari\'et\'e 
(resp. une pseudovari\'et\'e \`a bord, voir \cite[Chap. 3, 11]{Sp}), tout prisme de
dimension $d-1$  est face de deux prismes (resp. un ou deux prismes) de dimension $d$. 
Le bord, quand il existe, est le sous ensemble prismal constitu\'e des prismes de dimension $d-1$ qui sont faces 
d'exactement un prisme de dimension $d$. \libel{abord} 
Un ensemble prismal $\Pi$ qui est une pseudovari\'et\'e de dimension
$d$ est {\sl orient\'e} s'il est muni d'une orientation locale telle que les orientations induites sur tout
prisme de dimension $d-1$ par deux prismes dont il est face soient oppos\'ees.

Rappelons qu'un morphisme simplicial d'un complexe simplicial dans
un autre est une application dont la restriction \`a chaque simplexe de la source a
pour image un simplexe du but et qui est lin\'eaire en les coordonn\'ees
barycentriques.

\begin{definition} \libel{APPLPR}
Une application $\varphi$ d'un ensemble prismal $\Pi$ dans un ensemble prismal
$\Sigma$ est un morphisme prismal si pour tout prisme $\pi =\sigma_0 \times \sigma_1
\times \cdots \times \sigma_p$ de $\Pi$, l'image $\varphi (\pi )$ est un prisme
$\tau_0 \times \tau_1 \times \cdots \times \tau_s$ et l'application induite $\varphi$ sur 
$\pi$ est lin\'eaire dans les coordon\'ees barycentriques des prismes. 
\end{definition}

Un morphisme simplicial entre complexes simpliciaux est un morphisme prismal.

Un morphisme prismal entre des ensembles prismaux orient\'es est dit orient\'e s'il
res\-pecte les orientations des prismes. 

\section{Faisceau prismal associ\'e \`a un morphisme simplicial}\label{section2}

L'exemple suivant est fondamental pour ce qui suit :

Soit $f\colon \sg \to \tau$ un morphisme simplicial entre deux simplexes ; pour chaque
sommet $y_j$ de $\tau$, ($0\leq j\leq s =\hbox{dim}\tau$), posons
$\sg_j=f^{-1}(y_j)$ et 
$\pi_f(\sg)=\tau\times\sg_0\times\cdots\times\sg_s$. 
On se propose de d\'efinir des
morphismes prismaux  
$$ \psi_f^\sg : \pi_f(\sg) \to \sg \quad \hbox{ et } \quad \theta_f^\sg :
 f^{-1}(\rond \tau) \to  \pi_f(\sg)$$ 
dont les restrictions au dessus de $f^{-1}(\rond \tau)$ et $\rond\tau\times\sg_0\times\cdots\times\sg_s$
sont des isomorphismes inverses l'un de l'autre. 

Lorsqu'il n'y a pas d'ambigu\"\i t\'e sur le morphisme $f$, nous omettrons l'\'ecriture 
des indices $f$ dans ce qui suit. Ainsi nous noterons
\begin{equation}\label{pidesigma}
\pi(\sg)=\tau\times\sg_0\times\cdots\times\sg_s.
\end{equation}

Chaque sommet $a_i$ de $\sg$ est dans l'un des simplexes $\sg_j$ et un seul. Pour chaque $j$, notons
$I(j)$ l'ensemble des indices $i$ des sommets $a_i$ de $\sg_j$. Etant
donn\'e un point $x$ de $\sg$, de coordonn\'ees barycentriques $\lambda_i$ 
relativement aux sommets $a_i$, les coordonn\'ees de $f(x)$ relativement aux sommets $y_j$ de $\tau$ sont
\begin{equation}\label{coordotau}
t_j = \sum_{i \in I(j)} \lambda_i.
\end{equation}
Si l'image $f(x)$ est dans l'int\'erieur de $\tau$, les sommes $\sum_{i \in I(j)}
\lambda_i$ sont donc non nulles ;  pour chaque $j$, notons $x_j$ le point de $\sg_j$ de  
coordonn\'ees barycentriques 
\begin{equation}\label{decoord}
\mu_{j,k} =\frac{\lambda_k}{\sum_{i \in I(j)}
\lambda_i}, \qquad k\in I(j)
\end{equation}
relativement aux sommets $a_k$ de $\sg$ situ\'es dans $\sg_j$. 
Notons $\theta^{\sg_j}$ l'application $f^{-1}(\rond
\tau)\to \sg_j$ qui \`a $x$ associe le point $x_j$, et $\theta^\sg$ l'application de 
$f^{-1}(\rond \tau)$ dans le prisme $\pi(\sg)=\tau\times\sg_0\times\cdots\times\sg_s$ d\'efinie par
\begin{equation}\label{detheta}
\theta^\sg \colon x \mapsto (f(x), x_0,\cdots, x_s).
\end{equation} La projection $f$ induit un isomorphisme 
simplical du simplexe enveloppe
convexe des points $x_j$  dans $\sg$ sur $\tau$. Le point $x$ est le point de ce
simplexe   dont les coordonn\'ees barycentriques sont celles du point $f(x)$ dans
$\tau$. Ceci nous donne une description de $\sg$ comme {\it joint it\'er\'e} des
simplexes $\sg_j$.

On peut en effet d\'efinir le joint it\'er\'e
$\sg_0*\cdots * \sg_s$ de $s+1$ simplexes $\sg_j$ lin\'eairement ind\'ependants dans
un espace euclidien comme la r\'eunion des simplexes de dimension $s+1$ qui sont les
enveloppes convexes d'ensembles de points de la forme $(x_j \in \sg_j);\ \
j=0,\ldots , s$. Le simplexe obtenu est simplicialement isomorphe au r\'esultat de la construction classique 
du joint it\'er\'e comme  quotient de $\sg_0 \times \ldots \times \sg_s \times [0,1]^{s}$. 
Rappelons cette construction pour $s=1$ :
$$\sg_0* \sg_1 = \left( \sg_0 \times \sg_1 \times [0,1] \right) / \langle (a,b,0) \sim (a',b,0), (a,b,1) \sim (a,b',1) \rangle.$$

Par construction, l'ensemble des sommets du joint co\"\i ncide avec l'ensemble \libel{reunion} des 
sommets des $\sg_j$. Une orientation du joint peut donc s'interpr\'eter comme la donn\'ee d'un ordre
sur l'ensemble d'indices $j$ et d'une orientation de chacun des simplexes $\sg_j$. 

Si $f$ est un morphisme simplicial orient\'e et si on munit chaque $\sg_j$ de
l'orientation induite par celle de $\sg$, alors l'orientation naturelle du joint it\'er\'e d\'eduite
de celles de $\tau$ et des $\sg_j$ n'est autre que celle de $\sg$.

On d\'efinit le morphisme $\psi^\sg \colon
\pi (\sg )\to \sg= \sg_0*\cdots * \sg_s$ qui, au point de coordonn\'ees 
$$(t_j)_{j=0,\ldots,s},(\mu_{0,i_0})_{{i_0}\in I(0)},\ldots, 
(\mu_{s,i_s})_{{i_s}\in I(s)}$$  de $\pi(\sg)$, associe le point de $\sg$ dont la
coordonn\'ee barycentrique relative au sommet $a_i$ est $\lambda_i=t_j\mu_{j,i}$,
o\`u $j$ est l'indice tel que $i\in I(j)$. Ce morphisme est un morphisme prismal.

Remarquons que, si $\sg$ et $\tau$ sont orient\'es ainsi que le morphisme $f$, 
cela d\'etermine une unique orientation du joint  it\'er\'e $\sg_0*\cdots * \sg_s$ et une unique orientation du
prisme $\pi(\sg)$ telles que les morphismes $\pi(\sg) \to \sg $ et $\pi(\sg) \to \tau $ soient orient\'es.

\begin{proposition}\libel{JOINT}
Les restrictions au dessus de $f^{-1}(\rond \tau)$ et $\rond\tau\times\sg_0\times\cdots\times\sg_s$
des morphismes prismaux $ \psi^\sg : \pi(\sg) \to \sg $  et $ \theta^\sg :
 f^{-1}(\rond \tau) \to  \pi(\sg)$
sont des isomorphismes inverses l'un de l'autre.
\end{proposition}
 
\begin{preuve}
En effet tout point $x$ de $\sg$ s'\'ecrit 
$$x = \sum_{j=0}^s t_j x_j = \sum_{j=0}^s t_j \bigg(\sum_{i\in I(j)} \mu_{j,i} a_i \bigg) = \sum_{i=0}^p \lambda_i a_i.$$
\end{preuve} 

\begin{definition} \libel{FAISPR} Soit $\hbox{\rm P}$ un
ensemble prismal. Un {\it faisceau prismal} $\cal F$
sur  $\hbox{\rm P}$ est la donn\'ee pour chaque prisme ferm\'e $\rho$ de $\Rho$ d'un 
ensemble prismal ${\cal F}(\rho)$ dot\'e d'un 
morphisme prismal $ e_\rho \colon {\cal F}(\rho) \rightarrow \rho$ et pour chaque
face $\rho'$ de $\rho$ d'un morphisme prismal
$h_{\rho',\rho}\colon {\cal F}(\rho) \rightarrow {\cal F}(\rho')$ de telle fa\c con
que $h_{\rho,\rho}=Id_{{\cal
F}(\rho)}$ et que  si $\rho''$ est une face de $\rho'$, alors  $h_{\rho'',\rho}
=h_{\rho'',\rho'}\circ h_{\rho',\rho}$. On dit que le faisceau prismal ${\cal F}$ est propre si les ensembles
prismaux ${\cal F}(\rho)$ sont compacts.
\end{definition} 

\begin{remarque} En fait, la notion de  faisceau prismal s'apparente davantage \`a la notion
de carapace (ce qui fut la premi\`ere d\'efinition 
des faisceaux, voir S\'eminaire Cartan \cite{Car2}) qu'\`a la notion
de faisceau. Plus pr\'ecis\'ement, munissons l'ensemble des prismes d'un ensemble prismal 
de l'ordre partiel donn\'e par les inclusions des faces, puis de
 la topologie  dont une base de ferm\'es est constitu\'ee des intervalles ferm\'es $\pi_1 \le \pi \le \pi_2$. 
Les sections du faisceau prismal au dessus d'un ferm\'e $F$ sont les \'el\'ements du produit 
$\Pi_{\rho \in F} {\cal F}(\rho)$ compatibles avec  les homomorphismes de restriction.
\end{remarque}

Soient $\varphi \colon \Rho \rightarrow \Sigma$ un morphisme prismal, $\cal
F$ un faisceau prismal sur $\Rho$ et $\cal G$ un faisceau prismal sur  $\Sigma$. Un
{\it morphisme de faisceaux prismaux} de $\cal F$ dans $\cal G$ est la donn\'ee pour
chaque prisme ferm\'e $\rho$ de $P$ d'un morphisme prismal ${\cal F} (\rho) \rightarrow
{\cal G}(\varphi (\rho))$ compatibles avec  les homomorphismes de restriction.  

\begin{lemma} \libel{PRODU} Etant donn\'es un faisceau prismal ${\cal F}$
sur $\Rho$ et $\rho$ un prisme de $\Rho$, au dessus de l'int\'erieur 
$\rond\rho$ de $\rho$ l'ensemble ${\cal F}(\rho)$ est une r\'eunion de produits de simplexes. 
\end{lemma}
\begin{preuve}\libel{produit} Supposons dans un premier temps que $\rho$ soit un simplexe $\tau$. 
Tout prisme de l'ensemble ${\cal F}(\tau)$ est 
le produit d'un simplexe d'image $\tau$ par un nombre, \'eventuellement nul, de simplexes. En effet, 
dans le cas contraire, il existe un 
prisme $\sg_0 \times \sg_1$ de ${\cal F}(\tau)$ d'image $\tau$ mais tel que ni $\sg_0$, ni $\sg_1$ n'ait pour image $\tau$. 
Cela implique qu'il existe deux sommets $a_0$ et $a_1$ de $\sg_0$, dont on note 
$\lambda_{a_0}$ et $ \lambda_{a_1}$ les coordonn\'ees 
barycentriques correspondantes  et deux sommets $b_0$ et $b_1$ de $\sg_1$, dont on note 
$\lambda_{b_0}$ et $\lambda_{b_1}$ les coordonn\'ees 
barycentriques correspondantes, tels  que les trois sommets 
$(a_0,b_0), (a_0,b_1)$ et $(a_1,b_1)$ du produit $\sg_0 \times \sg_1$ aient pour images des sommets distincts de $\tau$.
Alors, les coordonn\'ees barycentriques d'un point de l'image de $\sg_0 \times \sg_1$ 
sont fonction du produit $\lambda_{a_0} \lambda_{b_1}$, 
ce qui contredit l'hypoth\`ese de lin\'earit\'e.

Si $\sg$ est un simplexe de ${\cal F}(\tau)$ d'image  $\tau$, alors $\sg$ est le joint 
des simplexes $\sg_i = e_{\cal F}^{-1}(y_i)\cap \sg$ du bord de $\sg$ situ\'es au dessus des sommets $y_i$ de $\tau$. 
La fibre de $e_{\cal F}$ au dessus d'un point de l'int\'erieur de $\tau$ 
est donc hom\'eomorphe au produit des simplexes $\sg_i$, d'o\`u le 
r\'esultat dans ce cas. On en d\'eduit le r\'esultat pour tout prisme de  ${\cal F}(\tau)$.

Dans le cas d'un prisme $\rho$, produit de simplexes, le r\'esultat provient de ce qu'il est v\'erifi\'e au dessus de 
chacune des composantes du produit. Une autre mani\`ere  de le montrer est de subdiviser tout prisme $\rho$ de
la base $\Rho$ en simplexes, ceci par r\'ecurrence en se fixant un barycentre dans chaque simplexe 
composante du prisme $\rho$, puis \`a d\'ecomposer le morphisme  ${\cal F}(\rho) \to \rho$ 
au dessus de ces simplexes.
\end{preuve}

\begin{lemma} \libel{LORIENT} 
Pour toute orientation de $\rho$, il est \'equivalent
de se donner une orientation de la fibre  ${\cal F}(b(\rho))$ au dessus du
barycentre  $b(\rho)$ de $\rho$, une orientation compatible pour  toutes les fibres ${\cal
F}(y)_{y\in \rond \rho}$, ou une orientation de l'ensemble prismal ${\cal
F}(\rho)$.
\end{lemma}
\begin{preuve}
Cela provient de ce que la donn\'ee d'une orientation d'un espace
fibr\'e (orientable) \'equivaut \`a la donn\'ee d'une orientation de la base suivie d'une
orientation de la fibre. 
\end{preuve}
Le bord (orient\'e) d'un produit orient\'e $\sg_0
\times \sg_1$ est 
$$\partial  (\sg_0 \times \sg_1) = \partial \sg_0 \times \sg_1 +
(-1)^{\vert \sg_0 \vert} \sg_0 \times \partial \sg_1,$$ 
d'o\`u, par r\'ecurrence, le
bord orient\'e d'un prisme $\sg_0 \times \cdots \times \sg_s$ est~:  
\begin{equation}\label{bordpri}
\partial
(\sg_0 \times \cdots \times \sg_s) = \sum_{j=0}^s (-1)^{\vert \sg_0 \vert + \cdots
+\vert \sg_{j-1} \vert } \sg_0 \times \cdots \times \partial \sg_j \times \cdots
\times \sg_s.
\end{equation} 

\begin{definitions}\libel{ORIENT} a) Etant donn\'e un simplexe orient\'e $\sg$ et une face 
$\sg'$ de codimension 1,
alors $\sg'$ h\'erite de deux orientations. La premi\`ere est celle induite par 
l'orientation de $\sg$, c'est-\`a-dire la restriction de l'ordre correspondant des 
sommets de $\sg$, la seconde est celle qu'il a en tant que face de $\sg$. Ces deux orientations
diff\`erent d'un signe  appel\'e nombre d'incidence et not\'e $[\sg;\sg']$. \hfill\break 
\noindent b) Etant donn\'e un ensemble
prismal orient\'e $\Rho$ qui est une vari\'et\'e topo\-lo\-gi\-que, pour tout couple
$(\rho ,\rho')$ de prismes orient\'es tel que $\rho'$ soit une face de codimension
$1$ de $\rho$, le nombre
d'incidence $[\rho;\rho']$ est \'egal \`a $+1$ si l'orientation de $\rho'$ 
co\"\i ncide avec l'orientation du bord de $\rho$, et \`a $-1$ sinon. 
\end{definitions}

\begin{lemma}\libel{INKCID} 
\begin{description}
\item{a)} Soit $\sg$ un simplexe orient\'e, notons $\sg'_i$ ses faces de codimension 1, alors 
$\partial \sg  = \sum_i [\sg ; \sg'_i] \sg'_i$,
\item{b)} Soit $\pi=\sigma_0\times  \cdots \times \sigma_j\times \cdots \times
\sigma_s$ un prisme orient\'e, et 
$\pi'=\sigma_0\times \cdots \times \sigma'_j\times \cdots \times \sigma_s$ une face de codimension 1. 
Le nombre d'incidence $[\pi;\pi']$ est \'egal \`a  
$(-1)^{\vert \sg_0\vert +\cdots + \vert \sg_{j-1}\vert }[\sg_j ; \sg'_j]$.
\end{description}
\end{lemma}

\begin{preuve} Le a) d\'ecoule de la d\'efinition du bord d'un simplexe. 
On en d\'eduit le b) en consid\'erant les permutations. 
\end{preuve}

\begin{definition} Un faisceau prismal ${\cal F}$ sur un ensemble prismal
orient\'e $\Rho$ qui est une vari\'et\'e topologique est orient\'e si pour chaque
prisme $\rho$ de dimension maxima de $\Rho$, on a une orientation de la fibre ${\cal
F}(b(\rho))$ au dessus du barycentre de $\rho$, de telle fa\c con que si l'on munit
${\cal F}(\rho)$ de l'orientation correspondante ({\it cf} Lemme \ref{LORIENT}), les
morphismes $h_{\rho', \rho}$ sont $[\rho;\rho']$-orient\'es.
\end{definition}

Les exemples de faisceaux prismaux qui suivent sont fondamentaux pour ce travail.

\begin{example}\libel{EXEMSF} Soit $f \colon S
\rightarrow T$ un morphisme simplicial surjectif d'ensembles simpliciaux, on d\'efinit 
un faisceau prismal ${\cal S}_{f}$ sur $T$ en posant, pour tout simplexe
$\tau$ de $T$, ${\cal S}_{f}(\tau) = f^{-1}(\tau)$. On continue de noter par $f$ la projection de 
${\cal S}_{f}$ sur $T$. Si $\tau'$ est une face de
$\tau$ et $\sg$ un simplexe de $f^{-1}(\tau)$, le simplexe $\sg' = \sg \cap f^{-1}(\tau')$ est
une face de $\sg$. On peut \'ecrire $\sg$ comme le joint de $\sg'$ et de sa face
oppos\'ee $\sg''$. On d\'efinit alors $h_{\tau',\tau} \colon {\cal S}_{f}(\tau)
\rightarrow {\cal S}_{f}(\tau')$ en prenant pour $h_{\tau',\tau}\colon \sg  \to
\sg'$ la projection simpliciale de $\sg$ sur $\sg'$ selon les fibres du joint.
\end{example}

\begin{example}\libel{EXEMPF} Soit $\Delta$ un complexe simplicial
fini d'une subdivision simpliciale lin\'eaire de $\bf R \rm ^n \times \bf R \rm ^m $ 
tel que la restriction $f$ \`a $\Delta$ de la premi\`ere projection $p$ soit
une application simpliciale sur un complexe simplicial $f(\Delta)$ d'une
subdivision simpliciale lin\'eaire de $\bf R \rm ^n$. 

Supposons que $\sigma^{12}$ soit une face commune des simplexes  $\sigma^1$ et $\sigma^2$ 
de $\Delta$ et que ces trois
simplexes aient la m\^eme image $\tau$, simplexe de sommets 
$y_0, \ldots, y_s$. Alors le prisme  $\pi (\sigma^{12}) = \pi_f (\sigma^{12})$ est
un sous-prisme de $\pi(\sigma^1)$ et de $\pi (\sigma^2)$.  En fait, 
avec la notation de l'exemple \ref{JOINT}, on voit que
l'intersection de $\pi (\sigma^1)$ et $\pi (\sigma^2)$ est le prisme $\tau \times
\prod (\sigma_j^1 \cap \sigma_j^2)$, o\`u $j=0,\ldots ,s$.  

Remarquons que la composition $\theta^\sg \circ  \psi^\sg$ (Proposition \ref{JOINT}) 
donne un plongement naturel de $\pi(\sg)$ dans ${\displaystyle{{\bf R}^n \times \Pi_0^s {\bf R}^m}}$. 
On en d\'eduit un plongement de $\pi(\sigma^1)$ et de $\pi (\sigma^2)$
dans ${\bf R}^n \times \Pi_0^s {\bf R}^m$, on v\'erifie aussit\^ot que ces deux plongements 
co\"\i ncident sur $\sg^{12}$ et que l'on a donc d\'efini un plongement de 
$\pi(\sigma^1) \cup \pi (\sigma^2)$ dans ${\bf R}^n \times \Pi_0^s {\bf R}^m$. 
Cela montre
que, si $\tau\in f(\Delta)$ est fix\'e,  la r\'eunion des prismes $\pi (\sigma)$
tels que $f(\sg )=\tau$, plong\'ee de la fa\c con naturelle que l'on vient de
d\'ecrire dans ${\bf R}^n \times \Pi_0^s {\bf R}^m$, est un
sous-ensemble prismal  $\cal F(\tau)$ muni d'un morphisme prismal
surjectif $e_\tau \colon \cal F(\tau) \rightarrow \tau$. L'ensemble prismal 
$e_\tau^{-1}(\rond \tau)$ est naturellement isomorphe \`a 
$\rond \tau \times \cup_k (\Pi_{j=0}^s (\sg^k_j))$ o\`u les $\sg^k$ sont les simplexes 
de $e_\tau^{-1}(\rond \tau)$ d'image $\tau$ et $\sg^k_j = \sg^k \cap f^{-1}(y_j)$. 
On peut donc appeler fibre type $F_\tau$ de $e_\tau$ la r\'eunion 
des $\Pi_{j=0}^s (\sg^k_j)$. 

Ce morphisme prismal a la propri\'et\'e que l'image inverse d'un simplexe {\it
ferm\'e} de $f(\Delta)$ est le produit de ce simplexe par un ensemble
prismal ; nous pourrons donc y d\'efinir des formes de Whitney relatives. 

Supposons que $\tau '$ soit une face de $\tau$ et posons  $\sigma'=\sigma \cap f^{-1}(\tau')$.
D'apr\`es ce qui pr\'ec\`ede, chaque prisme de $\cal F(\tau)$ est de la forme $\tau
\times \sigma_0\times \cdots \times \sigma_s$ o\`u $\sigma_j = f^{-1}(y_j)\cap
\sigma$. L'homomorphisme $h_{\tau',\tau}$ de $\cal F (\tau)$ dans $\cal F(\tau')$ est l'homomorphisme
de d\'eg\'enerescence qui associe au prisme $\tau \times \sigma_0\times \cdots
\times \sigma_s$ le prisme $\tau' \times \sigma'_0\times \cdots \times \sigma'_s$
o\`u $\sigma'_j = f^{-1}(y_j)\cap \sigma'$ si $y_j$ est un sommet de $\tau'$ et un
point sinon. 

Notons $e_{{\cal P}_f} \colon {\cal P}_f\to f(\Delta)$ le faisceau prismal
ainsi obtenu. 
\end{example}

Le th\'eor\`eme qui suit  montre  l'existence d'un morphisme canonique $\psi: {\cal
P}_f\rightarrow S_f$. 


\begin{figure}
\begin{tikzpicture} \label{figure1}
\draw (-5,-2) -- (-5,2) -- (-7,1) -- (-7,-1) -- (-5,-2) -- (-3,0) -- (-5,2);
\draw (-3,0) -- (-5,0) -- (-7,1);
\draw (-5,0) -- (-7,-1);

\draw (0,-3) -- (0,3) -- (2,3) -- (2,-3) -- (0,-3);
\draw (0,1) -- (2,1);
\draw (0,-1) -- (2,-1);

\draw (3.5,-2) -- (3.5,2);
\filldraw [black] (3.5,-2) circle (1.5pt);
\filldraw [black] (3.5,0) circle (1.5pt);
\filldraw [black] (3.5,2) circle (1.5pt);

\draw (5,-2) -- (5,2) -- (7,2) -- (7,-2) -- (5,-2);
\draw (5,0) -- (7,0);

\draw (-7,-4.5) -- (-3,-4.5);
\filldraw [black] (-7,-4.5) circle (1.5pt);
\filldraw [black] (-5,-4.5) circle (1.5pt);
\filldraw [black] (-3,-4.5) circle (1.5pt);
\draw (-7,-4.8) node{$y_0$};
\draw (-6,-4.2) node{$\tau_1$};
\draw (-5,-4.1) node{$\tau'$};
\draw (-5,-4.8) node{$y_1$};
\draw (-4,-4.2) node{$\tau_2$};
\draw (-3,-4.8) node{$y_2$};

\draw (0,-4.5) -- (2,-4.5);
\filldraw [black] (0,-4.5) circle (1.5pt);
\filldraw [black] (2,-4.5) circle (1.5pt);
\draw (1,-4.8) node{$\tau_1$};

\filldraw [black] (3.5,-4.5) circle (1.5pt);
\draw (3.5,-4.8) node{$\tau'$};

\draw (5,-4.5) -- (7,-4.5);
\filldraw [black] (5,-4.5) circle (1.5pt);
\filldraw [black] (7,-4.5) circle (1.5pt);
\draw (6,-4.8) node{$\tau_2$};

\draw (-5.7,1) node{$\sigma_1$};
\draw (-6.1,0) node{$\sigma_2$};
\draw (-5.7,-1) node{$\sigma_3$};
\draw (-4.3,0.6) node{$\sigma_4$};
\draw (-4.3,-0.6) node{$\sigma_5$};

\draw (1,2) node{$\pi(\sigma_1)$};
\draw (1,0) node{$\pi(\sigma_2)$};
\draw (1,-2) node{$\pi(\sigma_3)$};

\draw (6,1) node{$\pi(\sigma_4)$};
\draw (6,-1) node{$\pi(\sigma_5)$};

\draw[->] (-4,-2) -- (-4,-3.5);
\draw (-4.3,-2.7) node{$e_{{\cal S}_f}$};

\draw[->] (1,-3.4) -- (1,-4.1);
\draw (0.7,-3.7) node{$e_{\tau_1}$};

\draw[->] (3.5,-3.2) -- (3.5,-3.9);
\draw (3.2,-3.5) node{$e_{\tau'}$};

\draw[->] (6,-3.2) -- (6,-3.9);
\draw (5.7,-3.5) node{$e_{\tau_2}$};

\draw[->] (-1,0.8) -- (-2.5,0.8);
\draw (-1.8,0.5) node{$\psi$};

\end{tikzpicture}

\medskip
${\cal S}_f$ au dessus de $\tau_1 \cup \tau_2$ \hskip 3 truecm  ${\cal P}_f$ au dessus de $\tau_1$, de $\tau'$ et de $\tau_2$
\medskip

\caption{Exemples de ${\cal S}_f$ et ${\cal P}_f$}
\end{figure}
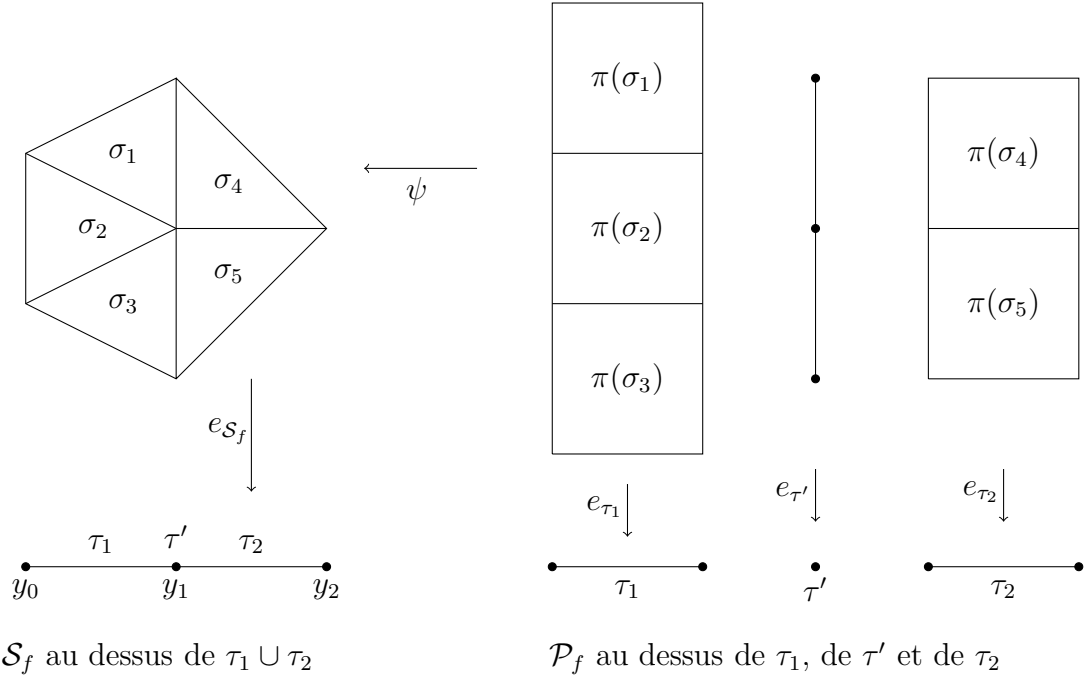



\begin{figure}
\begin{tikzpicture}
\draw (-6,10.5) -- (-2,10) -- (0,11) -- (-1,13) -- (-6,10.5);
\draw (-2,10) -- (-1,13);
\draw [dashed] (-6,10.5) --  (0,11);
\draw (-3.66,10.2) -- (-3.1,11.95);
\draw [dashed] (-3.66,10.2) -- (-2.68,10.75) -- (-3.1,11.95);

\draw (3.5,10) -- (5.5,11) --(4.5,13) -- (3.5,10) ;


\draw (-6,6) -- (-1,6) -- (0,8) -- (-6,9) -- (-6,6);
\draw (-6,9) -- (-1,6);
\draw [dashed] (-6,6) --  (0,8);
\draw (-3.5,6) -- (-3.5,7.5) -- (-3,8.5);
\draw [dashed] (-3.5,6) -- (-3,7) -- (-3,8.5);

\draw (4,6) -- (5,8) ;


\draw (-6,2) -- (0,3) -- (-5,5) -- (-6,2);
\draw [dashed] (-6,2) --  (-4,3) -- (0,3);
\draw [dashed] (-4,3) --  (-5,5);
\draw (-4,2.33) -- (-3.33,4.33);
\draw [dashed] (-4,2.33) -- (-2.66,3) -- (-3.33,4.33);

\filldraw [black] (4.5,3) circle (1.5pt);


\draw (-6,0) -- (0,0);
\filldraw [black] (-6,0) circle (1.5pt);
\filldraw [black] (0,0) circle (1.5pt);
\filldraw [black] (4.5,0) circle (1.5pt);

\draw (-4.5,0.3) node{$\tau$};
\draw (-3.6,-0.4) node{$y$};
\draw (-3.6,0) node{$\times$};
\draw (0,0.3) node{$\tau'$};
\draw (4.5,0.3) node{$\tau'$};

\draw[->] (-2.3,1.8) -- (-2.3,0.8);
\draw (-2.7,1.5) node{$e_{{\cal S}_f}$};

\draw[->] (4.5,2) -- (4.5,1);
\draw (4.2,1.5) node{$e_{{\cal S}_f}$};

\end{tikzpicture}
\bigskip
\caption{${\cal S}_f$ au dessus de $\tau$ et de $\tau'$} 
\end{figure}
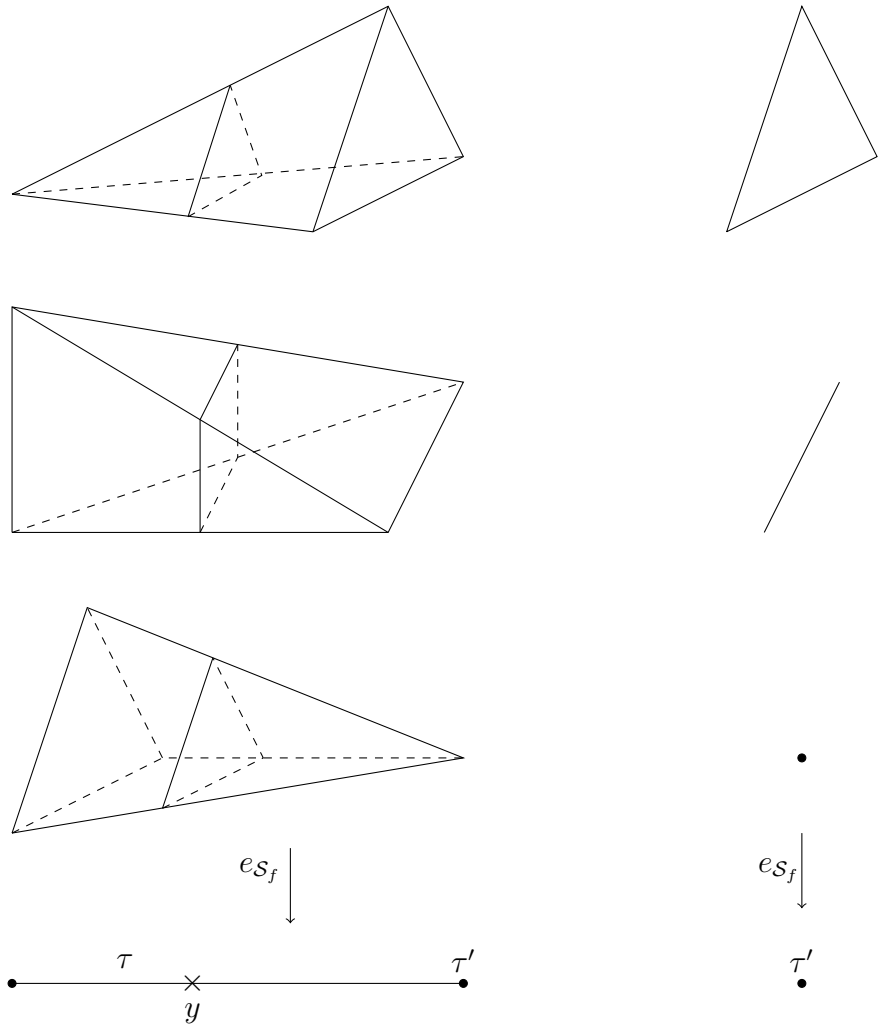


\begin{figure}
\begin{tikzpicture}
\draw (-8,1) -- (-4,0) -- (-2,1) -- (-3,3) -- (-8,3) -- (-8,1);
\draw (-8,3) -- (-4,0) -- (-3,3);
\draw (-6,0.5) -- (-6,1.5) -- (-5.5,3);
\draw [dashed] (-8,1) -- (-2,1) -- (-8,3);
\draw [dashed] (-6,0.5) -- (-5,1) -- (-5,2) -- (-5.5,3);
\draw [dashed] (-5,2) -- (-6,1.5);

\draw (0,0) -- (5,0) -- (7,1) -- (7,3) -- (5,2) -- (0,2) -- (1,4.5) -- (6,4.5) -- (7,3);
\draw (5,0) -- (5,2) -- (6,4.5);
\draw (0,0) -- (0,2);
\draw [dashed] (0,0) -- (2,1) -- (7,1);
\draw [dashed] (2,1) -- (2,3) -- (7,3);
\draw [dashed] (0,2) -- (2,3) -- (1,4.5);
\draw (-8,-2) -- (-3,-2);
\draw (1,-2) -- (6,-2);

\draw (-5.5,-2.3) node{$\tau$};
\draw (3.5,-2.3) node{$\tau$};

\draw[->] (-5.5,-0.5) -- (-5.5,-1.5);
\draw (-5.9,-1) node{$e_{{\cal S}_f}$};

\draw[->] (3.5,-0.5) -- (3.5,-1.5);
\draw (3.1,-1) node{$e_{\tau}$};

\draw[->] (-0.5,1.8) -- (-1.7,1.8);
\draw (-1.1,1.5) node{$\psi$};

\draw[->] (-0.5,-1.8) -- (-1.7,-1.8);
\draw (-1.1,-1.5) node{id};

\end{tikzpicture}

\medskip
${\cal S}_f$ au dessus de $\tau$ \hskip 6truecm ${\cal P}_f$ au dessus de $\tau$
\medskip

\caption{${\cal S}_f$ et  ${\cal P}_f$ au dessus de $\tau$}
\end{figure}



\begin{figure}
\begin{tikzpicture}
\draw (-9,1) -- (-6,-0.5) -- (-4,1.5) -- (-4,3.5) -- (-9,3) -- (-9,1);
\draw (-9,3) -- (-6,-0.5) -- (-4,3.5);
\draw [dashed] (-9,1) -- (-4,1.5) -- (-9,3);

\draw (0,4) -- (0,0) -- (3,-1.5) -- (5,0.5) -- (5,4.5) -- (0,4) -- (3,2.5) -- (3,-1.5);
\draw (0,2) -- (3,0.5) -- (5,2.5);
\draw (3,2.5) -- (5,4.5);
\draw [dashed] (0,2) -- (5,2.5);
\draw [dashed] (0,0) -- (5,0.5);

\draw (0,-3) -- (3,-4.5) -- (5,-2.5) --  (0,-3);
\draw (-9,-3) -- (-6,-4.5) -- (-4,-2.5) --  (-9,-3);

\draw (-6.5,-3.3) node{$\tau$};
\draw (2.7,-3.3) node{$\tau$};

\draw[->] (-6.5,-1) -- (-6.5,-2);
\draw (-6.9,-1.5) node{$e_{{\cal S}_f}$};

\draw[->] (2.5,-1.5) -- (2.5,-2.5);
\draw (2.2,-2) node{$e_{\tau}$};

\draw[->] (-1.5,1.8) -- (-2.7,1.8);
\draw (-2.1,1.5) node{$\psi$};

\draw[->] (-1.5,-3.2) -- (-2.7,-3.2);
\draw (-2.1,-3.5) node{id};

\end{tikzpicture}

\medskip
${\cal S}_f$ au dessus de $\tau$ \hskip 6truecm ${\cal P}_f$ au dessus de $\tau$
\medskip

\caption{${\cal S}_f$ et  ${\cal P}_f$ au dessus de $\tau$}
\end{figure}



\begin{figure}
\begin{center}
\begin{tikzpicture}

\draw (0.6,5) -- (1.1,5.5) -- (1.1,6.5) -- (0.6,6) -- (-0.4,6.5) -- (-0.4,5.5) -- (0.6,5);
\draw (-0.4,6.5) -- (1.1,6.5);
\draw(0.6,5) -- (0.6,6);
\draw [dashed]  (-0.4,5.5) --  (1.1,5.5);
\draw (2.4,8) -- (3.2,7.5) -- (3.6,8) -- (2.4,8);
\draw (2.2,6) -- (2.8,5.7) -- (3.1,6) -- (3.1,6.6) -- (2.8,6.3) -- (2.2,6.6) -- (2.2,6);
\draw (2.2,6.6) --  (3.1,6.6);
\draw (2.8,5.7) -- (2.8,6.3);
\draw [dashed]  (2.2,6) --  (3.1,6);
\draw (2.2,3.9) --  (2.2,4.7);
\draw (-0.7,4.5) --  (-0.7,5.1);
\filldraw [black] (-3.5,5.5 ) circle (1.5pt);
\draw (-0.4,7.2) -- (-1.3,7.2) -- (-0.7,6.9) -- (-0.4,7.2);

\draw[->] (-1.5,6) -- (-2.5,5.8);
\draw[->] (-0.2,5.3) -- (-0.5,5);
\draw[->] (1.2,5.3) -- (1.8,4.6);
\draw[->] (1.4,6.1) -- (1.8,6.1);
\draw[->] (1.4,6.9) -- (1.9,7.2);
\draw[->] (-0.2,6.7) -- (-0.5,6.9);

\draw (0.1,0.1) -- (1.1,0.6) -- (0.1,1.6) -- (-0.4,0.6) -- (0.1,0.1);
\draw (0.1,0.1) -- (0.1,1.6);
\draw [dashed] (-0.4,0.6) --  (1.1,0.6);
\draw (2.4,0.2) -- (3,0.5) -- (2.4,1.2) -- (2,0.5) -- (2.4,0.2);
\draw (2.4,0.2) -- (2.4,1.2);
\draw [dashed] (2,0.5) --  (3,0.5);
\draw (2.8,1.8) -- (3.6,2.2) -- (2.8,3) -- (2.4,2.2) -- (2.8,1.8);
\draw (2.8,1.8) -- (2.8,3);
\draw [dashed] (2.4,2.2) --  (3.6,2.2);
\draw (-1,1.4) -- (-0.4,1.7) -- (-1,2.4) -- (-1.4,1.7) -- (-1,1.4);
\draw (-1,1.4) -- (-1,2.4);
\draw [dashed] (-0.4,1.7) --  (-1.4,1.7);

\filldraw [black] (2,-1.3) circle (1.5pt);
\filldraw [black] (-3.5,0 ) circle (1.5pt);
\filldraw [black] (-0.7,-0.75) circle (1.5pt);

\draw[->] (-1.5,0.5) -- (-2.5,0.3);
\draw[->] (-0.2,-0.2) -- (-0.5,-0.5);
\draw[->] (0.7,-0.2) -- (1.3,-0.7);
\draw[->] (1.4,0.6) -- (1.8,0.6);
\draw[->] (1.4,1.4) -- (1.9,1.7);
\draw[->] (-0.4,1.1) -- (-0.7,1.4);

\draw (-0.2,-4.5) -- (0.8,-4.5) -- (1.3,-4) -- (1.3,-3) -- (0.8,-3.5) -- (-0.2,-3.5) -- (-0.2,-4.5);
\draw (0.8,-4.5) -- (0.8,-3.5);
\draw (-0.2,-3.5) -- (0.3,-3) -- (1.3,-3);
\draw [dashed] (-0.2,-4.5) -- (0.3,-4) -- (1.3,-4);
\draw [dashed] (0.3,-4) -- (0.3,-3);
\draw (2.2,-4.2) -- (2.9,-4.2) -- (2.9,-3.5) -- (2.2,-3.5) -- (2.2,-4.2);
\draw (2.4,-2.1) -- (3.4,-2.1);
\draw (2.1,-5.5) -- (2.1,-6.5);
\draw (-3.7,-4.7) -- (-3.2,-4.2);
\draw (-1,-5.9) -- (-0.7,-5.5) -- (-0.7,-4.8) -- (-1,-5.2) --  (-1,-5.9);
\draw (-1.7,-3.2) -- (-1,-3.2) -- (-0.7,-2.8) -- (-1.4,-2.8) --   (-1.7,-3.2);

\draw[->] (-1.5,-4) -- (-2.5,-4.2);
\draw[->] (-0.2,-4.7) -- (-0.5,-5);
\draw[->] (0.7,-4.7) -- (1.3,-5.4);
\draw[->] (1.5,-3.8) -- (1.9,-3.8);
\draw[->] (1.5,-3.1) -- (2.2,-2.5);
\draw[->] (-0.4,-3.4) -- (-0.7,-3.1);

\draw (3,-7.5) -- (-3.5,-10) -- (2,-11) --  (3,-7.5);
\draw (0.5,-9.5) node{$\tau$};

\end{tikzpicture}

\caption{Sp\'ecialisations possibles d'un
5-simplexe au dessus d'un 2-simplexe $\tau$}
On a dessin\'e les fibres
au dessus des points g\'en\'eriques des diff\'erentes faces de $\tau$.
\end{center}
\end{figure}


\begin{theorem}\libel{COMPFA} Soit $f\colon \Delta\to T$ un
morphisme simplicial. Consid\'erons le faisceau prismal ${\cal S}_f$ de l'exemple
\ref{EXEMSF} et le faisceau prismal ${\cal P}_f $ de l'exemple \ref{EXEMPF}.
Il existe un
morphisme surjectif $\psi \colon {\cal P}_f  \rightarrow {\cal S}_f$ de faisceaux
prismaux sur $f(\Delta)$ tel que \par 
\noindent a) Pour chaque simplexe $\tau$ de $f(\Delta)$ les morphismes $\theta^{\sg}$ et  
$\psi^{\sg}$ de la Proposition \ref{JOINT} d\'efinissent un isomorphisme $\theta_\tau$ de 
$f^{-1}({\rond \tau})$ sur ${\rond \tau}\times F_{\tau}$. Le compos\'e $\theta_\tau \circ \psi$ 
s'\'etend en un isomorphisme prismal de ${\cal P}_f (\tau)$ sur $\tau \times F_\tau$. 

\noindent b) La formation du faisceau ${\cal P}_f$ est fonctorielle et universelle. Plus
pr\'ecis\'ement~:
\begin{enumerate}
\item pour tout morphisme prismal $ \phi^{f} \colon \Rho' \rightarrow
\Delta$, il existe un faisceau prismal ${\cal F}(\Rho')$ sur $f(\Delta)$ et un
morphisme de faisceaux prismaux ${\cal F}(\phi^{f}) \colon {\cal F}(\Rho') \rightarrow {\cal P}_f$,
\item \'etant donn\'e un morphisme surjectif $\psi'
\colon {\cal F}' \rightarrow {\cal S}_f$ de faisceaux prismaux sur $f(\Delta)$ tel
que l'image inverse de tout simplexe $\tau$ de $f(\Delta)$ par $f\circ \psi'$ soit
r\'eunion de produits de $\tau$ par des prismes, il existe un unique morphisme
$\chi \colon {\cal F}' \rightarrow {\cal P}_f$ rendant commutatif le diagramme
de faisceaux prismaux au dessus de $f(\Delta)$ : 
\end{enumerate}
$$
\xymatrix{ 
{\cal F}' \ar[rr] ^{\chi} \ar[rd]_{\psi'} && {\cal P}_f \ar[ld]^\psi \\
&{\cal S}_f \; .
 }
$$
\end{theorem}

\begin{preuve}
 Montrons l'existence d'un morphisme de faisceaux prismaux $\psi \colon {\cal
P}_f \rightarrow {\cal S}_f$. Un point d'un prisme est d\'etermin\'e
par ses coordonn\'ees barycentriques dans chaque simplexe. Pour tout $\sg$ on a
d\'efini dans l'exemple \ref{JOINT} un morphisme prismal $\psi^\sg\colon \pi (\sg )\to
\sg$. Les morphismes $\psi ^\sigma$ et $ \psi ^{\sigma'}$ co\"\i ncident sur $\pi
(\sigma \cap \sigma')$, et nous avons donc d\'efini un morphisme prismal $\psi_\tau
\colon {\cal P}_f(\tau) \rightarrow f^{-1}(\tau) ={\cal S}_f(\tau)$. 

Entre autres, si $\tau'$ est une face d'un
simplexe $\tau$, on a un diagramme commutatif : 
$$
\xymatrix@C=1.3cm { 
{\cal P}_f(\tau)  \ar[r]^-{\psi_\tau}  \ar[d]^{h_{\tau',\tau}} & {\cal S}_f(\tau) \ar[d]^{h_{\tau',\tau}} \\ 
{\cal P}_f(\tau')   \ar[r]^-{\psi_{\tau'}}  & {\cal S}_f(\tau').}
$$

Le a) r\'esulte alors de la Proposition \ref{JOINT} et de 
 l'exemple \ref{EXEMPF}. 

Prouvons b); il suffit de v\'erifier l'\'enonc\'e restreint \`a un simplexe
$\sigma$ de ${\cal S}_f(\tau)$. La d\'efinition de ${\cal P}_f$ 
et la structure de produit de $\pi(\sg)$ impliquent que pour
tout prisme $\pi'$ de l'image inverse de $\sigma$ par $\psi'$, on a une application
naturelle de $\pi'$ dans $\pi (\sigma)$, d'o\`u le r\'esultat.
\end{preuve}

\begin{remarques}\libel{REM} \par 1) Le morphisme $\psi$ est
essentiellement un \'eclatement comme le montre l'\'ecriture locale $\lambda_i=t_j\mu_{j,i}$. 
En particulier, le morphisme $\psi^\sg : \pi(\sg) \to \sg$ n'est un isomorphisme que si 
$\sg$ est isomorphe \`a son image par $f$. 
Soulignons que le morphisme $\pi(\sg) = \tau\times\sg_0\times\cdots\times\sg_s \to \sg_j$
qui rend commutatif le diagramme
$$
\xymatrix@C=0.4cm{ 
\sg_0 * \cdots * \sg_s  \ar[rd]_{\theta^{\sigma_j}} && \tau\times\sg_0\times\cdots\times\sg_s 
\ar[ll]_{\psi^\sg}    \ar[ld]^{{\rm pr}_j} \\
&\sg_j }
$$
o\`u $\theta^{\sigma_j}$ d\'esigne le morphisme introduit en (\ref{detheta}), est
bien la
$j$-\`eme projection ${\rm pr}_j$ du produit.\par\noindent

2) Le jacobien de l'application $\psi^\sg$ est \'egal \`a :
$$t_0^{\vert \sg_0 \vert} t_1^{\vert \sg_1 \vert} \cdots t_s^{\vert \sg_s \vert}.$$

3) La cat\'egorie des ensembles et
morphismes prismaux est ``la plus petite" cat\'egorie contenant celle des ensembles et
morphismes simpliciaux et dans laquelle on a existence et unicit\'e \`a isomorphisme pr\`es du
produit fibr\'e.\par
 L'existence du produit fibr\'e d\'ecoule du fait que, puisque les applications
prismales sont lin\'eaires sur chaque simplexe, le sous-ensemble d'un produit
$\pi_1\times \pi_2$ de prismes d\'efini par la condition $f_1(x_1)=f_2(x_2)$, o\`u
$f_1\colon \pi_1\to \tau$ et $f_2\colon \pi_2\to
\tau$ sont des morphismes prismaux, est un prisme. La v\'erification de la
propri\'et\'e universelle est imm\'ediate.
\end{remarques}

\noindent{\sl Caract\'erisation des faisceaux prismaux provenant de
morphismes simpliciaux}

\begin{definitions}\libel{DIMREL} a) Soit ${\cal F}$ un
faisceau  prismal sur un complexe simplicial $T$, on dit que le prisme $\pi \in {\cal
F}(\tau)$ est trivial s'il s'\'ecrit $\pi=\tau \times\sigma_0 \times \cdots \times 
\sigma_s$;

b) On appelle {\it dimension relative}  d'un prisme $\pi\in {\cal F}(\tau)$ et on
note $\hbox{dim}_{rel}(\pi)$,  la diff\'erence $\dim \pi - \dim \tau$.

c) On dit que le prisme $\pi \in {\cal F}(\tau)$ est {\it \'equidimensionnel} au
dessus d'une face $\tau'$ de $\tau$ si la dimension relative de $\pi\vert \tau'$
est \'egale \`a celle de $\pi$. 
\end{definitions}

Un prisme de dimension relative nulle est \'equidimensionnel au dessus de toutes les
faces de $\tau$.   Le morphisme $\sg_1
\to \tau_1$ de la figure 1 n'est pas \'equidimensionnel au dessus du  sommet 
$y_0$ de $\tau_1$. 

\begin{proposition}
[Caract\'erisation des faisceaux de la forme ${\cal S}_f$] \libel{PROCAR}  Un faisceau prismal
${\cal F}$  sur un complexe simplicial $T$ est de la forme  ${\cal S}_f$, 
pour un morphisme simplicial $f:X \to T$, si et
seulement si~:
\begin{enumerate}
\item Tous les prismes de ${\cal F}$ sont des simplexes,
\item Pour tout couple $\tau'<\tau$ de simplexes de $T$, et tout simplexe $\sigma$ de
${\cal F}(\tau)$, le morphisme $h_{\tau',\tau} : {\cal F}(\tau) \to {\cal F}(\tau')$ est surjectif et  on a un
isomorphisme simplicial $h_{\tau',\tau}(\sigma)\cong \sigma\vert_{\tau'}$.
\end{enumerate}
\end{proposition}

\begin{preuve} Le fait que les conditions 1) et 2) soient n\'ecessaires
r\'esulte aussit\^ot de la construction de l'exemple \ref{EXEMSF}. Montrons qu'elles
sont suffisantes; nous construisons l'espace $X$ de la d\'efinition \ref{ANALYT}  par
recollement. Notons $X_{\tau}= e_\tau^{-1}(\tau)$, alors pour toute face $\tau'$ de
$\tau$, on a $X_{\tau}\vert_ {\tau'}=X_{\tau'}=h_{\tau',\tau}(X_{\tau})$.
D\'efinissons l'espace $X$ comme quotient  de la r\'eunion des $X_\tau$ par la
relation d'identification des restrictions au dessus des faces des simplexes de
$T$. L'application $f\colon X \to T$  est naturellement d\'efinie. 
\end{preuve}

\begin{lemma}\libel{LEMEQU} Consid\'erons le faisceau
prismal  ${\cal S}_f$ de base $T$. Pour tout couple $\tau'<\tau$ de simplexes de $T$,  et
tout simplexe
$\sigma$ de ${\cal S}_f(\tau)$, on a \'equivalence des propri\'et\'es  suivantes~:
\begin{description}
\item[{\sl (i)}] $\hbox{dim}_{rel}h_{\tau',\tau}(\sigma)=\hbox{dim}_{rel}(\sigma)$ 
\item[{\sl (ii)}] $\sigma$ est \'equidimensionnel au dessus de la face $\tau'$.
\item[{\sl (iii)}]  Si $\tau''$
est la face oppos\'ee de $\tau'$ dans $\tau$, la projection  $\sg\vert_{\tau''}\to
\tau''$ est un isomorphisme.
\end{description}
\end{lemma}

\begin{preuve}
Consid\'erons $\tau$
comme le joint de $\tau'$ et de sa face oppos\'ee $\tau''$. Puisque  $$\dim \sigma =\dim
\sigma\vert_{\tau'} + \dim \sigma\vert_{\tau''}+1$$ et $h_{\tau',\tau}(\sigma)\cong
\sigma\vert_{\tau'}$, l'assertion (i) est \'equivalente \`a dire que $\sigma\vert_{\tau''}$ est
isomorphe \`a $\tau''$ d'o\`u le r\'esultat (voir Figures 1 et 3). 
\end{preuve}

\begin{proposition} [Caract\'erisation des faisceaux du type
${\cal P}_f$]\libel{PROPPF}   Un faisceau prismal ${\cal F}$ sur un complexe simplicial $T$ est de type
${\cal P}_f$ o\`u $f\colon X\rightarrow T$ est un morphisme simplicial si et
seulement si  les conditions suivantes sont r\'ealis\'ees~:
\begin{description}
\item[{\sl a)}] Tout prisme $\pi \in {\cal F}(\tau)$ est le produit de simplexes $\pi=
\tau\times \sigma_0 \times \cdots \times \sigma_s $ avec $s=\dim \tau$, 
\item[{\sl b)}] Si $\tau'<\tau$ le morphisme $h_{\tau',\tau}  : {\cal F}(\tau) \to {\cal F}(\tau')$ est surjectif et si un prisme $\pi
\in {\cal F}(\tau)$ s'\'ecrit $\pi= \tau\times \sigma_0 \times \cdots \times
\sigma_s$, alors on a $h_{\tau',\tau}(\pi) =\tau'\times \sigma'_0 \times \cdots
\times \sigma'_k $ o\`u chaque  $\sigma'_{\beta}$ est l'un des
$\sigma_\alpha$.
\end{description}
\end{proposition}

\begin{preuve}
Si un faisceau prismal ${\cal F}$ est de type  ${\cal
P}_f$, il v\'erifie a) et b) par construction. Montrons la r\'eciproque. Soit donc
un  faisceau prismal ${\cal F}$ satisfaisant a) et b).  On lui associe un faisceau
prismal  ${\cal S}$ de base $T$ de la fa\c con suivante : les  prismes de ${\cal
S}(\tau)$ sont  les joints it\'er\'es $\sigma_0 * \cdots *  \sigma_s$ des simplexes
apparaissant dans  les prismes $\pi= \tau\times \sigma_0  \times \cdots \times
\sigma_s $ de ${\cal F}$.  Ce sont donc des simplexes.  La propri\'et\'e b)
implique que, pour tout simplexe $\sigma$ de ${\cal S}(\tau)$ et pour toute face
$\tau'$ de $\tau$, on a  $h_{\tau',\tau}(\sigma)\cong \sigma\vert_{\tau'}$.
D'apr\`es la proposition \ref{PROCAR}, le faisceau ${\cal S}$ est de la forme ${\cal
S}_f$ et le faisceau ${\cal F}$ est  le faisceau ${\cal P}_f$ qui lui est associ\'e par
la construction de l'exemple \ref{EXEMPF}.
\end{preuve}

\begin{remarque} \label{quisuit} Pour tout prisme  $\pi$ de  ${\cal P}_f(\tau)$ 
et pour tout sommet $\{y\}$ de $\tau$, on a 
$h_{\{y\},\tau}(\pi)=\{y\}\times  \sigma_{i(y)}$ o\`u $\sigma_{i(y)}$ est le
simplexe  $\psi(\pi)\cap (e_\tau)^{-1} (\{ y\})$ de ${\cal S}_f$.
\end{remarque}

\begin{corollary}\libel{COROPF} Consid\'erons le faisceau prismal
${\cal P}_f$ de base $T$. Pour tout couple $\tau'<\tau$ de simplexes de $T$, et tout
prisme  $\pi$ de ${\cal P}_f(\tau)$, on a l'\'equivalence des propri\'et\'es suivantes~:
\begin{description}
\item[{\sl i)}] $\hbox{dim}_{rel}h_{\tau',\tau}(\pi)=\hbox{dim}_{rel}(\pi)$ 
\item[{\sl ii)}] $h_{\tau',\tau}(\pi)=\pi\vert_{\tau'}$ 
\item[{\sl iii)}] $\psi(\pi)$ est \'equidimensionnel au dessus de la face  $\tau'$. 
\end{description}
\end{corollary}

\begin{preuve} La proposition \ref{PROPPF} et la remarque \ref{quisuit} 
impliquent l'\'equivalence de (i) et (ii), la proposition \ref{PROCAR}
et le fait que $\dim \psi(\pi)=\dim \pi$ impliquent l'\'equivalence de (i) et
(iii).
\end{preuve}

\begin{corollary}\libel{COROSF} Consid\'erons le
faisceau prismal ${\cal P}_f$ de base $T$. Pour tout 
simplexe $\tau$ de $T$, et tout prisme  $\pi$ de ${\cal
P}_f(\tau)$, on a l'\'equivalence des propri\'et\'es suivantes~:
\begin{description}
\item[{\sl i)}] $\psi(\pi)$ est \'equidimensionnel en un sommet $y_0$ de
$\tau$. 
\item[{\sl ii)}] $\pi$ est isomorphe au dessus de $\tau$ \`a un produit $\tau
\times\sigma_\pi$. 
\end{description}
\end{corollary}

\begin{preuve} Le corollaire est une cons\'equence directe de la
proposition \ref{PROPPF} et de la remarque \ref{quisuit}. Le simplexe $\sigma_\pi$ du (ii)
est $\psi(\pi)\cap (e_\tau)^{-1} (\{ y_0\})$ de ${\cal S}_f$ o\`u $\{ y_0\}$ est le sommet
de $\tau$ du (i) (voir Figures 1 et 3). 
\end{preuve}

\section{Formes de Whitney}
\setcounter{equation}{0}

\subsection{Formes r\'eguli\`eres, h\"old\'eriennes et sous-analytiques}
Soient $N$ un entier et $U$ un ouvert de ${\bf R}^N$.
On dira qu'une fonction $g$ d\'efinie sur un ouvert $U$ de ${\bf R}^N$ 
est {\it h\"old\'erienne d'exposant $\alpha$} si tout point de $U$ poss\`ede un voisinage $V$  dans $U$ 
tel qu'il existe une constante positive $C_V$ telle que pour tous $(x,x')$ dans $V$, on ait l'in\'egalit\'e
$$\vert g(x) - g(x')\vert \le C_V \vert x - x' \vert ^\alpha.$$

Soit $\mu\geq 0$ un entier. Suivant Whitney \cite{Whi}, nous dirons qu'une r-forme diff\'erentielle sur $U$
est {\it $\mu$-r\'eguli\`ere} si elle est continue sur U et satisfait les conditions suivantes~: 

- Si $\mu =0$, il existe une (r+1)-forme $\xi$ continue sur $U$ telle que
l'on ait $$ \int_{\partial \sigma} \omega = \int_{\sigma} \xi$$ pour tout
(r+1)-simplexe $\sigma$ contenu dans $U$. D'apr\`es le lemme 16a de \cite[Ch.III, p.104]{Whi}, 
la forme $\xi$ est uniquement d\'etermin\'ee par cette condition ;
on la notera d$\omega$ et on l'appellera, suivant Whitney, {\it forme d\'eriv\'ee} de
$\omega$.

- Si $\mu$ est $>0$,
$\omega$ et $d\omega$ sont diff\'erentiables de classe ${\cal C}^{\mu}$. 

Une forme $\omega$ est {\it h\"old\'erienne} si
elle est 0-r\'eguli\`ere  et si les coefficients de $\omega$ et $\xi$ sont des
fonctions h\"old\'eriennes sur $U$. 
\begin{definition}
Soit $h$ une fonction sous-analytique continue sur le compact sous-ana\-ly\-ti\-que $K\subset \R^N$ d'int\'erieur non vide. Par hypoth\`ese le graphe de $h$ est un sous-ensemble sous-analytique ferm\'e $\Gamma$ de $K\times {\mathbf P}^1$ contenu dans $K\times \R$. D'apr\`es [D], [D-W], le transform\'e de Semple-Nash $S\Gamma$ de $\Gamma$ est sous-analytique dans $\Gamma\times {\mathbf P}^N$. Il y est ferm\'e et est donc compact. Au dessus des points d'analyticit\'e de $\Gamma$, l'espace $S\Gamma$ est le lieu des points $(x,h(x), [\frac{\partial h}{\partial x_1}:\cdots :\frac{\partial h}{\partial x_N}:-1])$.  Notons $H_\infty$ l'hyperplan de ${\mathbf P}^N$ correspondant \`a la derni\`ere coordonn\'ee de $K\times \R$, et $\R^N\subset {\mathbf P}^N$ l'espace affine compl\'ementaire. Si $S\Gamma$ se trouve \^etre contenu dans $\Gamma\times {\mathbf R}^N\subset \Gamma\times {\mathbf P}^N$, nous dirons que la fonction $h$ est \textit{\`a d\'eriv\'ees born\'ees sur K}.\end{definition}
\begin{proposition}
 Si $h$ est \`a d\'eriv\'ees partielles born\'ees, les d\'eriv\'ees partielles de $h$ au sens des distibutions sont repr\'esent\'ees par des fonctions sous-analytiques born\'ees sur $K$.
 \end{proposition}
 \begin{preuve} En effet l'image de $S\Gamma\subset\Gamma\times {\mathbf R}^N$ dans $\Gamma\times \R$ par la projection sur la $i$-i\`eme coordonn\'ee est encore sous-analytique comme projection d'un ensemble sous-analytique compact et c'est le graphe d'une fonction sous-analytique sur $Z$ qui est un repr\'esentant de la d\'eriv\'ee au sens des distributions $\frac{\partial h}{\partial u_i}$, puisque les deux co\"\i ncident sur l'ouvert d'analyticit\'e de la fonction $h$ dont le compl\'ementaire est de mesure nulle. 
\end{preuve}

\begin{definition}\libel{forsousan} Une forme $\omega$ est {\it sous-analytique} si
elle est 0-r\'eguli\`ere  et si les coefficients de $\omega$ et $\xi$ sont des
fonctions sous-analytiques born\'ees sur $U$. Elle est sous-analytique continue si ses coefficients sont continus (mais pas n\'ecessairement ceux de sa diff\'erentielle, car imposer cela emp\^echerait les formes de Whitney d'\^etre continues).
\end{definition}

\begin{proposition}\libel{hreg} Soit $U$ un ouvert de ${\bf R}^N$, pour une forme
diff\'erentielle $\omega$ d\'efinie sur $U$, les conditions suivantes sont
\'equivalentes :
\begin{description}
\item[{\sl i)}] la forme $\omega$ est 0-r\'eguli\`ere dans $U$, ses coefficients sont 
sous-analytiques et born\'es et ceux de sa forme d\'eriv\'ee sont sous-analytiques.
\item[{\sl ii)}] les coefficients de la forme $\omega$ sont sous-analytiques et born\'es et sa diff\'erentielle au sens
des distributions admet un repr\'esentant sous-analytique. 
\end{description} 
\end{proposition}

\begin{preuve} (i) implique (ii) : 
La formule de Stokes pour les courants implique que la dif\-f\'eren\-tielle au sens
des distributions satisfait
$$ \int_{\partial \sigma} \omega = \int_{\sigma} d\omega$$
On en d\'eduit que $d\omega$ admet comme repr\'esentant la forme d\'eriv\'ee de 
$\omega$ qui est sous-analytique.

(ii) implique (i): Pour les m\^emes raisons, un repr\'esentant sous-analytique de la 
dif\-f\'e\-ren\-tielle au sens des distributions doit co\"\i ncider avec la forme d\'eriv\'ee de $\omg$ ce qui
montre que $\omega$ est sous-analytique.
\end{preuve}

Rappelons que d'apr\`es {\it loc. cit.}, lorsque $\mu$ est \'egal \`a z\'ero, une forme
$\omega$ est r\'eguli\`ere et d$\omega$ admet pour repr\'esentant $\xi$ si et seulement s'il existe une suite
$\omega_i$ de formes de classe ${\cal C}^1$ sur $U$ telle que, uniform\'ement sur tout
compact, $\omega_i$ tende vers $\omega$ et des repr\'esentants de 
$d\omega_i$ tendent vers $\xi$ au sens des courants. Whitney en
d\'eduit que si $f$ est un morphisme ${\cal C}^1$ d'un ouvert $U$ de ${\bf R}^N$
dans ${\bf R}^p$ et $\omega$ une forme r\'eguli\`ere sur un voisinage de l'image de
$U$, alors $f^*\omega$ est r\'eguli\`ere dans $U$. 

\subsection{Formes diff\'erentielles sur les complexes simpliciaux et les 
ensembles  prismaux}

Soit $\Delta$ un complexe simplicial lin\'eaire dans ${\bf R}^N$. Une r-forme
diff\'erentielle sur $\Delta$ est la donn\'ee pour tout simplexe $\sigma$ de
$\Delta$ d'une r-forme $\omega_{\sigma}$ d\'efinie sur $\sg$, c'est-\`a-dire que
$\omega_{\sigma}$ est une section d\'efinie sur $\sigma$ du fibr\'e $\Lambda^r
T^*{\bf R}^N$ de telle mani\`ere que, pour toute face $ i : \sg' \hookrightarrow
\sg$ de tout simplexe $\sg$ de $\Delta$, on ait $\omg_\sigma\vert_{\sg'}=\omega_{\sigma'}$.
Soulignons qu'il s'agit ici de la restriction \`a $ \sg'$ des coefficients de $\omg(\sg)$ et non 
pas de l'image r\'eciproque $i^*(\omg(\sg))$ par l'inclusion de $\sg'$ dans $\sg$.

Le faisceau des formes diff\'erentielles sur ${\bf R}^N$ est mou (voir
\cite[Chap.II, Exemple 3.7.1]{Go}). Rappelons que ceci signifie que toute forme 
diff\'erentielle d\'efinie sur un ferm\'e se prolonge \`a l'espace tout entier. 
On peut donc supposer que toute forme diff\'erentielle sur $\Delta$ est obtenue en 
restreignant \`a chaque simplexe $\sg$ de $\Delta$ une forme diff\'erentielle d\'efinie sur
un voisinage ouvert de $\sg$ de telle fa\c con que les restrictions de  deux telles
formes sur leur ouvert de d\'efinition commun co\"\i ncident.

Nous consid\'ererons ici des formes diff\'erentielles $\omg$ dont les coefficients
sont des  fonctions sous-analytiques sur $\Delta$ analytiques dans l'int\'erieur de
chaque simplexe. \par \noindent La d\'efinition des formes diff\'erentielles
s'\'etend aussit\^ot aux ensembles prismaux~:

\begin{definition}\libel{FORDIF} Soit $\Pi$ un ensemble prismal. On
appelle {\sl $r$-forme diff\'erentielle $\mu$-r\'eguli\`ere} (resp. sous-analytique) sur ${\Pi}$ la donn\'ee pour 
chaque prisme $\pi$ de $\Pi$ d'une $r$-forme
diff\'e\-ren\-tiel\-le $\mu$-r\'eguli\`ere (resp. sous-analytique) sur un voisinage ouvert de $\pi$ dans l'un de
ses plongements affines, de telle fa\c con que les formes diff\'erentielles
correspondant \`a deux prismes co\"\i ncident dans un voisinage de leur
intersection dans un plongement affine commun. 
\end{definition}
    
Clairement, la d\'efinition ne d\'epend pas des plongements choisis.

\subsection{Formes de Whitney sur les complexes simpliciaux} 

Notons $a_i$ les sommets du complexe simplicial $\Delta$; tout point $x$ de
$\Delta$  s'\'ecrit $x=\sum \lambda_i (x)a_i$ o\`u les 
coordonn\'ees barycentriques $\lambda_i$ 
satisfont  $\lambda_i(x)\geq 0$ et
$\sum_i\lambda_i(x)=1$. Le support de $\lambda_i$ est l'\'etoile ouverte de $a_i$ dans $\Delta$.
Nous allons, comme Whitney, construire une partition de
l'unit\'e subordonn\'ee au recouvrement ouvert de $\Delta$ constitu\'e des
\'etoiles des sommets de $\Delta$. 

Pour tout $i$, on note $F_i$ l'ensemble des points de $\Delta$ dont la $i$-i\`eme
coordonn\'ee barycentrique $\lambda_i(x)$ est $\geq \frac{1}{N+1}$ et $G_i$ l'ensemble
des points de $\Delta$ dont la $i$-i\`eme coordonn\'ee barycentrique $\lambda_i(x)$
est $\leq \frac{1}{N+2}$ . 

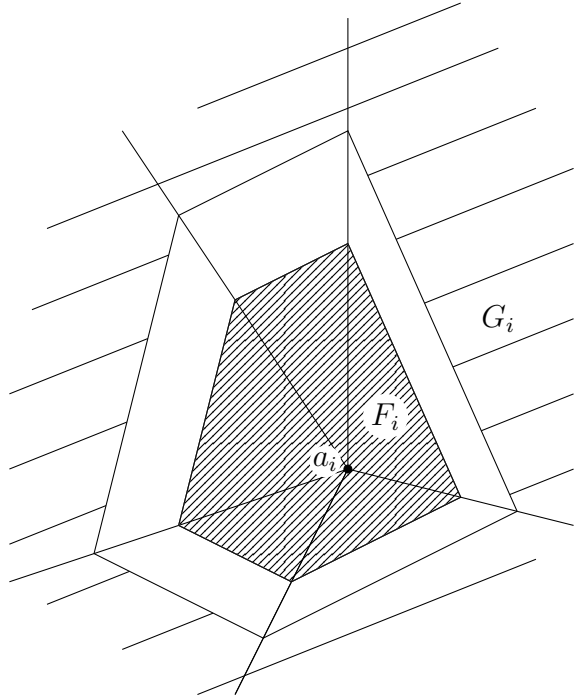
\begin{figure}[h]
\begin{center}
\begin{tikzpicture}

\node (A) at (2.25,-0.5625) [coordinate] {} ;
\node (B) at (0,4.5) [coordinate] {} ;
\node (C) at (-2.25,3.375) [coordinate] {} ;
\node (D) at (-3.375,-1.125) [coordinate] {} ;
\node (E) at (-1.125,-2.25) [coordinate] {} ;

\draw (2.25,-0.5625) -- (0,4.5) -- (-2.25,3.375) -- (-3.375,-1.125) -- (-1.125,-2.25)-- (2.25,-0.5625);
\draw (1.5,-0.375) -- (0,3) -- (-1.5,2.25) -- (-2.25,-0.75) -- (-0.75,-1.5)-- (1.5,-0.375);
\draw (0,0) -- (3,-0.75);
\draw (0,0) -- (0,6);
\draw (0,0) -- (-3,4.5);
\draw (0,0) -- (-4.5,-1.5);
\draw (0,0) -- (-1.5,-3);
\draw (0,0) -- (-1.5,-3);
\filldraw [black] (0,0) circle (1.5pt);

\filldraw[color=lightgray,pattern=north east lines] (1.5,-0.375) -- (0,3) -- (-1.5,2.25) -- (-2.25,-0.75) -- (-0.75,-1.5)-- (1.5,-0.375);
\draw (1.5,-0.375) -- (0,3) -- (-1.5,2.25) -- (-2.25,-0.75) -- (-0.75,-1.5)-- (1.5,-0.375);

\draw (2,2) node{$G_i$};
\fill[white] (0.5,0.7) circle (0.3);
\draw (0.5,0.7) node{$F_i$};
\fill[white] (-0.3,0.1) circle (0.21);
\draw (-0.3,0.1) node{$a_i$};

\node (a1) at (-4.5,0) [coordinate] {} ;
\node (b1) at (3,3) [coordinate] {} ;
\coordinate (c1) at (intersection of a1--b1 and D--C);
\coordinate (d1) at (intersection of a1--b1 and B--A);
\draw (a1)--(c1);
\draw (d1)--(b1);

\node (a2) at (-4.5,1) [coordinate] {} ;
\node (b2) at (3,4) [coordinate] {} ;
\coordinate (c2) at (intersection of a2--b2 and D--C);
\coordinate (d2) at (intersection of a2--b2 and B--A);
\draw (a2)--(c2);
\draw (d2)--(b2);

\node (a3) at (-4.5,2) [coordinate] {} ;
\node (b3) at (3,5) [coordinate] {} ;
\node (e3) at (2.5,0) [coordinate] {} ;
\node (f3) at (2.5,5) [coordinate] {} ;
\node (g3) at (-4.2,0) [coordinate] {} ;
\node (h3) at (-4.2,5) [coordinate] {} ;
\coordinate (c3) at (intersection of a3--b3 and D--C);
\coordinate (d3) at (intersection of a3--b3 and B--A);
\coordinate (k3) at (intersection of a3--b3 and e3--f3);
\coordinate (l3) at (intersection of a3--b3 and g3--h3);
\draw (l3)--(c3);
\draw (d3)--(k3);

\node (a4) at (-4.5,3) [coordinate] {} ;
\node (b4) at (3,6) [coordinate] {} ;
\node (e4) at (-4,0) [coordinate] {} ;
\node (f4) at (-4,5) [coordinate] {} ;
\node (g4) at (2,0) [coordinate] {} ;
\node (h4) at (2,5) [coordinate] {} ;
\coordinate (c4) at (intersection of e4--f4 and a4--b4);
\coordinate (i4) at (intersection of g4--h4 and a4--b4);
\draw (c4)--(i4);

\node (a5) at (-4.5,3.8) [coordinate] {} ;
\node (b5) at (3,6.8) [coordinate] {} ;
\node (e5) at (-2,0) [coordinate] {} ;
\node (f5) at (-2,5) [coordinate] {} ;
\node (g5) at (1.5,0) [coordinate] {} ;
\node (h5) at (1.5,5) [coordinate] {} ;
\coordinate (c5) at (intersection of e5--f5 and a5--b5);
\coordinate (i5) at (intersection of g5--h5 and a5--b5);
\draw (c5)--(i5);

\node (a6) at (-4.5,-1) [coordinate] {} ;
\node (b6) at (3,2) [coordinate] {} ;
\coordinate (c6) at (intersection of a6--b6 and D--C);
\coordinate (d6) at (intersection of a6--b6 and B--A);
\draw (a6)--(c6);
\draw (d6)--(b6);

\node (a7) at (-4.5,-2) [coordinate] {} ;
\node (b7) at (3,1) [coordinate] {} ;
\node (e7) at (-4,0) [coordinate] {} ;
\node (f7) at (-4,5) [coordinate] {} ;
\coordinate (c7) at (intersection of a7--b7 and D--E);
\coordinate (d7) at (intersection of a7--b7 and B--A);
\coordinate (g7) at (intersection of a7--b7 and e7--f7);
\draw (g7)--(c7);
\draw (d7)--(b7);

\node (a8) at (-4.5,-3) [coordinate] {} ;
\node (b8) at (3,0) [coordinate] {} ;
\node (e8) at (-3,0) [coordinate] {} ;
\node (f8) at (-3,5) [coordinate] {} ;
\coordinate (c8) at (intersection of a8--b8 and D--E);
\coordinate (d8) at (intersection of a8--b8 and B--A);
\coordinate (g8) at (intersection of a8--b8 and e8--f8);
\draw (g8)--(c8);
\draw (d8)--(b8);

\node (a9) at (-4.5,-4) [coordinate] {} ;
\node (b9) at (3,-1) [coordinate] {} ;
\node (e9) at (-2,0) [coordinate] {} ;
\node (f9) at (-2,5) [coordinate] {} ;
\node (m9) at (2.5,0) [coordinate] {} ;
\node (n9) at (2.5,5) [coordinate] {} ;
\coordinate (d9) at (intersection of a9--b9 and m9--n9);
\coordinate (g9) at (intersection of a9--b9 and e9--f9);
\draw (g9)--(d9);

\end{tikzpicture}
\end{center}
\caption{Les ensembles $F_i$ et $G_i$}
\end{figure}

Il existe des fonctions $\phi'_i$ , d\'efinies au voisinage de $\Delta$, de classe
${\cal C}^1$ au moins, telles que $\phi'_i$ soit positive dans $F_i$ et nulle dans
$G_i$. Les restrictions \`a $\Delta$ des fonctions $\phi_i = \frac{\phi'_i}{ \sum_j
\phi'_j}$ forment une partition de l'unit\'e $(U_i,\phi_i)$ o\`u $U_i$ d\'esigne le
compl\'ementaire de $G_i$.

Soit $\sigma$ un simplexe orient\'e de $\Delta$, de dimension $p$ et dont les
sommets sont $(a_0, \ldots, a_p)$; on d\'efinit la $p$-forme
diff\'erentielle de Whitney associ\'ee \`a $\sigma$ et \`a la partition de l'unit\'e
$(U_i,\phi_i)$, par la formule $$\tilde \omg (\sg ) = p! \sum_{i=0}^p
(-1)^{i}\phi_{i}d\phi_{0}\wedge
\cdots \wedge \widehat {d\phi_{i}}\wedge \cdots \wedge d\phi_{p}.$$ 

\begin{remarques} 
\begin{enumerate}
\item La forme diff\'erentielle $\tilde \omg (\sg )$ est de classe
\'egale \`a la classe des $\phi'_i$.
\item La forme de Whitney construite en utilisant les fonctions
coordonn\'ees barycentriques $\lambda_i$ sur $\sg$ \`a la place des $\phi_i$ est
la forme volume de $\sg$. Les $\lambda_i$ ne sont que
lipschitziennes sur $\Delta$; les $\phi_i$ servent \`a les lisser. 
\end{enumerate}
\end{remarques} 

Dans cet article, o\`u nous nous pla\c cons 
dans le cadre sous-analytique,  nous utiliserons les formes de Whitney ``non liss\'ees'' 
$$\omega(\sigma) = p!\sum_{i=0}^p (-1)^i
\lambda_i d\lambda_0 \wedge \cdots \wedge \widehat {d\lambda_i} \wedge \cdots \wedge
d\lambda_p.$$ 
Ceci est justifi\'e par la proposition : 

\begin{proposition}\libel{forsousa}
La forme de Whitney $\omg(\sg)$ a pour support l'\'etoile de $\sg$ dans $\Delta$, r\'eunion 
des \'etoiles des sommets de $\sg$. Elle est lin\'eaire par morceaux et continue, donc sous-analytique et continue. 
Sa diff\'erentielle au sens des
distributions admet un repr\'esentant sous-analytique born\'e.
\end{proposition}

\begin{preuve}
Chacune des fonctions $\lambda_i$ ayant pour support l'\'etoile du sommet $a_i$, la forme 
$\omg(\sg)$ a donc pour support la r\'eunion de ces \'etoiles, \`a savoir l'\'etoile de $\sg$. Les fonctions 
$\lambda_i$ sont lin\'eaires par morceaux et continues, les $1$-formes diff\'erentielles $d\lambda_i$, 
d\'efinies au sens des distributions admettent un repr\'esentant  sous-analytique et born\'e. 
Un calcul simple montre que 
\begin{equation}\label{diffomg}
d\omega(\sigma) = (p+1)!\ d\lambda_0 \wedge \cdots \wedge d\lambda_p.
\end{equation}
laquelle forme diff\'erentielle a pour support l'\'etoile de $\sg$ et jouit des m\^emes propri\'et\'es.
\end{preuve}

\subsection{Formes de Whitney sur les ensembles prismaux}

Soit $\rho = \sg_0 \times \cdots \times \sg_k$ un prisme orient\'e de l'ensemble prismal
$\Pi$. Notons $pr_j$ la projection de $\pi$ sur
$\sg_j$; on d\'efinit la forme de Whitney $\omg (\rho)$ par $$\omg
(\rho) = \ pr_0^*\omg (\sg_0) \wedge \cdots \wedge pr_k^*\omg (\sg_k).$$ 

Par la suite, nous omettrons les ``$pr_j$" sans que cela ne pr\^ete \`a ambigu\"\i t\'e.

En particulier, la forme de Whitney  du prisme $\pi(\sg) = \tau \times \sg_0 \times \cdots \times \sg_s$ s'\'ecrit 
\begin{equation}\label{omgpisg}
\omg(\pi(\sg)) = \omg(\tau) \wedge \omg(\sg_0)\wedge \cdots \wedge \omg(\sg_s)
\end{equation}

Repla\c cons nous dans la situation de la proposition \ref{JOINT} avec un morphisme $f: \Delta \to T$ 
que nous supposerons orient\'e. 
Soit $\sigma$ un $p$-simplexe orient\'e de $\Delta$, la forme de Whitney
$\omega(\sigma)$ s'\'ecrit : $$\omega(\sigma) = p!\sum_{i=0}^p (-1)^i \lambda_i d\lambda_0
\wedge \cdots \wedge \widehat {d\lambda_i} \wedge \cdots \wedge d\lambda_p$$ 
o\`u les $\lambda_i$ sont les coordonn\'ees barycentriques correspondant aux sommets $a_i$
de $\sigma$. Nous avons vu
que les composantes prismales d'un point de $\pi(\sigma)= \tau\times \sigma_0 \times \cdots \times \sigma_s$ 
s'\'ecrivent $\sum_{j=0}^s t_j y_j$ dans
$\tau$, et $\sum_{i \in I(j)}\mu_{j,i} a_i$ dans $\sigma_j$. L'application $\psi ^\sg $
associe \`a ce point le point de $\sigma$ ayant pour coordonn\'ees barycentriques les
$\lambda_{j,i} = t_j \mu_{j,i}$ (o\`u $j$ est tel que $i\in I(j)$). 
Les coordonn\'ees barycentriques $\lambda_i$ correspondant \`a la composante $\sg_j$ seront not\'ees 
$\lambda_{j,i_j}$ avec $i_j = 0, \ldots,  \alpha_j = \dim \sg_j$.

Regroupons les coordonn\'ees de $\sigma$ par coordonn\'ees correspondant \`a chacun des $\sigma_j$, 
$\omega(\sigma)$ devient :
\begin{eqnarray}
\omega(\sigma) = &(-1)^{a(\nu)} \  p! \ \sum_{j=0}^s 
\sum_{\ell =0}^{\alpha_j} (-1)^{\beta(j,\ell )} \ \lambda_{j,\ell }\ 
{\underbrace {d\lambda_{0,0}\wedge \cdots \wedge
d\lambda_{0,\alpha_0}}_{\sigma_0}}\wedge \cdots \wedge \nonumber \\
&{\underbrace {d\lambda_{j,0}\wedge \cdots \wedge \widehat
{d\lambda_{j,\ell }}\wedge \cdots \wedge
d\lambda_{j,\alpha_j}}_{\sigma_j}}\wedge \cdots \wedge {\underbrace
{d\lambda_{s,0}\wedge \cdots \wedge
d\lambda_{s,\alpha_s}}_{\sigma_s}}
\label{ex1}
\end{eqnarray} 
o\`u  l'ensemble des $p+1$ couples $(j,i_j)$ 
est en bijection avec $(0,1,\ldots,p)$ et 
$a(\nu)$ est la signature de la permutation 
$$\nu: (0,1,\ldots,p) \rightarrow ((0,0),\ldots
,(0,\alpha_0),(1,0),\ldots ,(1,\alpha_1),\ldots ,(s,0),\ldots ,(s,\alpha_s)) ;$$
enfin $\beta(j,\ell ) = \alpha_0 +\cdots +\alpha_{j-1} + j + \ell  $ ($0\le \ell \le \alpha_j$). Remarquons que 
$\nu$ exprime la compatibilit\'e de l'orientation de $\sg$ avec celle du joint $\sg_0 * \cdots *\sg_s$ des 
simplexes orient\'es $\sg_j$ (voir section \ref{section2}).

\begin{proposition}\libel{IMINVE} Reprenons les notations du
th\'eor\`eme \ref{COMPFA}.  Une orientation des simplexes $\sg_j$ et une orientation de $\tau$ 
induisent une orientation de $\sg$ comme joint it\'er\'e (voir section \ref{section2}).  On a l'\'egalit\'e   
$$(\psi ^\sg)^*(\omg (\sg)) = (-1)^{\alpha(\sigma,\nu)}\  
\frac{p!}{\vert \sigma_0\vert !\cdots \vert \sigma_s\vert ! \ s!} \ 
t_0^{\vert \sigma_0\vert}  \cdots t_s^{\vert \sigma_s\vert} \ \omg (\pi (\sg))$$  o\`u $\alpha(\sg,\nu) =
s\vert\sigma_0\vert + (s-1)\vert \sigma_1\vert + \cdots + \vert \sigma_{s-1}\vert + a(\nu)$ et les
$(t_j)_{j=0,\ldots ,s}$ sont les coordonn\'ees barycentriques de $\tau$.
\end{proposition}

\begin{preuve}
Tout d'abord, on remarque que, pour chacune des coordonn\'ees $\lambda_{j,i}$, on a :
\begin{equation}\label{pouf}
 (\psi ^\sg)^*(d\lambda_{j,i}) = dt_j\ \mu_{j,i} + t_j \ d\mu_{j,i}
 \end{equation}
d'autre part, pour tout $k=0,\ldots ,s$, on a :
$$\sum_{i\in I(k)} \mu_{k,i} =1\ \ \ {\rm donc}\ \ \ d\mu_{k,0}\wedge 
\cdots \wedge d\mu_{k,\alpha_k} =0, $$
o\`u $\alpha_k $ est $\dim \sg_k = \vert \sg_k \vert = \text{card}(I(k))$.

Reprenons l'expression de $\omg(\sg)$ donn\'ee par la formule (\ref{ex1}). 
Pour chaque terme de la somme, c'est-\`a-dire pour $j$ fix\'e, 
la contribution des produits des diff\'erentielles provenant de $\sg_j$ se calcule de fa\c con
diff\'erente de celle des produits des diff\'erentielles provenant des autres simplexes $\sg_k$. 
Plus pr\'ecis\'ement, au vu de (\ref{pouf}), nous v\'erifions la formule suivante 
pour $k\neq j$  :
\begin{equation}\label{paevi}
(\psi ^\sg)^*d\lambda_{k,0}\wedge \cdots \wedge d\lambda_{k,\alpha_k} = (t_k)^{\alpha_k}
\sum_{i_k=0}^{\alpha_k} (-1)^{i_k} \ \mu_{k,i_k} \  dt_k \wedge d\mu_{k,0}\wedge \cdots
\wedge \widehat {d\mu_{k,i_k}}\wedge 
\cdots \wedge d\mu_{k,\alpha_k}.
\end{equation}
Donc, pour $j$ fix\'e, et pour
$k$ diff\'erent de $j$,  tous les produits de diff\'erentielles associ\'es aux simplexes
$\sigma_k$ contiennent le terme $dt_k$ correspondant.

Pour d\'eterminer la
contribution des termes de la formule (\ref{ex1}) provenant de $\sg_j$, remarquons que 
\begin{equation}\label{paaf}
\sum_{j=0}^s t_j = 1,\ \ \  {\rm donc\  on\  a }\ \ \ dt_0\wedge 
\cdots \wedge dt_s =0.
\end{equation} 
Comme on vient de le voir, tous les termes $dt_k$, sauf $dt_j$,  apparaissent d\'ej\`a dans
l'expression  de $(\psi ^\sg)^*(\omg(\sg))$. 
La formule (\ref{paaf}) montre que les termes provenant de $\sg_j$ et qui 
contiennent $dt_j$ ont une contribution nulle. 
La seule contribution non
nulle du d\'eveloppement de  $$(\psi ^\sg)^*(d\lambda_{j,0}\wedge \cdots \wedge \widehat
{d\lambda_{j,\ell }}\wedge \cdots \wedge d\lambda_{j,\alpha_j}).$$ est donc~:
\begin{equation} 
(t_j)^{\alpha_j} d\mu_{j,0}\wedge \cdots \wedge \widehat {
d\mu_{j,\ell }}\wedge \cdots \wedge d\mu_{j,\alpha_j}. \label{ex3} 
\end{equation}

En reportant (\ref{paevi}) et (\ref{ex3}) dans (\ref{ex1}), il vient :
\[    \begin{array} {rl}
(\psi^\sg)^*  (\omg (\sg)) &= 
(-1)^{a(\nu)}  p! \ 
(t_0)^{\alpha_0}\cdots (t_s)^{\alpha_s} 
\sum_{j=0}^s (-1)^{\alpha_0 +\cdots +\alpha_{j-1}+j} t_j \\
&\left(\sum_{i_0=0}^{\alpha_0}(-1)^{i_0}\ \mu_{0,i_0}\ 
dt_0 \wedge d\mu_{0,0}\wedge \cdots \wedge \widehat {d\mu_{0,i_0}}\wedge 
\cdots \wedge d\mu_{0,\alpha_0}\right) \wedge \cdots \wedge \\
&\left(\sum_{\ell =0}^{\alpha_j}(-1)^{\ell }\ \mu_{j,\ell }\ 
d\mu_{j,0}\wedge \cdots \wedge \widehat {d\mu_{j,\ell }}\wedge 
\cdots \wedge d\mu_{j,\alpha_j}\right) \wedge \cdots \wedge \\
&\left(\sum_{i_s=0}^{\alpha_s}(-1)^{i_s}\ \mu_{s,i_s}\ 
dt_s \wedge d\mu_{s,0}\wedge \cdots \wedge \widehat {d\mu_{s,i_s}}\wedge 
\cdots \wedge d\mu_{s,\alpha_s}\right) 
\end{array}     \]

Chacune des parenth\`eses, sauf la $j$-\`eme, est 
\'egale \`a $\frac{1}{\alpha_k!} dt_k\wedge  \omg (\sg_k)$. 
Arriv\'e ici, nous savons que nous pouvons \'ecrire :
\[    \begin{array} {rl}
(\psi ^\sg)^*  (\omg (\sg))
& = (-1)^{a(\nu)}\ \frac{p!}{\alpha_0!\cdots \alpha_s!} \ 
(t_0)^{\alpha_0}\cdots (t_s)^{\alpha_s}\\
&\left(\sum_{j=0}^s (-1)^{\alpha (\sg)+j}t_j dt_0\wedge \cdots
\wedge \widehat {dt_j}\wedge \cdots \wedge dt_s \right) \wedge 
\omg (\sg_0)\wedge \cdots \wedge \omg (\sg_s), 
\end{array}     \]
o\`u nous avons not\'e $\alpha(\sg)= s\,\alpha_0 + (s-1) \,\alpha_1+ \cdots + \alpha_{s-1}$. 
Le signe est donn\'e par le nombre de permutations qui permettent cette \'ecriture de
$(\psi ^\sg)^*  (\omg (\sg))$. De fa\c con pr\'ecise : 
le terme $dt_0$ est en premi\`ere
place, il faut faire $\alpha_0$ permutations pour ramener $dt_1$ en
deuxi\`eme place, $\alpha_0+\alpha_1$ permutations pour ramener $dt_2$ en
troisi\`eme place, ainsi jusqu'\`a $dt_{j-1}$, lequel
n\'ecessite $\alpha_0+\cdots +\alpha_{j-2}$ permutations pour venir en
$j$-\`eme place. Ensuite, pour ramener $dt_{j+1}$ en $(j+1)$-\`eme place, il
faut $\alpha_0+\cdots +\alpha_j$ permutations, ainsi jusqu'\`a
$dt_s$ lequel n\'ecessite $\alpha_0+\cdots +\alpha_{s-1}$ permutations
pour venir en $s$-\`eme place.  On en d\'eduit  le r\'esultat.
\end{preuve}

\subsection{Formes de Whitney et formes relatives}

Consid\'erons le faisceau prismal ${\cal P}_f$ associ\'e \`a un morphisme
simplicial $f$. Soit $\pi \in {\cal P}_f(\tau)$, de la forme $\pi =\tau \times 
\sg_0 \times \cdots \times \sg_s$.\par
 Montrons maintenant comment $\omg (\sg)$ s'exprime en
fonction de $\omg (\sg')$ pour toute face $\sg'$ de $\sg$ et plus g\'en\'eralement $\omg (\rho)$ en fonction de $\omg
(\rho')$ pour toute face $\rho'$ d'un prisme $\rho$.\par
Fixons d'abord quelques notations : Soit $\sg$ un $p$-simplexe orient\'e de
sommets $a_0, \ldots ,a_p$  et notons $(\lambda_0,\ldots, \lambda_p)$ les coordonn\'ees
barycentriques correspondantes.  A  toute face 
$\sg'$ de $\sg$, de sommets $(a_{i_0},\ldots, a_{i_q})$, on associe la forme  diff\'erentielle 
suivante d\'efinie sur $\sg$ :
$$
\omg(\sg';\sg) =  q  !\sum_{k=0}^q (-1)^k \lambda_{i_k} d\lambda_{i_0}\wedge \cdots
\wedge \widehat{d\lambda_{i_k}}\wedge \cdots \wedge d\lambda_{i_q}. $$ 
La restriction de $\omg(\sg';\sg)$ \`a $\sg'$ est la forme de Whitney de $\sg'$ ; la
forme diff\'erentielle $\omg(\sg';\sg)$ n'est autre que l'extension \`a $\sg$ tout entier de 
l'\'ecriture de la forme de Whitney de $\sg'$. En ce sens elle constitue une extension canonique de 
$\omg(\sg')$ \`a $\sg$. 

\begin{remarque} On peut d\'efinir une telle extension pour tout simplexe de l'\'etoile ${\rm St}_\Delta\sigma'$ de
$\sigma'$, et les formes ainsi d\'efinies co\"\i ncident sur l'intersection de deux des
simplexes de cette \'etoile. On obtient ainsi une forme diff\'erentielle
$\omega(\sigma';\hbox{\rm St}_\Delta\sigma')$.
\end{remarque}

\begin{lemma}\libel{LEMCOD}
\begin{description}
\item [{\sl a)}] Pour une face $\sg'$ de codimension 1 d'un simplexe orient\'e $\sg$, on a:
$$d \omega(\sg' ;\sigma) = [\sg ; \sg'] \omg(\sigma).$$ 
\item [{\sl b)}] 
La forme de Whitney d'une face $\pi'$ de
codimension 1 d'un prisme orient\'e $\pi$ est la  restriction \`a cette face d'une forme
diff\'erentielle  $\omg (\pi';\pi)$ canoniquement d\'efinie sur ce prisme et dont la diff\'erentielle est
\'egale \`a la forme de Whitney du prisme, au facteur $[\pi ;\pi ']$ pr\`es, autrement dit :
$$d\omg(\pi' ;\pi)=[\pi ;\pi ']\ \omg (\pi ).$$  
\item [{\sl c)}] Notons $L(\pi)$ l'espace vectoriel des formes diff\'erentielles  de degr\'e
$\vert
\pi\vert -1$ dont les coefficients sont des fonctions lin\'eaires en les coordonn\'ees
barycentriques des simplexes $\sg_i$ constituant $\pi$ et dont la restriction \`a chaque
face de codimension 1 est un multiple scalaire de la forme de Whitney de cette face. La
collection des formes $\omg (\pi';\pi)$ pour toutes les faces $\pi'$  de codimension 1  de
$\pi$ forme une base de $L(\pi)$. 
\end{description}
\end{lemma}

\begin{preuve}
D\'emontrons d'abord l'\'egalit\'e du a). 
Sans perte de g\'en\'eralit\'e, on
peut supposer que les coordonn\'ees barycentriques de $\sigma$ sont
$\lambda_0,\lambda_1,\ldots ,\lambda_p$ et celles de $\sigma'$ sont $\lambda_0,\lambda_1,\ldots ,\lambda_{p-1}$. 
D'une part  la
forme de Whitney de $\sigma'$ s'\'ecrit  
$$\omega(\sigma')=(p-1)!\sum_{i=0}^{p-1} (-1)^i \lambda_i\  d\lambda_0 \wedge \cdots \wedge 
\widehat{d\lambda_i} \wedge \cdots \wedge d\lambda_{p-1}$$   et cette m\^eme \'ecriture
d\'efinit la forme $\omg (\sg' ;\sg)$ sur $\sg$. Il vient (voir (\ref{diffomg}))
$$
d \omg (\sg' ;\sg)= p ! d\lambda_0 \wedge \cdots \wedge d\lambda_{p-1}.
$$
D'autre part, 
en rempla\c cant $\lambda_p$ par $1 - \sum_{i=0}^{p-1} \lambda_i $ et $d\lambda_p$ par $- \sum_{i=0}^{p-1} d\lambda_i $, 
 dans l'expression de $\omg(\sg)$, on a
\begin{equation} \label{voleuc}
\omg(\sigma) = p!\sum_{i=0}^{p} (-1)^i \lambda_i \ d\lambda_0 \wedge \cdots \wedge \widehat{d\lambda_i} \wedge \cdots \wedge d\lambda_{p}
= (-1)^p p!\  d\lambda_0
\wedge \cdots \wedge  d\lambda_{p-1}.
\end{equation} 
On en d\'eduit le r\'esultat
puisque, dans ce cas $[\sg ; \sg']=(-1)^p$.

Prouvons b) et supposons que
$\pi'=\sigma_0\times \cdots \times \sigma'_i\times \cdots \times \sigma_s$ soit
une face de codimension 1 de $\pi=\sigma_0\times  \cdots \times \sigma_i\times
\cdots \times \sigma_s$, alors  
$$\omg(\pi';\pi) = \, \omg (\sg_0)
\wedge \cdots \wedge \omg (\sg'_i;\sg_i) \wedge \cdots \wedge \omg (\sg_s)$$ 
et donc  
$$d\omg(\pi';\pi) = \sum_{j=0}^s (-1)^{\vert \sg_0\vert  +
\cdots +\vert \sg_{j-1}\vert} \omg (\sg_0) \wedge \cdots \wedge d\omg
(\sg_j) \wedge \cdots \wedge \omg (\sg'_i;\sg_i) \wedge \cdots
\wedge \omg(\sg_s)$$    
(o\`u, bien entendu, pour $j=i$, $d\omg(\sg_j) = d\omg
(\sg'_i ;\sg_i)$).  Sur $\pi$, toutes les formes diff\'erentielles $d\omg
(\sg_j)$, pour $j\ne i$,  sont nulles car la somme des coordonn\'ees
barycentriques intervenant dans  $\sg_j$ est \'egale \`a 1 dans $\pi$. D'autre
part, le calcul pr\'ec\'edent  montre que $d\omg (\sg'_i ;\sg_i)= [\sg_i ;
\sg'_i]\omg (\sg_i)$, d'o\`u  le r\'esultat d'apr\`es le Lemme \ref{INKCID}. 

Prouvons l'assertion c). Puisqu'une forme de degr\'e maximum sur un simplexe, 
d\'e\-pen\-dant lin\'eairement des coefficients, est un multiple scalaire de la forme de Whitney, on a l'\'egalit\'e~:
$$L(\pi) = \bigoplus_{i=0}^s {\bf R} \omg(\sg_0) \wedge \ldots \wedge L(\sg_i) \wedge \ldots \wedge {\bf R} \omg(\sg_s).$$
Nous sommes donc ramen\'es au  cas d'un simplexe $\sg$. Soit $\omg \in L(\sg )$ ; sur chacun des simplexes $\sg'_k$ de
codimension 1 du bord de
$\sg$, on a $\omg\vert_{\sg'_k} = \lambda_k \, \omg(\sg'_k)$ et donc $\omg - \sum \lambda_k 
\, \omg(\sg'_k) =0$ sur le bord de $\sg$. Puisque, pour $k' \ne k$, on a 
$\omg(\sg'_k;\sg)\vert_{\sg'_{k'}} =0$, il vient $$(\omg - \sum_{\sg'_k \subset \partial \sg} \lambda_k \, 
\omg(\sg'_k;\sg))\Big\vert_{\partial \sg} =0 .$$
Puisque cette derni\`ere forme est lin\'eaire en les coordonn\'ees barycentriques et nulle sur
le bord de $\sg$, elle est nulle. Les formes de Whitney \'etendues sont lin\'eairement
ind\'ependantes sur $\sg$ puisque leurs restrictions au bord le sont.
\end{preuve}
\begin{rema}\libel{bord} Etant donn\'e un simplexe orient\'e $\sg$ de dimension $r$, notons $\sg'_i$ les composantes de son bord. Le a) du r\'esultat pr\'ec\'edent implique que si l'on consid\`ere la forme diff\'erentielle $\int \omega (\sg)=\frac{1}{r+1}\sum_{i=1}^{r+1}[\sg ;\sg'_i]\omega (\sg'_i;\sg)$, on a l'\'egalit\'e $d(\int \omega (\sg))=\omega (\sg)$. L'application qui \`a $\sg$ associe $\int \omega (\sg)$ peut \^etre vue comme la version "formes de Whitney" du bord.
\end{rema}
 
Consid\'erons maintenant le cas d'un simplexe $\gamma$ face de  $\sigma$ et appelons 
$\phi_h$ les faces de $\sigma$ admettant $\gamma$ pour face de codimension 1.
Plus pr\'ecis\'ement, notons $(\lambda_0, \lambda_1, \ldots , \lambda_\ell)$ les 
coordonn\'ees barycentriques du simplexe $\gamma$  
 dans $\sigma$ de coordonn\'ees barycentriques 
$(\lambda_0, \lambda_1, \ldots , \lambda_p)$. Pour tout 
$h = \ell+1, \ldots , p$, notons $\phi_h$ le simplexe de $\sigma$
de coordonn\'ees barycentriques 
$(\lambda_0, \lambda_1, \ldots , \lambda_\ell, \lambda_h)$.
 
 \begin{lemma}\libel{satrap} Avec les notations pr\'ec\'edentes, il vient  :
\begin{equation}\label{satrap1}
\sum_{h=\ell+1}^p \omega(\phi_h;\sigma) = (-1)^{\ell+1} (\ell+1)!\;
d\lambda_0
\wedge \cdots \wedge d\lambda_\ell. 
\end{equation}
\end{lemma}

\begin{preuve} On a :
$$\omega(\phi_h;\sigma) = (\ell+1)! \left(\sum_{j=0}^\ell (-1)^j \lambda_j
d\lambda_0
\wedge \cdots 
\wedge \widehat{d\lambda_j} \wedge \cdots \wedge d\lambda_\ell\wedge
d\lambda_h + (-1)^{\ell+1} \lambda_h d\lambda_0
\wedge \cdots \wedge d\lambda_\ell\right)$$
et donc
$$\sum_{h=\ell+1}^p \omega(\phi_h;\sigma) = (\ell+1)! 
\sum_{j=0}^\ell (-1)^j \lambda_j d\lambda_0 \wedge \cdots 
\wedge \widehat{d\lambda_j} \wedge \cdots \wedge d\lambda_\ell\wedge
\left(\sum_{h=\ell+1}^p d\lambda_h \right) 
$$
$$
+ (-1)^{\ell+1}(\ell+1)!  \left(\sum_{h=\ell+1}^p \lambda_h \right)  d\lambda_0
\wedge \cdots \wedge d\lambda_\ell 
$$
o\`u $\sum_{h=\ell+1}^p \lambda_h = 1 - \sum_{h=0}^\ell \lambda_h$, donc 
$\sum_{h=\ell+1}^p d\lambda_h = - \sum_{h=0}^\ell d\lambda_h$. Il vient
$$\sum_{h=\ell+1}^p \omega(\phi_h;\sigma) = (\ell+1)!
\sum_{j=0}^\ell (-1)^{\ell+1} \lambda_j d\lambda_0 \wedge \cdots 
\wedge d\lambda_\ell + (-1)^{\ell+1}(\ell+1)! \left(\sum_{h=\ell+1}^p
\lambda_h \right) d\lambda_0 \wedge \cdots \wedge d\lambda_\ell 
$$
et le r\'esultat.
\end{preuve}
Le lemme suivant est crucial pour la suite. Il \'enonce une \'egalit\'e de formes diff\'erentielles au sens des distributions ou des courants. Comme nous le verrons, la preuve du th\'eor\`eme repose sur des constructions explicites de repr\'esentants sous-analytiques des solutions d'\'equations impliquant des distributions.
 \begin{lemma}\libel{satrapaz} Avec les m\^emes notations
  que le lemme pr\'ec\'edent, \'etant donn\'ee une fonction $E$ d\'efinie, sous-analytique et born\'ee sur $\sigma$, notant toujours  
$\omega(\gamma;\sigma)$ l'extension \`a $\sigma$ de la forme de Whitney 
de $\gamma$, on l'\'egalit\'e au sens des distributions  :
\begin{equation}\label{satrap2}
d\left( E\; \omega(\gamma;\sigma)\right) =
(-1)^{\ell+1} \sum_{h=\ell+1}^p \left( 
E_h + \frac{1}{\ell +1} \sum_{i=0,i\ne h}^p \lambda_i \frac{\partial E_h}{\partial \lambda_i} \right)
\omg(\phi_h;\sg), \end{equation}
o\`u $E_h$ est la fonction de $p$ variables d\'efinie sur $\sg$ par 
$$E_h(\lambda_0, \ldots,\widehat{\lambda_h},\ldots,\lambda_p) =
E(\lambda_0, \ldots,1 - \Sigma_{i=0, i\ne h}^p \lambda_i, \ldots,\lambda_p).$$

\end{lemma}

\begin{preuve} 
On a 
$$d\left( E\; \omega(\gamma;\sigma)\right) =
dE\wedge \omega(\gamma;\sigma) + E\; d(\omega(\gamma;\sigma))$$
o\`u le premier terme s'\'ecrit
$$dE\wedge \omega(\gamma;\sigma)  = 
\left(\sum_{i=0}^p \frac {\partial E}{\partial \lambda_i} 
d\lambda_i \right) \wedge\left(
\ell !\;  \sum_{j=0}^\ell (-1)^j \lambda_j d\lambda_0 \wedge \cdots 
\wedge \widehat{d\lambda_j} \wedge \cdots \wedge d\lambda_\ell
\right).$$
S\'eparons cette formule en deux sommes relativement \`a l'indice $i$, la premi\`ere en sommant de $i
= 0$ \`a $\ell$ et la seconde de $\ell+1$ \`a $p$. 

Dans la premi\`ere somme et pour $i$ fix\'e, ou bien $i\ne j$, alors $d\lambda_i$
appara{\^\i t} dans le produit  $d\lambda_0 \wedge \cdots 
\wedge \widehat{d\lambda_j} \wedge \cdots \wedge d\lambda_\ell$ et dans ce
cas, sa contribution est nulle, ou bien $i=j$ et alors $d\lambda_i=d\lambda_j$ 
compl\`ete le produit avec $j$ permutations pour retrouver la $j$-\`eme place. 
La premi\`ere somme est donc \'egale \`a :
\begin{equation}\label{dun}
\ell !\;\sum_{i=0}^\ell \lambda_i \left(\frac {\partial E}{\partial
\lambda_i} \right) d\lambda_0 \wedge \cdots \wedge d\lambda_\ell.
\end{equation}

La deuxi\`eme somme, pour $i =\ell +1$ \`a $p$, est \'egale \`a
$$\ell ! (-1)^\ell\;\sum_{i=\ell+1}^p \left(\frac {\partial E}{\partial
\lambda_i} 
\right) \sum_{j=0}^\ell (-1)^j \lambda_j d\lambda_0 \wedge \cdots 
\wedge \widehat{d\lambda_j} \wedge \cdots \wedge d\lambda_\ell\wedge
d\lambda_i.$$
Comme on a 
$$\omega(\phi_i;\sigma) = (\ell+1)! \left(\sum_{j=0}^\ell (-1)^j \lambda_j
d\lambda_0
\wedge \cdots 
\wedge \widehat{d\lambda_j} \wedge \cdots \wedge d\lambda_\ell\wedge
d\lambda_i + (-1)^{\ell+1} \lambda_i d\lambda_0
\wedge \cdots \wedge d\lambda_\ell\right)$$
autrement dit
$$\left(\sum_{j=0}^\ell (-1)^j \lambda_j d\lambda_0 \wedge \cdots 
\wedge \widehat{d\lambda_j} \wedge \cdots \wedge d\lambda_\ell\wedge
d\lambda_i\right) = \frac{1}{(\ell+1)!} \omega(\phi_i;\sigma) - 
(-1)^{\ell+1} \lambda_i d\lambda_0
\wedge \cdots \wedge d\lambda_\ell,$$
la deuxi\`eme somme est donc \'egale \`a 
$$\ell ! (-1)^\ell\;\sum_{i=\ell+1}^p \left(\frac {\partial E}{\partial
\lambda_i} 
\right) \left(
\frac{1}{(\ell+1)!}\omega(\phi_i;\sigma) - (-1)^{\ell+1} \lambda_i d\lambda_0
\wedge \cdots \wedge d\lambda_\ell
\right)
$$
que nous d\'ecomposons comme suit :
\begin{equation}\label{ddedeux}
(-1)^\ell \frac{1}{\ell +1}
\; \sum_{i=\ell+1}^p \left(\frac {\partial E}{\partial
\lambda_i} 
\right)  \omega(\phi_i;\sigma)
\end{equation}
\begin{equation}\label{detrois}
+ \ell ! \;\sum_{i=\ell+1}^p\left(\frac {\partial E}{\partial
\lambda_i} \right)  \lambda_i d\lambda_0 \wedge \cdots \wedge d\lambda_\ell.
\end{equation}
En utilisant le lemme \ref{satrap}, la somme de (\ref{dun}) et (\ref{detrois}) est \'egale \`a
\begin{equation}\label{dequatre}
\ell ! \;\left(\sum_{i=0}^p  \lambda_i\frac {\partial E}{\partial
\lambda_i} \right)  d\lambda_0 \wedge \cdots \wedge d\lambda_\ell = 
(-1)^{\ell+1} \frac{1}{\ell +1} \;\left(\sum_{i=0}^p  \lambda_i\frac{\partial E}{\partial
\lambda_i} \right) \left(\sum_{h=\ell+1}^p \omega(\phi_h;\sigma) \right).
\end{equation}

Enfin,  par sommation de (\ref{ddedeux}) et (\ref{dequatre}), on obtient
$$d E\wedge \omega(\gamma;\sigma) = (-1)^{\ell+1} \frac{1}{\ell +1}
\;\sum_{h=\ell+1}^p \left(\sum_{i=0}^p  
\lambda_i\frac {\partial E}{\partial \lambda_i} -  
\frac {\partial E}{\partial \lambda_h}\right) \omega(\phi_h;\sigma).
$$
Maintenant 
$$E \; d\omega(\gamma;\sigma)  = E \;
\ell !\; d\left( \sum_{j=0}^\ell (-1)^j \lambda_j d\lambda_0 \wedge \cdots 
\wedge \widehat{d\lambda_j} \wedge \cdots \wedge d\lambda_\ell \right)
= E \;
\ell !\; (\ell+1) \; d\lambda_0 \wedge \cdots \wedge d\lambda_\ell 
$$
$$=E\;(\ell+1)! \; d\lambda_0 \wedge \cdots \wedge d\lambda_\ell
= (-1)^{\ell+1} E \sum_{h=\ell+1}^p \omega(\phi_h;\sigma)
$$
en utilisant (\ref{satrap1}). On obtient donc finalement
$$d\left( E\; \omega(\gamma;\sigma)\right) =
(-1)^{\ell+1} \sum_{h=\ell+1}^p \left( E + \frac{1}{\ell +1}
\left(\sum_{i=0}^p  
\lambda_i\frac {\partial E}{\partial \lambda_i} -  
\frac {\partial E}{\partial \lambda_h}\right)
\right)\omega(\phi_h;\sigma).
$$
On peut encore \'ecrire diff\'eremment cette formule, en remarquant que, sur $\sg$, 
les variables $\lambda_i$ sont d\'ependantes : on a 
$$\lambda_h = 1 - \sum_{i=0, i\ne h}^p \lambda_i.$$
Notant 
$$E_h(\lambda_0, \ldots,\widehat{\lambda_h},\ldots,\lambda_p) =
E(\lambda_0, \ldots,1 - \Sigma_{i=0, i\ne h}^p \lambda_i, \ldots,\lambda_p), $$
on a, pour $i\ne h$,  
$$\frac{\partial E_h}{\partial \lambda_i} =
\frac{\partial  E}{\partial \lambda_i} - \frac{\partial  E}{\partial \lambda_h}.$$
On a  alors 
\begin{align*}
E & + \frac{1}{\ell +1}  \left(\sum_{i=0}^p \lambda_i\frac {\partial E}{\partial \lambda_i} -  
\frac {\partial E}{\partial \lambda_h}\right)\\
& = \,  E_h + \frac{1}{\ell +1} \left(\sum_{i=0,i\ne h}^p \lambda_i \left( 
\frac{\partial E_h}{\partial \lambda_i} + \frac {\partial E}{\partial \lambda_h} \right) 
+ \lambda_h \frac {\partial E}{\partial \lambda_h} - \frac {\partial E}{\partial \lambda_h}\right)
=  \, E_h + \frac{1}{\ell +1} \sum_{i=0,i\ne h}^p \lambda_i \frac{\partial E_h}{\partial \lambda_i} .
\end{align*}
et le r\'esultat.
\end{preuve}

 Remarquons que la donn\'ee d'une orientation sur $\sg$ \'equivaut \`a la donn\'ee d'une 
 orientation sur $\sg'$, d'une orientation sur la face oppos\'ee $\sg''$ et d'un  ordre sur le couple 
 $(\sg',\sg'')$. 
La fonction $u_{ \sg'}= \sum_{k=0}^q \lambda_{i_k}$  est d\'efinie  sur
tout le simplexe $\sg$ et est \'egale \`a 1 sur $\sg'$. On a :
$$\frac{du_{\sg'}}{u_{\sg'}(1-u_{\sg'})}  = \frac{du_{\sg'}}{u_{\sg'}} - \frac{du_{\sg''}}{u_{\sg''}}
= - \frac{du_{\sg''}}{u_{\sg''}(1-u_{\sg''})} $$
o\`u la fonction $u_{ \sg''}= 1- u_{\sg'}$ est   d\'efinie  sur
tout le simplexe $\sg$ et est \'egale \`a 1 sur $\sg''$. 

\begin{proposition}\libel{faceface} Soient $\sg'$ une face de  dimension $q$ du simplexe 
orient\'e $\sg$ de  dimension $p$ et $\sg''$ la face oppos\'ee, on a l'\'egalit\'e :
$$
\omg (\sg) =  
 {\frac{(-1)^{p}}{p-q}}
\begin{pmatrix}
p \\
q
 \end{pmatrix}
\omg (\sg';\sg)\wedge \omg
(\sg'';\sg)\wedge \left( {\displaystyle{\frac{du_{\sg'}}{u_{\sg'}(1-u_{\sg'})} 
}   }\right).
$$
\end{proposition}

Remarquons que cette formule traduit le fait que, \`a un coefficient pr\`es, le volume du simplexe
est \'egal au volume de deux faces oppos\'ees par le volume normalis\'e d'un segment qui les joint.

\begin{preuve}
Nous pouvons supposer, sans perte de g\'en\'eralit\'e que les simplexes orient\'es 
$\sg'$ et $\sg''$ ont pour sommets respectifs 
les points $a_0, a_1, \ldots, a_q$ et $a_{q+1}, \ldots,a_p$ pris dans cet ordre. 
Ecrivons, en fonction des coordonn\'ees barycentriques correspondantes l'expression
$$\omg(\sg';\sg) \wedge \omg(\sg'';\sg) \wedge 
\left( \frac {d(\lambda_0 + \cdots + \lambda_q)} {\lambda_0 + \cdots + \lambda_q} 
- \frac{d(\lambda_{q+1} + \cdots + \lambda_p)}{\lambda_{q+1} + \cdots + \lambda_p} \right)$$
On obtient:
$$  \begin{array} {rl}
q ! \displaystyle \sum_{i=0}^q (-1)^i \lambda_i  & d\lambda_0 \wedge  \ldots \wedge {\widehat {d\lambda_i}} \wedge \ldots \wedge d\lambda_q \\
 \bigwedge &
(p-q-1)!   \displaystyle \sum_{j=q+1}^p (-1)^j \lambda_j d\lambda_{q+1} \wedge \ldots \wedge {\widehat {d\lambda_j}} \wedge \ldots \wedge d\lambda_p \\
&  \bigwedge \displaystyle
\left( \frac {d\lambda_0 + \cdots + d\lambda_q} {\lambda_0 + \cdots + \lambda_q} - \frac{d\lambda_{q+1} 
+ \cdots + d\lambda_p}{\lambda_{q+1} + \cdots + \lambda_p} \right) = \\
{}\\
q !  (p-q-1)! &  \displaystyle \sum_{i=0}^q  \sum_{j=q+1}^p (-1)^{i+j}  \frac{\lambda_i \lambda_j  } {\lambda_0 + \cdots + \lambda_q} \\
 (d\lambda_0 \wedge  \ldots & \wedge {\widehat {d\lambda_i}} \wedge \ldots \wedge d\lambda_q)
\wedge (d\lambda_{q+1} \wedge \ldots \wedge {\widehat {d\lambda_j}} \wedge \ldots \wedge d\lambda_p) \wedge 
(d\lambda_0 + \cdots + d\lambda_q)\\ {}\\
+ \ 
 q !  (p-q-1)! &  \displaystyle \sum_{i=0}^q  \sum_{j=q+1}^p (-1)^{i+j+1}  \frac{\lambda_i \lambda_j  } {\lambda_{q+1} + \cdots + \lambda_p} \\
 (d\lambda_0 \wedge  \ldots \wedge& {\widehat {d\lambda_i}} \wedge \ldots \wedge d\lambda_q)
\wedge (d\lambda_{q+1} \wedge \ldots \wedge {\widehat {d\lambda_j}} \wedge \ldots \wedge d\lambda_p) \wedge 
(d\lambda_{q+1} + \cdots + d\lambda_p).\\
\end{array}    $$

Etudions les deux termes obtenus: dans le premier terme, et pour $i$ fix\'e, le produit ext\'erieur  de 
$ (d\lambda_0 \wedge  \ldots \wedge {\widehat {d\lambda_i}} \wedge \ldots \wedge d\lambda_q)$ et $(d\lambda_0 + \cdots + d\lambda_q)$ 
est nul sauf \`a produire un \'el\'ement $ (d\lambda_0 \wedge  \ldots \wedge { {d\lambda_i}} \wedge \ldots \wedge d\lambda_q)$. 
Le nombre de permutations pour ramener le terme $d\lambda_i$ de la derni\`ere \`a la $i$-\`eme place est \'egal 
\`a $(q-i+1) + (p-q-1) = p-i$. On en d\'eduit que ce premier 
terme est \'egal \`a:
$$ q !  (p-q-1)!  \sum_{i=0}^q  \sum_{j=q+1}^p (-1)^{p+j}  \frac{\lambda_i \lambda_j  } {\lambda_0 + \cdots + \lambda_q} 
 d\lambda_0 \wedge  \ldots \wedge {{d\lambda_i}} \wedge \ldots \wedge d\lambda_q
\wedge d\lambda_{q+1} \wedge \ldots \wedge {\widehat {d\lambda_j}} \wedge \ldots \wedge d\lambda_p, $$
c'est-\`a-dire
$$q !  (p-q-1)!  \sum_{j=q+1}^p (-1)^{p+j}  \lambda_j  \, 
 d\lambda_0 \wedge  \ldots \wedge  {\widehat {d\lambda_j}} \wedge \ldots \wedge d\lambda_p. $$
 
 De la m\^eme fa\c con, et cette fois en fixant $j$, on v\'erifie que le second terme est \'egal \`a:
 $$q !  (p-q-1)!  \sum_{i=0}^q (-1)^{p+i}  \lambda_i \,
 d\lambda_0 \wedge  \ldots \wedge  {\widehat {d\lambda_i}} \wedge \ldots \wedge d\lambda_p. $$

Notre expression est donc finalement \'egale \`a
$$q !  (p-q-1)!  (-1)^{p} \sum_{i=0}^p (-1)^{i}  \lambda_i \, 
 d\lambda_0 \wedge  \ldots \wedge  {\widehat {d\lambda_i}} \wedge \ldots \wedge d\lambda_p = 
 \frac{q !  (p-q-1)! }{p!} (-1)^{p} \omg(\sg). $$
\end{preuve}

\begin{corollary}\libel{facepri}
Dans le cas particulier o\`u $\sg'$ est une face de codimension 1 et $\sg''$ est le point $a_p$ 
la formule s'\'ecrit
\[   
\omg (\sg) = (-1)^{p+1} p\; 
\omg (\sg';\sg)\wedge   {\displaystyle{\frac{d\lambda_p}{1-\lambda_p} }}
   \]
\end{corollary}

\begin{corollary}\libel{facepro}Plus g\'en\'eralement,on peut exprimer, pour tout prisme $\pi=\sg_0 \times \cdots \times
\sg_s$, la forme $\omg (\pi)$  en fonction de la forme de Whitney d'un prisme
$\pi'= \sg'_0 \times \cdots \times \sg'_s$ du bord de $\pi$ :
$$\omg (\pi) =  (-1)^a 
\omg(\pi';\pi)\wedge \omg(\pi'';\pi) \bigwedge_{j=0}^s 
 \left( {\displaystyle{\frac{du_{\sg_j'}}{u_{\sg_j'}(1-u_{\sg_j'})} 
}   }\right)
$$  
o\`u $\pi'' = \sg_0'' \times \sg_1'' \times \cdots \times \sg_s''$ est le prisme produit des faces oppos\'ees 
aux $\sg_i'$ dans les $\sg_i$ et  
orient\'es par l'ordre des leurs sommets, $u_{\sg_i'}$ est la fonction  somme des coordonn\'ees barycentriques
relatives \`a $\sg_i'$ et 
$a= \displaystyle\sum_{i=1}^s \vert \sg'_i\vert \left(\displaystyle \prod_{j=0}^{i-1} \vert \sg''_j\vert +1 \right).$
\end{corollary}

\begin{preuve}
La preuve consiste \`a appliquer la proposition \ref{faceface} pour chacun des simplexes $\sg_i$ du prisme
$\pi$.
\end{preuve}

\section{Formes diff\'erentielles sur les faisceaux prismaux}
\setcounter{equation}{0}

\subsection{D\'efinitions}

\begin{definition}\libel{FODIRE} Soient $\Rho$ un
ensemble prismal et ${\cal F}$ un faisceau prismal sur $\Rho$. Comme en 
\ref{FAISPR}, notons $e_\rho\colon {\cal F}(\rho )\to \rho$ la projection associ\'ee \`a chaque 
prisme $\rho$ de $\Rho$. On appelle $r$-forme diff\'erentielle  {\sl $\mu$-r\'eguli\`ere} (resp. sous-analytique) sur ${\cal F}(\rho)$ la donn\'ee
d'une $r$-forme diff\'erentielle $\mu$-r\'eguli\`ere
(resp. sous-analytique) $\omega_\pi$ sur chaque prisme $\pi$ de ${\cal F}(\rho)$ de telle fa\c con que les restrictions 
de ces formes diff\'erentielles aux faces communes \`a deux prismes de ${\cal F}(\rho)$ 
co\"\i ncident. 

On appelle $r$-forme diff\'erentielle {$\mu$-r\'eguli\`ere} sur le faisceau prismal 
${\cal F}$ la donn\'ee d'une forme diff\'erentielle $\mu$-r\'eguli\`ere sur chacun des ${\cal F}(\rho)$
de fa\c con \`a ce que, pour chaque face
$\rho'$ de $\rho$, la forme $\omega_\rho$ restreinte \`a l'image r\'eciproque
de $\rho'$ par $e_\rho$ co\"\i ncide avec $\omega_{\rho'}$. 
Rappelons que d'apr\`es la d\'efinition \ref{FAISPR}, si
$\rho'$ est une face de $\rho$, alors $e_\rho^{-1}(\rho')$ est une r\'eunion de
prismes de ${\cal F}(\rho)$.  Nous noterons $\Omega^r_{\mu}({\cal F})$ l'espace
vectoriel des $r$-formes diff\'erentielles $\mu$-r\'eguli\`eres sur le faisceau
prismal ${\cal F}$.\end{definition} \par\medskip

En particulier, lorsque ${\cal F}$ est le faisceau ${\cal T}$ tel que ${\cal
T}(\rho)=\rho$ pour tout $\rho$, l'espace vectoriel $\Omega^0_{\mu}({\cal
T})$ muni de la multiplication des fonctions est appel\'e {\it alg\`ebre des fonctions
$\mu$-r\'eguli\`eres sur $P$} et not\'e ${\cal O}_{P,\mu}$. 
Pour chaque entier $r$, l'espace
vectoriel $\Omega^r_{\mu}({\cal T})$ est en fait un ${\cal O}_{P,\mu}$-module.

Dor\'enavant, nous ne consid\'ererons plus que des formes 
sous-analytiques, c'est-\`a-dire sa\-tis\-faisant dans un voisinage
ouvert de chaque prisme $\rho$ de $\Rho$ dans son affine ambiant les conditions de la 
Proposition \ref{hreg} et
supprimerons donc l'indice $\mu$. Pour chaque prisme
$\rho\subset \Rho$, les espaces vectoriels $\Omega^r({\cal F}(\rho))$, munis des morphismes induits
par la diff\'erentielle au sens des distributions, forment un complexe not\'e
$\Omega^{\bullet}({\cal F(\rho)})$. Puisque la diff\'erentielle ne commute pas aux morphismes de sp\'ecialisation, 
on ne peut pas mettre une structure de complexe diff\'erentiel sur la famille des espaces vectoriels $\Omega^{r}({\cal F})$. 

Notons $e^*_{\cal F} : \Omega^{\bullet} ({\cal T}) \to 
\Omega^\bullet({\cal F})$   le morphisme associ\'e \`a $e$ qui pour
tout  prisme $\rho$ co\"\i ncide avec le morphisme naturel 
 $e^*_{{\cal F},\rho} : \Omega^{\bullet} ({\cal T}(\rho)) \to 
\Omega^\bullet({\cal F}(\rho))$. Remarquons que, \'etant donn\'e un morphisme 
$\chi : {\cal F}' \to {\cal F}$ de faisceaux prismaux sur $\Rho$, on a 
\begin{equation}\label{estar}
e^*_{{\cal F}'} = 
\chi^* \circ e^*_{\cal F}.
\end{equation}

D\'efinissons la forme diff\'erentielle de degr\'e
$\vert \rho\vert$
$$de_\rho = de_{{\cal F},\rho} =e^*_{{\cal F},\rho} \omg(\rho).$$ 

Puisque la forme $de_\rho$ est ferm\'ee, les noyaux $K^\bullet_\rho = {\rm Ker}(\wedge de_\rho)$
forment un sous-complexe   
de  $\Omega^\bullet({\cal F}(\rho))$.

D'apr\`es le corollaire \ref{facepro}, pour toute face $\rho'$ de $\rho$, la restriction de
$de_\rho\vert_{e_\rho^{-1}(\rho')}$ est de la forme $de_{\rho'}\wedge \theta$. Par cons\'equent,
l'image de  $K^r_\rho$ par l'homomorphisme de sp\'ecialisation $\Omega^r({\cal F}(\rho))\to
\Omega^r({\cal F}(\rho'))$  est contenue dans $K^r_{\rho'}$. 
D'o\`u nous obtenons  par passage aux quotients $\Omega^r_{rel}({\cal F}(\rho))= \Omega^r({\cal F}(\rho)) / K^r_\rho$, 
un morphisme  
$$\Omega^r_{rel}({\cal F}(\rho))  \longrightarrow 
\Omega^r_{rel}({\cal F}({\rho'})).$$

Nous appellerons $r$-forme diff\'erentielle relative
la donn\'ee pour chaque simplexe $\rho$  d'un \'el\'ement de l'espace vectoriel
quotient $\Omega^r_{rel}({\cal F}(\rho))$, de mani\`ere compatible avec les
homomorphismes de sp\'ecialisation pr\'ec\'edents. 
\begin{lemma}\libel{LECARE} Soient   $\omg$ et $\omg'$ deux $r$-formes
diff\'erentielles sur un faisceau prismal
${\cal F}$. Pour tout prisme $\rho$ et pour $t\in \rond\rho$ notons $i_t$ l'inclusion de la fibre $e_\rho^{-1}(t)$ dans 
${\cal F}(\rho)$, alors 
$d e_\rho \wedge \omg=0$ si et seulement si, pour
tout
$t\in
\rond\rho$, on a : $$i_t^*\omg = 0.$$
\end{lemma}
\begin{preuve}
Cela r\'esulte de l'\'ecriture des formes diff\'e\-ren\-tiel\-les dans
les coordonn\'ees ba\-ry\-cen\-tri\-ques du prisme. D'une part,   les compos\'ees avec $e_\rho$ des coordonn\'ees 
barycentriques de $\rho$ font partie d'un syst\`eme de coordonn\'ees barycentriques du prisme $\pi$,
et d'autre part, pour tout facteur $\tau$ de $\rho$, la forme de Whitney $\omg(\tau)$ s'\'ecrit 
$(-1)^s s!\  dt_0
\wedge \cdots \wedge  dt_{s-1}$ (voir (\ref{voleuc})).
\end{preuve}

\begin{remarque} \libel{LECARA} 
Toute r\'etraction $\chi$ de $\rond\rho$ sur un point $t$ de $\rond\rho$ induit une r\'etraction de ${\cal F}(\rond\rho)$ sur la fibre $e_\rho^{-1}(t)$ 
et donc une r\'etraction lin\'eaire  $\Lambda^r(d\chi): \Lambda^r(T^*(e_\rho^{-1}(t))) \to  \Lambda^r(T^*{\cal F}(\rond\rho))$. Avec ces notations, on a 
$$(\omg\vert_{e_\rho^{-1}(t)} - \chi(i_t^*\omg)) \wedge de =0.$$
Choisissons une r\'etraction  $\chi$ dans chaque simplexe, nous nous permettrons lorsqu'il s'agira 
de formes relatives, d'identifier la restriction \`a une fibre et l'image inverse sur une fibre. 
\end{remarque}

On appellera diff\'erentielle relative 
et on notera $d_e$ la diff\'erentielle induite par la diff\'erentielle $d$ dans les quotients 
$\Omega^r_{rel}({\cal F}(\rho))$.

\begin{definition} \libel{formvert}
Appelons forme diff\'erentielle {\it verticale} sur un prisme $\pi$ au dessus de $\tau$ une forme
diff\'erentielle qui est l'image r\'eciproque par la projection $
\tau\times\sigma_0\times  \cdots \times \sigma_s\to \sigma_0\times  \cdots \times \sigma_s$ d'une forme diff\'erentielle
sur $\sigma_0\times  \cdots \times \sigma_s$. Toute forme diff\'erentielle verticale d\'efinit naturellement une forme
diff\'erentielle relative.
\end{definition}

\subsection{Etude des formes de Whitney relatives dans ${\cal P}_f$}

Soit ${\cal P}_f$ le faisceau prismal associ\'e \`a un morphisme simplicial orient\'e 
$f\colon  \Delta \to T$ comme dans le th\'eor\`eme \ref{COMPFA}, et $\pi=\pi(\sg)$ un prisme de
${\cal P}_f$ d'images $\sg \subset\Delta$ et $\tau\subset T$. Notons
$y_0,\ldots,y_s$ les sommets de $\tau$ et $\sigma_j=\sg \cap e^{-1}(y_j)$. Le but de cette section est la construction de
formes de Whitney relatives  engendrant la cohomologie des fibres en tout
degr\'e. Le cas des formes de degr\'e \'egal \`a la dimension relative des simplexes est facile puisque la chute de
dimension au bord de $\tau$ implique la nullit\'e de la restriction de la forme ; c'est l'objet de cette section~:

\begin{lemma}\libel{OPICSF}
La donn\'ee pour chaque $\pi=\pi(\sg)\in {\cal P}_f$ de la forme diff\'erentielle verticale 
\begin{equation}\label{ex4}
\omg(\pi(\sg)/\tau)= \omg (\sg_0)\wedge\omg (\sg_1)\wedge \ldots \wedge\omg (\sg_s) 
\end{equation}
d\'efinit une forme diff\'erentielle relative de degr\'e
$r=\dim_{\rm rel}\sg$  sur  ${\cal P}_f$ au dessus de $T$.
\end{lemma}

\begin{preuve}
Etant donn\'ees des orientations de $\Delta$ et  de $T$  telles que le morphisme $f$ soit 
orient\'e, les formes diff\'erentielles $\omg(\sg)$ associ\'ees aux simplexes $\sg$ ayant m\^eme image
$\tau$ et m\^eme dimension relative se recollent en une forme diff\'erentielle sur $f^{-1}(\tau)$. 
L'orientation de $f$ induit une orientation de la fibre type $\sg_0 \times \cdots \times \sg_s$ de chaque 
simplexe et ces orientations sont compatibles entre elles. Au moyen de l'isomorphisme 
$\theta^\sg$, on en d\'eduit que les $\omg(\pi(\sg))$ se recollent au dessus de $\tau$, et donc 
que les formes diff\'erentielles $\omg(\pi/\tau)$ se recollent dans la fibre. 
\end{preuve}

\begin{lemma}\libel{OPICSH}
\begin{description}
\item [{\sl a)}] La donn\'ee sur chaque simplexe $\sg$ de $\Delta$ de la forme diff\'erentielle 
\begin{equation}
\label{equaOPI}
t_0^{\vert \sg_0\vert}\cdots t_s^{\vert 
\sg_s\vert}\omg(\sg_0)\wedge\omg (\sg_1)\wedge \ldots \wedge\omg (\sg_s) 
= t_0^{\vert \sg_0\vert}\cdots t_s^{\vert \sg_s\vert} \omg(\pi(\sg)/\tau)
\end{equation}
d\'etermine une forme diff\'erentielle relative sur ${\cal P}_f$. 
\item [{\sl b)}] La donn\'ee sur chaque simplexe $\sg$ de $\Delta$ d'une fonction sous-analytique 
$A_\sg$ prenant les m\^emes valeurs sur les faces communes \`a deux simplexes d\'etermine 
une forme diff\'e\-ren\-tiel\-le relative sur ${\cal P}_f$ dont la valeur sur $\pi(\sg)$ est
$$(A_\sg \circ \psi^\sg) \ t_0^{\vert \sg_0\vert}\cdots t_s^{\vert \sg_s\vert}  \ \omg(\pi(\sg)/\tau).$$
\item [{\sl c)}] 
D'apr\`es la Proposition \ref{IMINVE} 
et avec ses notations, on a 
\begin{equation}(\psi^\sg)^*(\omg (\sg ))=(-1)^{\alpha(\sigma,\nu)} 
\frac{p!}{\vert \sigma_0\vert !\cdots \vert \sigma_s\vert ! \ s!} \ 
t_0^{\vert \sg _0\vert}\cdots 
t_s^{\vert  \sg_s\vert}\omg (\tau)\wedge \omg (\pi (\sg)/\tau).\label{whitrel}
\end{equation} 
\end{description}
\end{lemma}
\begin{preuve}
Prouvons a). Le lemme \ref{OPICSF} montre que, au dessus de chaque simplexe $\tau$ la forme (\ref{equaOPI})
d\'etermine une forme diff\'erentielle relative au dessus de $\tau$. 
Soit maintenant $\tau'$ une face de $\tau$ au dessus de laquelle la dimension relative de $\sg$ ne varie pas,  
cela signifie que pour tout sommet $y_i$ de $\tau$ qui 
n'est pas dans $\tau'$, on a $\vert \sg_i \vert = 0$ (exemple du $2$-simplexe situ\'e \`a mi-hauteur dans la Figure 4). 
Posons $\pi' = \pi\vert_{\tau'}$, alors $\omg(\pi/\tau)\vert_{\pi'} = \omg(\pi'/\tau')$. Dans les autres cas, 
l'un au moins des $\sg_i$ n'est pas nul et donc notre forme diff\'erentielle s'annule au dessus de $\tau'$.

Le b) r\'esulte aussit\^ot du a) et le c) de la proposition \ref{IMINVE}.
\end{preuve}

\begin{lemma}\libel{OPICSG}
Posons $\pi = \pi(\sg)$, on a pour tout $t\in \tau$ l'\'egalit\'e 
$$\int_{\pi (t)}\omg(\pi /\tau) = 1.$$
\end{lemma}

\begin{preuve}
En effet $\omg(\pi /\tau)$ est la forme volume produit des formes volumes des simplexes 
verticaux du prisme $\pi$. Cela r\'esulte aussi du 2) de la Remarque \ref{REM}.
\end{preuve}

\begin{definition}\libel{DIFEXT}   {\sl Extension verticale
d'une forme diff\'erentielle relative}. Soit $\pi'$ une face de codimension 1 d'un prisme
$\pi$, notons $\tau'$ l'image de $\pi'$. De m\^eme que dans le cas absolu, on d\'efinit
l'extension \`a $\pi$ d'une  forme diff\'erentielle relative sur $\pi'$ comme la forme
diff\'erentielle $\omg(\pi'/\tau';\pi)$  donn\'ee par la m\^eme \'ecriture en coordonn\'ees barycentriques. 
C'est une forme diff\'erentielle relative au morphisme de  $\pi$ sur son image $\tau$ 
dont $\tau'$ est une face. 
\end{definition}

Remarquons qu'une face de codimension relative 1 d'un prisme $\pi=\tau\times\sg_0\times
\cdots\times \sg_s$ de ${\cal P}_f$ est n\'e\-ces\-sai\-re\-ment  obtenue de la fa\c con
suivante~: on remplace un des simplexes $\sg_j$ par une de ses faces de codimension 1 et on
fait le produit avec les autres. A chaque face de codimension relative 1 est donc associ\'e un
sommet bien d\'etermin\'e de $\tau$.

\begin{lemma}\libel{LEMROL} Soit $\pi'$ un prisme de ${\cal
P}_f$ d'image $\tau'$ et soit $\pi$ un prisme dont $\pi'$ est une face de codimension 1. On
a, avec les notations du lemme \ref{LEMCOD}~:  
$$d_e\omg (\pi' /\tau';\pi)=
\begin{cases}
[\pi;\pi'] \omg (\pi/\tau') & \text{ si $\pi$ a pour image $\tau'$,} \\
0& \text{ sinon}
\end{cases}$$
\end{lemma}

\begin{preuve}
Dans le cas o\`u $\pi$ a pour image $\tau'$, il s'\'ecrit 
$\pi=\tau'\times\sg_0\times \cdots\times \sg_j\times \cdots \times \sg_s$ et le
prisme $\pi'$ s'\'ecrit $\pi'=\tau'\times\sg_0\times \cdots\times \sg'_j\times 
\cdots \times \sg_s$. Alors on a
$$\omg (\pi' /\tau';\pi)= \omg(\sg_0)\wedge \cdots\wedge
\omg(\sg'_j;\sg_j)\wedge\cdots\wedge\omg(\sg_s)$$ 
d'o\`u le r\'esultat, d'apr\`es le lemme \ref{LEMCOD}. Sinon, 
les prismes $\pi'$ et
$\pi$ s'\'ecrivent sous la forme $\pi'=\tau'\times\sg_0\times \cdots\times \sg_s$ et $\pi
=\tau\times\sg_0\times \cdots\times \sg_s$ o\`u $\tau'$ est une face de codimension 
1 de $\tau$. Dans
ce cas, $\omg (\pi' /\tau';\pi)=\omg (\pi/\tau)$ et sa diff\'erentielle relative est nulle, 
puisque pour chaque $\sigma_i$, on a $d\omg({\sigma_i})  = 0$ dans $\sigma_i$.  
\end{preuve}

\subsection{Triangulations et formes sous-analytiques}

Soit $X\subset {\bf R}^n\times {\bf R}^m$ un sous-ensemble sous-analytique non singulier tel
que la restriction \`a $X$ de la premi\`ere projection soit propre et triangulable.
Soit  $\Delta \subset {\bf R}^n\times {\bf R}^m$ le complexe simplicial lin\'eaire
image r\'eciproque de $X$ par un hom\'eomorphisme sous-analytique 
$$\xymatrix@C=0.3cm{ 
X \subset  {\bf R}^n\times {\bf R}^m   \ar[rd]
&& \Delta \subset  {\bf R}^n\times {\bf R}^m \ar[ld]^f \ar[ll]^t\\&{\bf R}^n  }$$
de triangulation et notons $f$
la restriction \`a $\Delta$ de la premi\`ere projection. Reprenons les notations du
th\'eor\`eme \ref{COMPFA}  et en particulier notons $\psi \colon {\cal P}_f\to {\cal
S}_f$ le morphisme prismal de {\it loc. cit}.

\begin{lemma}\libel{LEIMRE} Soit $\omg$ une $r$-forme diff\'erentielle analytique (ou de
classe ${\cal C}^k$, $k\geq 2$) sur $X$. L'image r\'eciproque $t^*\omg$ est une forme
sous-analytique sur $\Delta$ et l'on a $d(t^*\omg )=t^*d\omg$.
\end{lemma}

\begin{preuve}
L'\'egalit\'e $d(t^*\omg )=t^*d\omg$ est v\'erifi\'ee \`a l'int\'erieur de chaque simplexe  de $\Delta$ puisque
l'hom\'eomorphisme $t$ y est analytique. Remarquons que les  coefficients de $t^*\omg$ et
$t^*d\omg$ sont  sous-analytiques. 

Montrons que la forme $t^*\omg$ est r\'eguli\`ere sur $\Delta$ au sens de Whitney
\cite[Chapitre III, p.104]{Whi}. Cela signifie qu'elle est continue et qu'il existe une
$(r+1)$-forme
$\omg'$ continue sur $\Delta$ telle que, pour tout $(r+1)$-simplexe singulier
$\zeta~:\Sigma^{r+1} \to\Delta$ o\`u $\Sigma^{r+1}$ d\'esigne le simplexe type de
dimension $r+1$, on ait~: $$\int_\zeta \omg' = \int_{\partial \zeta} \omg .$$ D'apr\`es
{\it loc cit.} il suffit de le montrer pour des simplexes singuliers qui sont des
immersions analytiques \`a l'int\'erieur de chaque face, sauf \'eventuellement sur un
ensemble de mesure nulle. 
Pour un tel $\zeta$, le simplexe $t\circ \zeta$ est un simplexe de $X$ poss\'edant les
m\^emes propri\'et\'es. On a alors, en notant encore $\zeta$ l'image de l'application
$\zeta$, les \'egalit\'es suivantes~:
$$\int_{\zeta}  t^*d\omg = \int_{t(\zeta)}  d\omg = \int_{\partial t(\zeta)}  \omg =
\int_{\partial (\zeta)}  t^*\omg $$ ce qui montre que la forme $t^*d\omg$ est la forme
$\omg'$ cherch\'ee. D'apr\`es la g\'en\'eralisation \`a notre situation du lemme 16a de
Whitney \cite[Ch.III]{Whi}, une telle forme, si elle existe, est unique, ce qui entra\^\i ne qu'elle est la
d\'eriv\'ee au sens des distributions de $t^*\omg$ et le r\'esultat.
\end{preuve}

Il r\'esulte de la d\'emonstration de Whitney du th\'eor\`eme de de Rham
\cite[Chapitre IV, \S 29]{Whi} que les formes de Whitney engendrent la cohomologie des formes
$0$-r\'eguli\`eres. Nous verrons plus bas dans le Lemme \ref{LEMENG} la version sous-analytique de ce r\'esultat.

\begin{lemma}\libel{SUBAN}1) Soit $h$ une fonction sous-analytique continue sur un domaine sous-analytique compact $K$ de $\R^k$ contenant l'origine et tout chemin joignant l'origine \`a l'un de ses points. Supposons que $h$ est analytique dans l'int\'erieur de $K$. Pour toute bijection sous-analytique continue $\phi\colon [0,1]\to [0,1]$  v\'erifiant $\phi(0)=0,\ \phi(1)=1$, notant $\phi(s)u$ pour $\phi (s)u_1,\ldots ,\phi(s)u_k$, la fonction 
$$u\mapsto\int_0^1h(\phi(s)u)ds$$
est sous-analytique et continue sur $K$.\par\noindent
2) Si de plus la fonction $h$ est arc-analytique, (voir \cite{B-M1}), il en est de m\^eme de la fonction $u\mapsto\int_0^1h(\phi(s)u)ds$
\end{lemma} 
\begin{preuve}La continuit\'e de la fonction r\'esulte du fait que $h$ est continue. Nous allons utiliser le th\'eor\`eme de rectilin\'earisation des fonctions sous-analytiques (voir \cite{B-M2} Cor. 4.9, Parusi\'nski [Pa]. Appliquons le Th\'eor\`eme 3.4 de [K-P], qui reprend le Th\'eor\`eme 4.1 de [Pa], \`a la fonction $ h\colon K\to \R$ . Nous obtenons une collection finie de morphismes $\pi_\alpha\colon W_\alpha\to \R^k$ telle que:\par\noindent
 $\bullet$ Chaque $W_\alpha$ est analytiquement isomorphe \`a $\R^k$ et contient un compact sous-analytique $K_\alpha$ de telle mani\`ere que $\bigcup_\alpha \pi_\alpha(K_\alpha)=K$,\par\noindent
 $\bullet$ La fonction $h\circ\pi_\alpha$ est analytique dans $W_\alpha$, et dans des coordonn\'ees $y_1,\ldots ,y_k$ sur $\R^k$ les fonctions $u_i\circ\pi_\alpha$ s'\'ecrivent 
 $$u_i\circ\pi_\alpha=
 B_\alpha\, y_1^{e^{(i)}_1}\ldots y_k^{e^{(i)}_k}\ 1\leq i\leq k,$$ o\`u les $e^{(i)}_j$ sont des entiers non n\'egatifs et $B_\alpha$ est une fonction analytique ne s'annulant pas dans $K_\alpha$. Nous pouvons par un changement des coordonn\'ees faire dispara\^\i tre les $B_\alpha$ en absorbant chacune d'elles dans l'une des variables $y_j$.

 Le morphisme $\pi_\alpha$ \'etant g\'en\'eriquement fini, le d\'eterminant de la matrice $({e^{(i)}_j})$ est un entier non nul, 
 disons $d\in {\mathbf Z}$. Il en r\'esulte que nous pouvons trouver des puissances $\psi_{j}(s)=\phi(s)^{R_j}$ de $\phi(s)$, avec des exposants $R_j\in \frac{1}{d}{\mathbf Z}$ uniquement d\'etermin\'es,  telles que l'on ait  \par\noindent
  $$(\phi (s)u_i)\circ\pi_\alpha=(\psi_{1}(s)y_1)^{e^{(i)}_1}\ldots (\psi_{k}(s)y_k)^{e^{(i)}_k}\ 1\leq i\leq k.$$
Choisissons maintenant un point $z$ appartenant \`a  $K_\alpha$, de coordonn\'ees $y_1,\ldots ,y_k$ et tel que $\pi_\alpha(z)=u\in K$. Nous allons calculer la compos\'ee $(\int_0^1h(\phi(s)u)ds)\circ\pi_\alpha$ au voisinage de $z$. \par\noindent
Nous pouvons supposer que pour $s\in [1-a_1,1]$ le chemin $(\psi_1(s)y_1)^{e^{(i)}_1},\ldots ,(\psi_k(s)y_k)^{e^{(i)}_k}$ reste dans $K_\alpha$. Comme certains des $\psi_j(s)$ vont devenir tr\`es grands lorsque $s$ devient petit, le chemin va en g\'en\'eral sortir de $K_\alpha$ avant d'atteindre $s=0$, mais se prolongera analytiquement dans un autre compact $K_\beta$ pour des valeurs $s\in [1-a_2,1-a_1]$, et puisque notre chemin et les $K_\alpha$ sont sous-analytiques, ceci se reproduira un nombre fini de fois. Ainsi notre int\'egrale, au voisinage du point $z$, est la somme d'un nombre fini d'int\'egrales de fonctions analytiques, les $h\circ\pi_\alpha$, le long de chemins analytiques d\'ependant analytiquement du point $z$; c'est donc une fonction analytique sur chacun des $W_\alpha$. Ceci prouve 1). \par\noindent
 Il reste \`a d\'emontrer le point 2). Pour cela, \'etant donn\'e un arc analytique $u_i=u_i(t)$ nous devons v\'erifier que $E=\int_0^1h(\phi(s)u(t))ds$ est analytique en $t$ si la fonction $h(u(t))$ l'est. Par d\'erivation it\'er\'ee sous l'int\'egrale on calcule le d\'eveloppement de Taylor, et l'on constate que ses termes sont le produit des termes correspondants du d\'eveloppement de Taylor de $\int_0^1h(u(t))ds$ par des puissances de $\phi (s)$, qui sont $\leq 1$. Or $\int_0^1h(u(t))ds$ est analytique en $t$ parce que $h(u(t))$ l'est, d'o\`u le r\'esultat.

\end{preuve}
\begin{lemma}\libel{LEMENG} Sur un complexe simplicial, ou plus g\'en\'eralement  prismal, qui est
une vari\'et\'e, les formes de Whitney engendrent la cohomologie  des formes
sous-analytiques.
\end{lemma}

\begin{preuve} Faisons d'abord la d\'emonstration dans le cas d'une d\'ecomposition simpliciale. 
Il suffit de v\'erifier que le lemme
d'int\'egration des formes exactes dans un domaine \'etoil\'e de centre $p_0$ et
les  lemmes de prolongement de Whitney \cite[Chapitre IV, \S\S 25--26]{Whi} sont valables dans le cas
sous-analytique.  Prouvons par exemple que le lemme 25a, p. 136  est
valable dans le cas  sous-analytique. Avec des notations analogues \`a celles de {\it
loc. cit.}, apr\`es avoir pris pour origine le centre $p_0$, il s'agit  de v\'erifier que la forme   
$$\omg_1(p)=\int_{{\bf I}\times p}g^*\omg(s p)ds$$  
est  sous-analytique, o\`u $g\colon {\bf I}\times {\bf R}^n \to 
{\bf R}^n$ est l'application affine d\'efinie par $g(s , p)= sp$. 
Cela r\'esulte du Lemme \ref{SUBAN} en prenant $\phi (s)=s$. Les lemmes de prolongement se v\'erifient de mani\`ere
analogue.\par
 Dans le cas d'une d\'ecomposition prismale d'une vari\'et\'e, 
chaque prisme $\sg_1 \times \ldots \times \sg_s$ admet des d\'ecompositions simpliciales ``standard" 
dont chaque sommet est un $s$-uple $(v_{i_1}, \ldots, v_{i_s})$ de sommets ne 
faisant intervenir que des sommets des $\sg_i$. 
Sur un prisme  la forme $\omg(\sg_1) \wedge \cdots \wedge \omg(\sg_s)$ est la forme volume 
du prisme. Sa restriction \`a chaque simplexe d'une d\'ecomposition ``standard" du prisme
est la forme volume du simplexe. On est alors ramen\'e au cas des d\'ecompositions simpliciales.

\end{preuve}

A une forme sous-analytique ferm\'ee $\omg$ sur $X$ on associe donc une forme
sous-analytique ferm\'ee ferm\'ee $\omg_1$, combinaison de  formes de Whitney de simplexes de
$\Delta$, telle que $t^*\omg = \omg _1 +d\alpha$, o\`u $\alpha$ est une forme
sous-analytique ferm\'ee sur $\Delta$.

\section{Formes sous-analytiques et formes de Whitney relatives}

Pour une application prismale  associ\'ee \`a une application simpliciale $f\colon
\Delta \to T$ comme dans l'exemple \ref{EXEMPF}, \'etant donn\'e un point $x\in \sg$ et $\varepsilon \in ]0,1]$, 
on note  $\sg_{(x,\varepsilon)}$ l'homoth\'etique de $\sg$ par l'homoth\'etie de centre $x$ et de rapport 
$\varepsilon$ (voir Figure 7). 
De m\^eme, si $\phi$ est une face de $\sg$, on note $\phi_{(x,\varepsilon)}$ 
l'homoth\'etique de $\phi$ par l'homoth\'etie de centre $x\in \sg$ et de rapport $\varepsilon$. Enfin, dans le 
faisceau prismal ${\mathcal P}_f$, on note $\pi(\sg)_{(x,\varepsilon)}$ et $\pi(\phi)_{(x,\varepsilon)}$ 
respectivement les homoth\'etiques de $\pi(\sg)$ et $\pi(\phi)$ par homoth\'etie de centre $x$ et de rapport 
$\varepsilon$. 

\begin{figure}
\begin{center}
\begin{tikzpicture}\label{defig7}

\node (M) at (-2,2.5) [coordinate] {} ;
\node [right] at (M) {C} ;
\node (N) at (-2,7.5) [coordinate] {} ;
\node [right] at (N) {D} ;
\node (O) at (-8,2.5) [coordinate] {} ;
\node [left] at (O) {A} ;
\node (P) at (-5.4,0.5) [coordinate] {} ;
\node [below left] at (P) {B} ;

\node (on) at (-5,5) [coordinate] {} ;
\node (pn) at (-3.7,4) [coordinate] {} ;
\node (pm) at (-3.7,1.5) [coordinate] {} ;
\node (om) at (-5,2.5) [coordinate] {} ;

\node (u1) at (-4.15,2.696)  [coordinate] {} ;
\node (u2) at (-4.55,3.004)  [coordinate] {} ;
\node (v1) at (-4.15,3.496)  [coordinate] {} ;
\node (v2) at (-4.55,3.804)  [coordinate] {} ;

\node (v1a) at (-5,1.746)  [coordinate] {} ; 
\node (u1a) at (-1.7, 4.137)  [coordinate] {} ; 

\node (o) at (-5.515,3.004) [coordinate] {} ;
\node [below left] at (o) {a} ;
\coordinate (m) at (intersection of u1--u1a and o--u2);
\node [right] at (m) {c} ;
\coordinate (n) at (intersection of o--v2 and v1--v1a);
\node [above] at (n) {d} ;
\coordinate (p) at (intersection of m--u1 and v1--v1a);
\node [below left] at (p) {b} ;

\coordinate (b1) at (intersection of o--n and on--pn);
\coordinate (b2) at (intersection of p--m and pn--pm);
\coordinate (b3) at (intersection of p--n and on--pn);

\filldraw[fill=gray!30!white, draw=black!50!black] (o) -- (v2) -- (u2) -- (u1) -- (p) -- (o);
\filldraw[fill=gray!30!white, draw=black!50!black] (b1) -- (n) -- (m) -- (b2) -- (pn) -- (b1);
\filldraw[fill=gray!10!white, draw=black!50!black] (v2) -- (b1) -- (pn) -- (b2) -- (u1) -- (v1) -- (v2);

\filldraw[fill=gray!60!white, draw=black!50!black] (u1) -- (v1) -- (v2) -- (u2) -- (u1);

\draw [line width=2pt]  (O) -- (P) -- (M) -- (N) -- (O);
\draw [line width=2pt] (P) -- (N);
\draw [line width=1.5pt]  (on) -- (pn) -- (pm);
\draw[dashed] (O) -- (M);
\draw  [line width=1.3pt]  [densely dashed] (on) -- (om) -- (pm);

\draw (p) -- (v1);
\draw (b3) -- (n);
\draw [densely dashed] (v2) -- (b1);
\draw [densely dashed] (v1) -- (b3);
\draw [densely dashed] (u1) -- (b2);
\draw [densely dashed] (o) -- (m);

\node (A) at (0,1.5) [coordinate] {} ;
\node (B) at (5,1.5) [coordinate] {} ;
\node (C) at (5,6.5) [coordinate] {} ;
\node (D) at (0,6.5) [coordinate] {} ;
\node (E) at (2,1.5) [coordinate] {} ;
\node (F) at (2,6.5) [coordinate] {} ;

\node (G) at (1,2.5) [coordinate] {} ;
\node (H) at (1,7.5) [coordinate] {} ;
\node (J) at (3,0.5) [coordinate] {} ;
\node (K) at (3,5.5) [coordinate] {} ;

\node (a) at (1.5,3.9) [coordinate] {} ;
\node (b) at (2.5,2.9) [coordinate] {} ;
\node (c) at (1.5,5.1) [coordinate] {} ;
\node (d) at (2.5,4.1) [coordinate] {} ;

\draw (G) -- (J) -- (K) -- (H) -- (G);
\filldraw[fill=gray!60!white, draw=black!50!black] (a) -- (b) -- (d) -- (c) -- (a);

\draw (-0.5,-1) -- (4.5,-1);
\draw (3,6) node{$\pi(\sigma) \cap  f^{-1}(f(x))$};
\draw (1.8,4.2) node{$x$};
\filldraw [black] (2,4 ) circle (1.5pt);
\draw (3,3.5) node{$\pi(\sg)_{(x,\varepsilon)} \cap f^{-1}(f(x))$};

\draw[->] (2.5,0.3) -- (2.5,-0.2);
\draw (2.8,0.05) node{$e_{\cal P}$};
\draw (2.5,-0.7) node{$\tau_{(f(x),\varepsilon)}$};
\draw (3.5,-1.3) node{$\tau$};
\draw (2,-1.3) node{$f(x)$};
\filldraw [black] (2,-1) circle (1.5pt);
\draw [line width=1.5pt]  (1,-1) -- (3,-1);
\filldraw [black] (-0.5,-1) circle (1.5pt);
\filldraw [black] (4.5,-1) circle (1.5pt);


\draw (0.5,7.5) -- (1.5,7.5);
\draw[dashed] (0,7.5) -- (0.5,7.5);
\draw[dashed] (1.5,7.5) -- (2,7.5);
\draw (2.5,5.5) -- (3.5,5.5);
\draw[dashed] (2,5.5) -- (2.5,5.5);
\draw[dashed] (3.5,5.5) -- (4,5.5);
\draw (2.5,0.5) -- (3.5,0.5);
\draw[dashed] (2,0.5) -- (2.5,0.5);
\draw[dashed] (3.5,0.5) -- (4,0.5);
\draw (0.5,2.5) -- (1,2.5);
\draw[dashed] (0,2.5) -- (0.5,2.5);
\draw[dashed] (1,2.5) -- (2,2.5);


\draw (-2.5,4.6) node{$\sigma$};
\draw (-5.12,3.2) node{$\sigma_{(x,\varepsilon)}$};
\coordinate (x1) at (intersection of u1--v2 and u2--v1);
\draw (-4.43,3.5) node{$x$};
\filldraw [black] (x1) circle (1.5pt);

\draw (-7,-1) -- (-2,-1);
\draw (-4,-0.7) node{$\tau_{(f(x),\varepsilon)}$};
\draw (-3,-1.3) node{$\tau$};
\draw[->] (-4.5,0.3) -- (-4.5,-0.2);
\draw (-4.8,0.05) node{$e_{\cal S}$};
\draw (-4.5,-1.3) node{$f(x)$};
\draw [line width=1.5pt]  (-5.5,-1) -- (-3.5,-1);
\filldraw [black] (-4.5,-1) circle (1.5pt);
\filldraw [black] (-7,-1) circle (1.5pt);
\filldraw [black] (-2,-1) circle (1.5pt);


\draw[->] (0,3.5) -- (-1.2,3.5);
\draw (-0.6,3.2) node{$\psi^\sigma$};

\filldraw[fill=gray!30!white, draw=black!50!black] (-8,-2) -- (-7,-2) -- (-7,-2.7) -- (-8,-2.7) -- (-8,-2);
\draw (-6.3,-2.5) node{$\sigma_{(x,\varepsilon)}$};
\filldraw[fill=gray!60!white, draw=black!50!black] (-5.2,-2) -- (-4.2,-2) -- (-4.2,-2.7) -- (-5.2,-2.7) -- (-5.2,-2);
\draw (-2.3,-2.5) node{$\sigma_{(x,\varepsilon)} \cap f^{-1}(f(x))$};

\filldraw[fill=gray!60!white, draw=black!50!black] (0.5,-2) -- (1.5,-2) -- (1.5,-2.7) -- (0.5,-2.7) -- (0.5,-2);
\draw (3.65,-2.5) node{$\pi(\sg)_{(x,\varepsilon)} \cap f^{-1}(f(x))$};

\end{tikzpicture}
\end{center}

Ici, $\sigma$ est le t\'etra\`edre $ABCD$ et $\sigma_{(x,\varepsilon)}$ est le t\'etra\`edre $abcd$.
 Pour une face $\phi$ de $\sigma$, par exemple le triangle $ACD$, alors $\phi_{(x,\varepsilon)}$ est le triangle $acd$. 
 
Dans ${\cal P}_f$, on n'a dessin\'e que la fibre de $\pi(\sg)$ au dessus de $f(x)$. 
Le dessin de $\pi(\sg)$ est le produit de $\pi(\sg) \cap  f^{-1}(f(x))$ par $\tau$. 
De m\^eme, on n'a dessin\'e que l'intersection de $\pi(\sg)_{(x,\varepsilon)}$ 
avec la fibre, {\it i.e.} $\pi(\sg)_{(x,\varepsilon)} \cap f^{-1}(f(x))$. Le dessin de $\pi(\sg)_{(x,\varepsilon)}$ est 
le produit de cette intersection par $\tau_{(f(x),\varepsilon)}$.

\caption{Les homoth\'etiques $\sg_{(x,\varepsilon)}$ et $\pi(\sg)_{(x,\varepsilon)}$}
\end{figure}

Si $\eta$ est une forme diff\'erentielle induisant une forme de degr\'e $r$ 
non nulle dans les fibres d'un morphisme simplicial ou prismal $f\colon \phi\to
\tau$, on notera $\eta^{f(x)}$, ou simplement $\eta^f$ s'il n'y a pas d'ambig\"uit\'e, la restriction, 
au sens des formes diff\'erentielles, de
$\eta$ \`a la fibre $f^{-1}(f(x))$. 
Posons $\phi^f_{(x,\varepsilon)} = \phi_{(x,\varepsilon)}\cap f^{-1}(f(x))$. 

D'apr\`es le lemme \ref{LECARE}, l'int\'egrale 
$$\int_{\phi^{f}_{(x,\varepsilon)}} \eta^f$$
ne d\'epend que des coefficients de $\eta$ dans la fibre $f^{-1}(f(x))$. D'autre part, le volume 
euclidien de  $\phi^{f}_{(x,\varepsilon)}$ est \'egal \`a
$$\vol (\phi^{f}_{(x,\varepsilon)} )=t_0^{\vert \phi_0\vert}\cdots t_s^{\vert 
\phi_s\vert}
\int_{\theta^\sg (\phi^{f}_{(x,\varepsilon)} )}\omg(\pi(\phi)/\tau).$$ 
o\`u $\theta^\sg$ est le morphisme d\'efini en \ref{decoord}.

\begin{construc} \libel{eclata}{\rm
Soit $\pi = \tau \times \pi/ \tau = \tau \times \sg_0 \times \cdots  \times \sg_s$ un prisme comme en section 2. 
Rappelons que l'on note $t_0, \ldots , t_s$ les coordonn\'ees barycentriques de $\tau$ et  $\mu_{j,k}$, 
$0\le j\le s$, $0\le k \le \dim \sg_j$ celles de $\sg_0\times\ldots \times\sg_s$ (voir (\ref{decoord})). Dans ce qui suit, nous prendrons 
comme coordonn\'ees cart\'esiennes de chacun des $\sg_j$, les coordonn\'ees 
$\mu_{j,k}$ pour $1\le  k \le \dim \sg_j$, c'est-\`a-dire que l'on pose 
$\mu_{j,0} = 1 - \sum_{k=1}^{\dim \sg_j}\mu_{j,k}$. 
On fait de m\^eme pour le simplexe $\tau$. 
Cela fournit un plongement de $\pi$ dans ${\bf R}^m =
{\bf R}^{\vert \tau \vert} \times {\bf R}^{\vert \sg_0 \vert} \cdots \times {\bf R}^{\vert \sg_s \vert}$ et un plongement
de $\pi / \tau$ dans ${\bf R}^{m - \vert \tau \vert}$. 

Prenons une autre copie de ${\bf R}^m$ dont on notera $u_{j,k}$ les coordonn\'ees correspondant aux 
$\mu_{j,k}$ et $\theta_\ell$ les coordonn\'ees correspondant aux coordonn\'ees $t_\ell$ de $\tau$ et consid\'erons le sous-ensemble $Z$ de $\tau \times \pi/\tau  \times {\bf I} \times {\bf R}^m$ d\'efini 
par les in\'egalit\'es 
\begin{align}\label{inegaga}
- \varepsilon \mu_{j,k} \le u_{j,k} -&  \mu_{j,k} \le \varepsilon (1-  \mu_{j,k})\\
\sum_k (u_{j,k} - \mu_{j,k})& \le \varepsilon (1-  \sum_k  \mu_{j,k})
\end{align}
pour $0\le j\le s$. 
Notons $p$ la restriction \`a $Z$ de la projection de $\pi \times {\bf I} \times {\bf R}^m$
sur $\pi \times {\bf I}$. la fibre de $p$ au dessus du point de coordonn\'ees $(t_\ell, \mu_{j,k}, \varepsilon)$ 
est le prisme $\pi_{(x,\varepsilon)}$. L'intersection de cette fibre avec l'espace lin\'eaire d\'efini par $\theta_\ell=t_\ell,\ 1\le \ell\le s$ est la fibre $\pi^f_{{x,\epsilon}}$ passant par $x$.\par
Notons $b\colon X\to \pi\times {\bf I}$ l'\'eclatement de l'id\'eal $(\varepsilon, (\mu_{j,k})_{j,k})$ dans $\pi\times {\bf I}$. Nous allons montrer que certaines int\'egrales d\'ependant de param\`etres $v\in \pi\times {\bf I}$ deviennent analytiques apr\`es composition avec $b$ et sont donc sous-analytiques sur $\pi\times {\bf I}$. Nous allons d\'etailler le calcul dans une carte particuli\`ere de l'\'eclatement $b$; les calculs dans les autres cartes sont analogues.\par\noindent

Consid\'erons donc l'application $b: \tau \times{\bf R}^{m - \vert \tau \vert}  \times {\bf I} \to 
\tau \times{\bf R}^{m - \vert \tau \vert}  \times {\bf I}$ 
d\'efinie par  
$ \mu_{j,k} \circ b = \varepsilon  \mu'_{j,k}$ pour $k\ge 0$ et $\varepsilon \circ b = \varepsilon$. C'est celle qui correspond \`a la carte de $X$ o\`u $\epsilon$ engendre l'id\'eal transform\'e dans $X$.\par
L'image par $b$ du sous-ensemble de $\tau \times {\bf R}^{m - \vert \tau \vert}_{\ge 0}  \times {\bf I} $ 
d\'efini par les in\'egalit\'es 
$$ 0 \le \varepsilon \mu'_{j,k} \le 1 \quad \text { et } \quad \sum \varepsilon \mu'_{j,k} \le 1$$ est 
$\tau \times \pi/\tau  \times {\bf I}$.

Consid\'erons le diagramme suivant :
$$
\xymatrix @C=0.5cm 
{ 
& \tau\times  {\bf R}^{m - \vert \tau \vert}\times  Y_\alpha  \ar[dd]^{\zeta_\alpha}\\
&& \\
& \tau\times  {\bf R}^{m - \vert \tau \vert}\times S_\alpha \ar[rr]^{} \ar[dd]^{{\widetilde{\pi}}_{\alpha,B}}
& &\tau\times  {\bf R}^{m - \vert \tau \vert}\times  {\bf I} \times W_\alpha \ar[rr]^{} \ar[dd]^{\pi_{\alpha,B}}
 && \pi \times  {\bf I} \times W_\alpha \ar[dd]^{\pi_{\alpha}}\\
 &&&&&\\
& \tau\times  {\bf R}^{m - \vert \tau \vert}\times {\bf I} \times{\bf R}^m  \ar[rr]^-{c} 
&& \tau\times  {\bf R}^{m - \vert \tau \vert}\times {\bf I} \times{\bf R}^m  \ar[rr]^-{B} \ar@{-}[d]^{P_b} 
&& \pi \times  {\bf I} \times {\bf R}^m \ar[dd]^{P}\\
{\widetilde Z}_b \ar[rr]^-{c}  \ar@{^{(}->}[ur]  \ar[rrrd] && Z_b \ar[rr] 
  \ar@{^{(}->}[ur] \ar[rd]^{p_b} 
& \ar[d]& Z \ar@{^{(}->}[ur] \ar[rd]^{p}\\
&&&  \tau\times  {\bf R}^{m - \vert \tau \vert}\times {\bf I} \ar[rr]^{b} &&\pi\times {\bf I}
}
$$
\medskip

Le produit fibr\'e $Z_b$ de $p$ par $b$ est le sous-espace de 
$\tau \times {\bf R}^{m - \vert \tau \vert}\times  {\bf I} \times {\bf R}^m$ 
d\'efini par les in\'egalit\'es 
$$- \varepsilon^2 \mu'_{j,k} \le u_{j,k} - \varepsilon\mu'_{j,k} \le \varepsilon (1-  \varepsilon\mu'_{j,k})$$
$$\sum_k (u_{j,k} - \varepsilon\mu'_{j,k}) \le \varepsilon (1-  
\varepsilon \sum_k \mu'_{j,k}).$$
Notons $p_b$ la projection de $Z_b$ sur $\tau \times {\bf R}^{m - \vert \tau \vert}\times  {\bf I}$. La dimension des fibres de $p_b$ chute encore pour $\varepsilon =0$. Nous allons consid\'erer l'\'eclatement $c\colon V\to{\bf R}^{m - \vert \tau \vert}\times {\bf I} \times{\bf R}^m$ de l'id\'eal $(\varepsilon, (u_{j,k})_{j,k})$. A nouveau nous allons d\'etailler le calcul seulement dans une carte de $V$. \par Consid\'erons donc le morphisme $c$ de $\tau \times {\bf R}^{m - \vert \tau \vert}\times  {\bf I} \times {\bf R}^m$
dans $\tau \times {\bf R}^{m - \vert \tau \vert}\times  {\bf I} \times {\bf R}^m$ d\'etermin\'e par 
$\mu_{j,k} \circ c = \mu'_{j,k}$ et $u_{j,k} \circ c = \varepsilon u'_{j,k}$. Alors $Z_b$ est l'image par le morphisme $c$ du sous-ensemble $\widetilde Z_b$
de $\tau \times {\bf R}^{m - \vert \tau \vert}\times  {\bf I} \times {\bf R}^m$ 
d\'efini par les in\'egalit\'es 
$$- \varepsilon \mu'_{j,k} \le u'_{j,k} - \mu'_{j,k} \le 1-  \varepsilon\mu'_{j,k}$$
$$\sum_k (u'_{j,k} - \mu'_{j,k}) \le 1-  \varepsilon \sum_k \mu'_{j,k}.$$

La fibre du morphisme $p_b \circ c$ au dessus du point de coordonn\'ees $(t_j, \mu'_{j,k},  \varepsilon)$ a pour
image dans $Z$ le simplexe $\pi_{(x,\varepsilon)}$ o\`u $x$ est le point de coordonn\'ees $(t_j, \mu_{j,k})$ 
avec $\mu_{j,k} = \varepsilon  \mu'_{j,k}$. 

Remarquons que toutes les fibres de $P_b \circ c$ restreint \`a $\widetilde Z_b$ ont la m\^eme 
dimension. 
\par
Soient maintenant $\omg$ une forme diff\'erentielle sous-analytique born\'ee sur le prisme $\pi$ et $\pi_\alpha\colon W_\alpha\to \R^m$ des morphismes analytiques compos\'es d'\'eclatements locaux tels que 
chaque forme diff\'erentielle $\pi_\alpha^*\omg$ soit analytique sur $W_\alpha$. Consid\'erons pour chaque $\alpha$ le morphisme induit $\pi\times{\bf I}\times W_\alpha\to\pi\times{\bf I}\times \R^m$, que nous noterons encore $\pi_\alpha$ par abus. Notons $\pi_{\alpha, B}\colon \colon \tau\times \R^{m-\vert\tau\vert}\times{\bf I}\times W_\alpha\to\tau\times \R^{m-\vert\tau\vert}\times{\bf I}\times \R^m$ le morphisme qui s'en d\'eduit par changement de base par le morphisme $B$.
Nous pouvons supposer que chaque $W_\alpha$ contient un compact sous-analytique $K_\alpha$ tel que la r\'eunion des images des $ K_\alpha$ recouvre l'image du prisme $\pi$
 qui se trouve dans $\R^m$. Nous supposerons de plus que les fonctions $u_{j,k}\circ\pi_\alpha$ sont des mon\^omes $n_{j,k}=y_1^{a^{k,j}_1}\ldots y_m^{a^{k,j}_m}$ en des coordonn\'ees $y_1,\ldots, y_m$ sur $W_\alpha$, en absorbant les unit\'es dans les coordonn\'ees et qu'il en est de m\^eme des fonctions $\theta_\ell\circ\pi_\alpha$.\par
 Le produit fibr\'e $\tau\times\R^{m-\vert\tau\vert}\times S_\alpha$ du diagramme ci-dessus a alors pour \'equations dans $\tau\times\R^{m-\vert\tau\vert}\times{\bf I} \times {\bf R}^m\times W_\alpha$ les $\varepsilon u'_{j,k}- n_{j,k}=0$ pour tous les $k,j$.\par
 Remarquons maintenant que les diff\'erences des exposants des mon\^omes de ces \'equations binomiales, qui sont de la forme $(1, 0\ldots,0, 1,0,\ldots,0, -a^{k,j}_1,\ldots, -a^{k,j}_m)$, engendrent un r\'eseau satur\'e. Les \'equations d\'efinissent donc une vari\'et\'e torique affine sur $\R$, qui peut \^etre r\'esolue par un morphisme torique, ou monomial (voir \cite{Te2}, 6.1).
 Cela signifie que l'on peut trouver un morphisme torique de vari\'et\'es toriques non singuli\`eres $Z_\alpha\to \tau\times{\bf I} \times {\bf R}^m\times \R^m$ tel que la transform\'ee stricte 
$Y_\alpha$ de $S_\alpha$ soit non singuli\`ere. Dans une carte locale de $Z_\alpha$ munie de coordonn\'ees $z_1,\ldots, z_{m+1}$ notre domaine d'int\'egration est d\'efini par des in\'egalit\'es de la forme:
 $$- z^E \mu'_{j,k} \le z^{A_{j,k}} - \mu'_{j,k} \le 1-  z^E\mu'_{j,k}$$
$$\sum_k (z^{A_{j,k}} - \mu'_{j,k}) \le  1-  
z^E \sum_k \mu'_{j,k}.$$
Nous int\'egrerons sur ce domaine, dont les in\'equations d\'ependent lin\'eairement des param\`etres $\mu'_{j,k}$ et analytiquement de $\varepsilon$, la restriction aux fibres au dessus de $t\in\tau$ de l'image r\'eciproque de notre forme diff\'erentielle, qui est \`a support dans l'espace analytique non singulier  $Y_\alpha$ et y est analytique.} \end{construc}
\begin{lemma}\libel{LEMETA} Gardons les notations introduites au d\'ebut de cette section. Soit $\sg$ un simplexe orient\'e de $\Delta$ ayant pour image 
$\tau$ et de dimension relative $d$. Soit $\eta$ une $r$-forme diff\'erentielle
sous-analytique sur $\sigma$ dont la restriction aux fibres de
$\sigma\to\tau$ est de degr\'e $r$, avec $r\leq d$. Alors, pour toute face
$\phi$ de $\sg$ de dimension relative $r$ et d'image $\tau$, l'application
$\widetilde A_\phi\colon \sigma
\to {\bf R}$  d\'efinie par  
\[
x\mapsto \widetilde A_{\phi}(x) = 
{\lim_{\varepsilon \to 0}\left(
\frac
{\displaystyle\int_{\phi^{f}_{(x,\varepsilon)}} \eta^f}
{\hbox{\rm vol} (\phi^{f}_{(x,\varepsilon)} )}
\right)} \]
qui ne d\'epend que de la classe de $\eta$ dans les formes relatives, est sous-analytique born\'ee sur $\sigma$.
\end{lemma}
\begin{preuve} 
On a l'\'egalit\'e:
$$\frac
{\displaystyle \int_{\phi^{f}_{(x,\varepsilon)}} \eta^f}
{\hbox{\rm vol} (\phi^{f}_{(x,\varepsilon)} )} =
\frac
{\displaystyle\int_{\pi(\phi)^f_{(x,\varepsilon)}} \psi^*\eta^f}
{t_0^{\vert \phi_0\vert}\cdots t_s^{\vert 
\phi_s\vert}
\hbox{\rm vol} (\pi(\phi)^f_{(x,\varepsilon)} )}$$
Nous allons \'etudier l'int\'egrale $\displaystyle\int_{\pi(\phi)^f_{(x,\varepsilon)}} \psi^*\eta^f $ en utilisant la construction \ref{eclata}. Avec les m\^emes notations, nous avons 
$$\int_{\pi(\phi)^f_{(x,\varepsilon)}} \psi^*\eta^f
 =\int_{p^{-1}(v)\cap \Theta^{-1}(\overline v)} \imath^*_{\overline v} (\eta)
$$
o\`u $v = (x, \varepsilon)$  est le point de $\pi \times {\bf I}$ 
qui a pour coordonn\'ees les $(t_i, \mu_{j,k}, \varepsilon)$, o\`u $\Theta$ d\'esigne le morphisme $\R^m\to \R^{\vert \tau\vert}$ d\'efini par les fonctions $\theta_\ell$ et $\overline v$ d\'esigne l'image de $v$ dans $\tau$. Enfin, $\imath^*_{\overline v} $ d\'esigne la retriction de la forme $\eta$ \`a la fibre de $Z$ au-dessus de $\overline v$ via la restriction \`a $Z$ du morphisme $\Theta$.
 Notons $I(v)$ cette int\'egrale. 

Pour tout point $\tilde v \in b^{-1} (v)$ nous avons
$$I(v)= I(\tilde v) = \int_{p_b^{-1}(\tilde v)\cap\Theta^{-1}(\overline v)} \imath^*_{\overline v}B^* (\eta).$$
De m\^eme, nous pouvons calculer notre int\'egrale $I(v)$ comme int\'egrale sur la fibre $\Theta^{-1}(\overline v)\cap \tilde Z_b$ de $\tilde Z_b$ au-dessus de $\overline v\in \tau$ de l'image inverse $c^*B^* (\eta)$, qui est encore sous-analytique. 
\par 
D'apr\`es la construction \ref{eclata}, pour un choix des $W_\alpha$ rendant analytiques les coefficients de $\eta$ l'int\'egrale $I({\tilde v})$ est la somme d'int\'egrales de formes analytiques restreintes aux fibres du morphisme $\tau \times \R^{m-\vert \tau\vert}\times\R^m\to \tau$ induit par le morphisme $\Theta$ et aux domaines qui sont les images r\'eciproques de $\tilde Z_b$ dans les cartes de $\tau \times \R^{m-\vert \tau\vert}\times  Y_\alpha$.\par Ces domaines d\'ependant analytiquement du param\`etre $v$, on en d\'eduit que la fonction $I(\tilde v)$ est analytique sur la carte de l'\'eclat\'e $X$ que nous avons \'etudi\'ee. Les calculs dans les autres cartes sont analogues et fournissent l'analyticit\'e de $I(v)\circ b$ sur $X$, et donc le fait que la fonction $I(v)$ est sous analytique sur $\pi\times {\bf I}$. 

Pour chaque $t\in \tau$, la forme $\imath^*_t(\eta)$ est un multiple sous-analytique de la forme volume 
$dV_m$ de
${\bf R}^m$. Donc, le quotient $$\frac
{\displaystyle\int_{\pi(\phi)^f_{(x,\varepsilon)}} \psi^*\eta^f}
{t_0^{\vert \phi_0\vert}\cdots t_s^{\vert 
\phi_s\vert}
\hbox{\rm vol} (\pi(\phi)^f_{(x,\varepsilon)} )}$$
est born\'e comme fonction de $x$. De m\^eme, l'int\'egrale
$$\int_{p^{-1}_b(\tilde v)} \imath^*_{\overline v} B^* (dV_m)$$
est sous-analytique. Comme le quotient des deux int\'egrales est born\'e, c'est une fonction sous-analytique
sur $p_b(Z_b)$; elle reste sous-analytique en restriction \`a $\varepsilon = 0$, ce qui montre que la fonction
$\widetilde A_\phi$ est sous-analytique sur $\sg$. 
\end{preuve}

En utilisant les notations de la proposition \ref{IMINVE}, pour toute face $\phi$ de $\sg$, de dimension relative $r$ et d'image $\tau$, 
telle que $\phi_0, \ldots ,\phi_s$ d\'esignent les faces de $\phi$ situ\'ees au dessus des sommets  $y_0, \ldots, y_s$ de $\tau$, 
nous noterons
\begin{equation}\label{nouvaphi}
A_{\phi}(x)=  (-1)^{\alpha(\phi,\nu)}\  
\frac{(r+s)!}{\vert \phi_0\vert !\cdots \vert \phi_s\vert ! \ s!}  \widetilde A_{\phi}(x) .
\end{equation}

\begin{proposition}\libel{LEMETB} Soit $\sg$ un simplexe orient\'e de $\Delta$ ayant
pour image 
$\tau$ et de dimension relative $d$. Soit $\eta$ une $r$-forme diff\'erentielle
sous-analytique sur $\sigma$, avec $r\leq d$. Notant
$(\phi^j)_{j\in J}$ les faces de $\sigma$ de dimension relative $r$  et d'image $\tau$, on a
sur le prisme $\pi (\sg)$ l'\'egalit\'e 
$$(\psi^*df) \wedge \left(\psi^*\eta- \sum_{j\in J}t_0^{\vert
\phi^j_0\vert}\cdots t_s^{\vert \phi^j_s\vert}(A_{\phi^j}\circ\psi) \omg(\pi(\phi^j) /\tau;\pi(\sg)
)\right)=0.$$ 
\end{proposition}

\begin{preuve} Soit $\phi$ une face de $\sg$ de dimension relative $r$  et d'image $\tau$. 

Notons $\overline {\eta} (x)$ la forme diff\'erentielle sous-analytique sur $\phi^{f}_{(x,\varepsilon)}$ dont la valeur en tout point $x'$ de 
$\phi^{f}_{(x,\varepsilon)}$ est \'egale \`a $\widetilde A_{\phi}(x) \cdot \omg (\phi^{f}_{(x,\varepsilon)})$. 

Par  d\'efinition de $\widetilde A_{\phi}$ on a, pour $\varepsilon$ suffisamment petit 
$$
\left\vert
\int_{\phi^{f}_{(x,\varepsilon)}}\eta^{f(x)} (x) - \int_{\phi^{f}_{(x,\varepsilon)}}\overline {\eta}(x) 
\right\vert \le C_\phi (\varepsilon) \varepsilon^r.
$$
o\`u $C(\varepsilon)$ tend vers 0 avec $\varepsilon$. 

En prenant l'image par le morphisme lin\'eaire $\theta^\sg$ d\'efini dans l'exemple \ref{JOINT}, au vu de la Proposition \ref{IMINVE}
et de (\ref{nouvaphi}), on obtient
$$\left|\int_{\theta^\sg (\phi^{f}_{(x,\varepsilon)})}\left(\psi^*\eta -t_0^{\vert \phi_0\vert}\cdots
t_s^{\vert \phi_s\vert} (A_{\phi}\circ\psi)\omg(\pi(\phi)/\tau;\pi(\sg)\right)\right|\le 
C_\phi (\varepsilon)  \varepsilon^r.$$ 

Notons 
$$\beta = \psi^*\eta -\sum_{j\in
J} t_0^{\vert \phi^j_0\vert}\cdots
t_s^{\vert \phi^j_s\vert} (A_{\phi^j}\circ\psi)\omg(\pi(\phi^j)/\tau;\pi(\sg))$$
o\`u $(\phi^j)_{j\in J}$ d\'ecrit l'ensemble des faces de $\sg$ de dimension
relative $r$ et d'image $\tau$.

Remarquons que si $\phi^i$ et $\phi^j$ sont deux faces distinctes de $\sg$ de dimension
relative $r$ et d'image $\tau$,  la restriction de $\omg (\phi^i;\sg)$ \`a ${\phi^j}$ est nulle, et donc  
la restriction de $\omg (\pi(\phi^i)/\tau;\pi(\sg))$  \`a $\pi({\phi^j})$ est nulle. 
Etant donn\'e un prisme $\phi^{f}$ de dimension $r$ contenu dans une fibre de $f$ et dans 
une face de dimension relative $r$ 
nous pouvons donc donner l'estimation suivante ind\'ependante de la face qui le contient  :
\begin{equation}\label{honnete}
\left|\int_{\theta^\sg (\phi^{f}_{(x,\varepsilon)})} \beta\right|\le C(\varepsilon) \varepsilon^r,
\end{equation}
o\`u $C(\varepsilon)$ est la plus grande des constantes $C_\phi (\varepsilon)$ correspondant 
aux   faces  $\phi$ de $\sg$ de dimension relative $r$  et d'image $\tau$.

L'in\'egalit\'e (\ref{honnete})  implique que la restriction de $\beta$ 
\`a chaque $\theta^\sg(\phi^{f}_{(x,\varepsilon)})$ est nulle. D'apr\`es le lemme \ref{LECARE} ceci \'equivaut 
\`a dire que, pour toute coordonn\'ee barycentrique $\lambda_i$ de $\sg$, on a $\beta \wedge d\lambda_i = 0$. 
D'o\`u le r\'esultat puisque les 
$\lambda_i$ forment un syst\`eme de coordonn\'ees dans $\pi(\sg)$. 
\end{preuve}

\section{Primitives relatives de formes diff\'erentielles sous-analytiques}

Notre r\'esultat principal appara\^\i tra en section \ref{princip} comme corollaire du r\'esultat suivant~:

\begin{theorem}\libel{THEOPR} 
Soient $g\colon X\to \R^n$ un
morphisme analytique orient\'e propre et triangulable entre vari\'et\'es analytiques, 
et soit $\omg$ une r-forme diff\'erentielle d\'efinie sur $X$, sous-analytique et continue, 
telle que sa restriction \`a chaque fibre non-singuli\`ere $g^{-1}(y)$ 
soit la diff\'erentielle d'une forme sous-analytique $\xi_y$. 
Il existe une forme diff\'erentielle sous-analytique $\Omega$ sur $X$, de degr\'e $r-1$ 
et,  pour chaque fibre non-singuli\`ere $g^{-1}(y) $, une forme $\alpha_y$ de degr\'e $r-2$
sous-analytique et continue,  telles que l'on ait~:
$$\xi_y - \iota^*_y\Omega= d\alpha_y,$$ o\`u $\iota_y$ d\'esigne le plongement dans $X$ de la fibre $g^{-1}(y)$.
\end{theorem}

Ce th\'eor\`eme exprime le fait que l'on peut remplacer les formes $\xi_y$, primitives de $\omg$ sur chaque
fibre non singuli\`ere, par une forme $\Omega$ qui a l'avantage d'\^etre d\'efinie sur tout $X$, 
mais ceci \`a une forme exacte pr\`es sur chaque fibre non singuli\`ere. L'existence de $\Omega$ sera 
montr\'ee au niveau du faisceau prismal ${\cal P}_f$, en la construisant d'abord dans chaque 
prisme $\pi(\sg)$ de dimension maximum, puis ``verticalement", au dessus de
chaque simplexe $\tau$ de la base, enfin en montrant que les formes ainsi obtenues se recollent ``horizontalement". 
\medskip

\begin{preuve}
Choisissons une triangulation
sous-analytique de $g$ et, reprenant les notations de \ref{debut}, notons $t\colon \Delta \to X$
l'hom\'eomorphisme  sous-analytique correspondant, et $f\colon \Delta \to T$ le
morphisme simplicial orient\'e d\'eduit de $g$.  La forme $t^*\omg$ est sous-analytique sur
$\Delta$.  D'apr\`es la Proposition \ref{LEMETB}, on peut \'ecrire sur le faisceau prismal ${\cal P}_f$ l'\'egalit\'e de formes diff\'erentielles  
\begin{equation} \label{eurostar}
(\psi^*de_{\cal S}^*)\wedge (\psi^*(t^*\omg )-\omg_1)=0,
\end{equation}
 o\`u sur chaque prisme $\pi(\sigma)$, 
notant $\tau$ son image,   la forme $\omg_1$
est \'egale \`a 
$$\omg_1=\sum_{\ell \in L} (A_{\phi^\ell}\circ\psi) t_0^{\vert \phi^\ell_0\vert}\cdots t_s^{\vert 
\phi^\ell_s\vert} \omg(\pi(\phi^\ell)/\tau;\pi (\sg)),$$ 
o\`u les $(\phi^\ell)_{\ell\in L}$ sont les faces de $\sigma$ de dimension relative $r$  et d'image $\tau$, et les  
$A_{\phi^\ell}$ sont des fonctions sous-analytiques d\'efinies sur $\sigma$.
Notons encore $A_{\phi^\ell}$ pour $A_{\phi^\ell} \circ \psi$; par construction, lorsque le prisme 
$\pi(\sigma)$ varie ainsi que son image $\tau$, la  donn\'ee 
$\pi (\sigma) \mapsto A_{\phi^\ell} \omg (\pi(\phi^\ell)/\tau;\pi)$ d\'efinit une
forme diff\'erentielle relative sur ${\cal P}_f$.

Cherchons les conditions pour que la forme diff\'erentielle $\omg _1$  soit exacte dans les 
fibres. Le lemme  \ref{LEMCOD}, b) sugg\`ere de l'\'ecrire comme 
diff\'erentielle d'une combinaison de formes de Whitney. 

Examinons d'abord ce qui se passe dans le prisme $\pi(\sg)=\tau\times \sg_0\times \sg_1\times\cdots\times\sg_s$. 
Les faces de $\pi(\sg)$  de
dimension relative $r$  et ayant pour image $\tau$ sont de la forme 
$$\pi(\phi) = \tau\times \sg_0'\times\sg_1'\times\cdots\times \sg_s'$$ 
o\`u chaque $\sg_i'$ est une face du simplexe $\sg_i$ et o\`u 
la somme des dimensions des simplexes $\sg_i'$ vaut $r$. Chaque face  de $\pi(\phi)$ de 
dimension relative $r-1$  et ayant pour image $\tau$ est obtenue par le proc\'ed\'e suivant : 
On choisit un sommet $y_j$ du simplexe $\tau$ et, dans le simplexe $\sg_j'$ au dessus de $y_j$, un sommet $x_{j,q}$. Si l'on note 
$\sg''_{j,q}$ la face de $\sg_j'$ oppos\'ee au sommet $x_{j,q}$, on obtient une telle face de $\pi(\phi)$ not\'ee
$$\pi(\gamma)=\tau\times \sg'_0\times \cdots \sg'_{j-1}\times \sg''_{j,q}\times \sg'_{j+1}\times\cdots\times\sg'_s.$$

\begin{figure}[h]
\begin{center}
\begin{tikzpicture}[scale=0.85] \label{figure8}

\node (A) at (0,0) [coordinate] {} ;
\node (B) at (-2,1) [coordinate] {} ;
\node (C) at (-0.6,3.7) [coordinate] {} ;
\node (D) at (5,2) [coordinate] {} ;
\node (E) at (-0.8,-2) [coordinate] {} ;
\node (F)  at (4.6,5.7) [coordinate] {} ;

\node (a) at (-1.5,-4) [coordinate] {} ;
\node (b) at (5,-4) [coordinate] {} ;

\draw (A) -- (B) --(C) ;
\draw  (D) -- (E) -- (B);
\draw (E) -- (A) --(C) ;
\draw (A) -- (F) ;
\draw (D) -- (F) --(C) ;
\draw  (A) --(D) ;
\draw  (a) --(b) ;
\draw [dashed] (B) -- (D) -- (C);

\filldraw [black] (a) circle (1.5pt);
\filldraw [black] (b) circle (1.5pt);
\filldraw [black] (A) circle (1.5pt);
\filldraw [black] (B) circle (1.5pt);

\filldraw[color=lightgray,pattern=vertical lines] (A) -- (B)-- (D);

\draw (0.6,1.8) node{$\sg$};
\draw (-1,1.7) node{$\sg_j$};
\draw (-1,0.13) node{$\sg'_j$};
\draw (0.3,-0.2) node{$\sg''_{j,q}$};
\fill[white] (1.5,1.1) circle (0.3);
\draw (1.5,1.1) node{$\phi$};
\draw (2.1,0.5) node{$\gamma$};
\draw (-2.2,0.8) node{$x_{j,q}$};

\draw (-1.5,-3.7) node{$y_j$};
\draw (5,-3.7) node{$y_0$};
\draw (1.7,-3.7) node{$\tau$};

\draw [line width=2pt] (A) -- (B);
\draw [line width=2pt] (A) -- (D);

\node (A1) at (8,0) [coordinate] {} ;
\node (B1) at (6,1) [coordinate] {} ;
\node (C1) at (7.4,3.7) [coordinate] {} ;
\node (E1) at (7.2,-2) [coordinate] {} ;


\node (F1) at (14.5,0) [coordinate] {} ;
\node (G1) at (12.5,1) [coordinate] {} ;
\node (H1) at (13.9,3.7) [coordinate] {} ;
\node (J1) at (13.7,-2) [coordinate] {} ;

\node (K1) at (9.5,1.8) [coordinate] {} ;
\node (L1) at (9.1,5.5) [coordinate] {} ;
\node (M1) at (16,1.8) [coordinate] {} ;
\node (N1) at (15.4,5.5) [coordinate] {} ;

\node (a1) at (6.5,-4) [coordinate] {} ;
\node (b1) at (13,-4) [coordinate] {} ;


\draw (A1) -- (B1) --(C1) ;
\draw (E1) -- (A1) --(C1) ;
\draw  (A1) --(F1) ;
\draw (B1) -- (E1) --(J1) -- (F1);
\draw (A1) -- (K1) --(L1) -- (C1);
\draw (F1) -- (M1) -- (N1) ;
\draw  (K1) --(M1) ;
\draw  (L1) --(N1) ;
\draw  (a1) --(b1) ;
\draw [dashed] (B1) -- (G1) -- (J1);
\draw [dashed] (C1) -- (H1) -- (F1);
\draw [dashed] (G1) -- (H1) -- (N1);

\filldraw [black] (a1) circle (1.5pt);
\filldraw [black] (b1) circle (1.5pt);
\filldraw [black] (A1) circle (1.5pt);
\filldraw [black] (B1) circle (1.5pt);

\filldraw[color=lightgray,pattern=vertical lines] (A1) -- (B1)-- (G1) -- (F1);

\draw (10.8,1.4) node{$\pi(\sg)$};
\draw (7.1,1.7) node{$\widetilde \sg_j$};
\draw (6.9,0.13) node{$\widetilde\sg'_j$};
\draw (8.3,-0.4) node{$\widetilde\sg''_{j,q}$};
\fill[white] (10.1,0.5) circle (0.42);
\draw (10.1,0.5) node{$\pi(\phi)$};
\draw (11.1,-0.3) node{$\pi(\gamma)$};
\draw (5.8,0.8) node{$\widetilde x_{j,q}$};

\draw (6.5,-3.7) node{$y_j$};
\draw (9.7,-3.7) node{$\tau$};

\draw [line width=2pt] (B1) -- (A1) -- (F1);

\end{tikzpicture}
\end{center}
 
\hglue2.5truecm ${\cal S}_f$ au dessus de $\tau$ \hskip 3truecm ${\cal P}_f$ au dessus de $\tau$
\medskip

Dans cette figure,  $\sg$ est le t\'etra\`edre de fibres $\sg_j$ au dessus de $y_j$ et un point au dessus de $y_0$, 
le triangle $\phi$ (hachur\'e) en est une face de dimension relative 1 
et sa fibre au dessus de $y_j$ est $\sg'_j$, face de $\sg_j$. Le segment $\gamma$ 
est une face de codimension 1 de $\phi$, diff\'erant de $\phi$ au dessus du sommet $y_j$ de $\tau$. 
Le sommet $x_{j,q}$  de $\sg'_j$ est le seul sommet de $\phi$ n'appartenant pas \`a $\gamma$.

Dans ${\mathcal P}_f$, le prisme $\pi(\sg)$ est le prisme triangulaire de base $\widetilde \sg_j$. Le ``rectangle" 
$\pi(\phi)$ en est une face.
Les \'el\'ements de  la fibre de ${\mathcal P}_f$ au dessus du sommet  $y_j$ de $\tau$ 
sont not\'es comme suit~: $\widetilde u = \sg_0 \times \cdots \times u \times \cdots \times \sg_s$ 
o\`u $u$ figure en $j$-\`eme position.

 On a \'egalement dessin\'e d'autres simplexes tels que $\sg$ et $\phi$ dont $\gamma$ est une face. 

\caption{Faces $\gamma$ de $\phi$ dans ${\mathcal S}_f$ et ${\mathcal P}_f$}
\end{figure}

La donn\'ee du couple  $(\phi,\gamma)$ d\'etermine donc l'indice $j=j(\phi,\gamma)$ du sommet $y_j$ au dessus duquel 
$\phi$ et $\gamma$ diff\`erent ainsi que l'unique sommet $x_{j,q}$ de $\phi$ non situ\'e dans $\gamma$. 

Ecrire que la forme diff\'erentielle $\omg_1$ est exacte dans les fibres revient \`a r\'esoudre 
dans le prisme $\pi(\sg)$, pour chaque face $\pi(\phi)$ comme ci-dessus, 
l'\'equation
\begin{equation}\label{equaf}
\begin{array} {rl}
t_0^{\vert \phi_0\vert}\cdots t_s^{\vert 
\phi_s\vert} & A_{\phi}  \omg(\pi(\phi)/\tau;\pi (\sg))=\\
& d_e\left(\sum_{j=0}^s t_0^{\vert \phi_0\vert}\cdots  t_{j-1}^{\vert\phi_{j-1}\vert} t_j^{\vert\phi_j\vert -1} 
t_{j+1}^{\vert\phi_{j+1}\vert}\cdots t_s^{\vert 
\phi_s\vert}\left(\sum_\gamma C^\phi_{\gamma}\omg(\pi(\gamma)/\tau;\pi(\sg))\right)\right),
\end{array}\end{equation} 
o\`u la seconde somme porte sur les faces $\gamma$ de codimension $1$ de $\phi$ diff\'erant de $\phi$ 
par un sommet $x_{j,q}$ situ\'e au dessus du sommet $y_j$ de $\tau$. 
Les inconnues $C^\phi_{\gamma}$ sont des fonctions sous-analytiques des coordonn\'ees barycentriques 
de $\tau$ et  de $\phi$, \`a support dans le prisme $\pi(\sg)$.

Le calcul de $d_e (C^\phi_{\gamma} \omg(\pi(\gamma)/\tau;\pi(\sg))$ est une version relative, et dans $\pi(\sg)$,  
du Lemme \ref{satrapaz}. Notons $(\phi_h)_{h\in H}$ l'ensemble des faces de 
dimension relative $r$ de $\sg$ ayant pour image $\tau$ et admettant $\gamma$ pour face de codimension 1. 
En utilisant les notations de Lemme \ref{satrapaz} et celles des coordonn\'ees barycentriques 
(\ref{decoord}), il vient:
\begin{equation}\label{diffaphi}
d_e ( C^\phi_{\gamma} \omg(\pi(\gamma)/\tau;\pi(\sg)) = 
(-1)^{r} \sum_{h} \left( 
(C^{\phi}_{\gamma})_h + \frac{1}{r} \sum_i \lambda_i \frac{\partial (C^{\phi}_{\gamma})_h}{\partial \lambda_i} \right)
\omg(\pi(\phi_h)/\tau;\pi (\sg)),
\end{equation} 
o\`u la somme sur $i$ porte sur les coordonn\'ees barycentriques $\lambda_i$ de $ \gamma$, consid\'er\'ees comme 
coordonn\'ees barycentriques de $\phi_h$, c'est-\`a-dire que la coordonn\'ee barycentrique de $\phi_h$ 
correspondant au sommet de $\phi_h \setminus \gamma$ n'appara\^\i t pas (voir le lemme \ref{satrapaz}).

D'apr\`es les formules (\ref{equaf}) et (\ref{diffaphi}), la d\'etermination des $C^\phi_{\gamma}$ se ram\`ene \`a la r\'esolution pour tout 
prisme $\sigma$ d'image $\tau$, pour toute face $\phi$  de $\sg$ d'image $\tau$ et de dimension relative $r$  
et pour tout sommet $y_j$ de $\tau$, de l'\'equation
$$ t_j A_\phi = \sum_{\gamma \subset \phi} 
\left( C_{\gamma}^\phi + \frac{1}{r} \sum_{\lambda_i \in I(\gamma)} \lambda_i 
\frac{\partial C_{\gamma}^\phi}{\partial \lambda_i } \right)$$
o\`u $\gamma$ d\'ecrit l'ensemble des faces de $\phi$ de codimension $1$ et d'image $\tau$, ne diff\'erant
de $\phi$ qu'au dessus de $y_j$ . Celles-ci 
sont de dimension relative $r-1$ et la seconde somme porte sur les coordonn\'ees barycentriques $\lambda_i$ de $\gamma$, dont l'ensemble est not\'e $I(\gamma)$. 

Notons $n(\phi/\tau)$ le nombre de telles faces $\gamma$ de $\phi$.
Nous allons chercher des solutions de la forme 
$C_{\gamma}^\phi = t_j  \widetilde{C_{\gamma}^\phi }$ 
o\`u  la fonction $\widetilde{C_{\gamma}^\phi }$
est solution de l'\'equation aux d\'eriv\'ees partielles associ\'ee au probl\`eme
\begin{equation}\label{ctilde}
\frac{1}{n(\phi/\tau)} A_\phi =  \widetilde{C_{\gamma}^\phi } + 
\frac{1}{r} \sum_{ \lambda_i\in  I(\gamma)} \lambda_i 
\frac{\partial  \widetilde{C_{\gamma}^\phi}}{\partial \lambda_i}
\end{equation}
co\"\i ncidant avec $A_\phi$ en restriction \`a $\gamma$.

Remarquons que, travaillant, \`a ce niveau 
de la d\'emonstration, notre construction nous assure de l'existence
d'une solution locale. Le fait que l'on puisse trouver une solution globale (dans les fibres au dessus de l'int\'erieur de $\tau$) 
viendra de l'hypoth\`ese d'exactitude de la restriction de la forme $\omega$ aux fibres lisses. 

\subsection{R\'esolution de l'\'equation aux d\'eriv\'ees partielles}

La r\'esolution de l'\'equation (\ref{ctilde}) proc\`ede de la proposition suivante~:

\begin{proposition}\libel{INTEGRATION}1) Sur le simplexe $\sigma=\{\underline u\in \R^k, u_i\geq 0 ,\sum_{i=1}^ku_i\leq 1\}$, l'\'equation
$$E+\frac{1}{r}\sum_{i=1}^ku_i\frac{\partial E}{\partial u_i}=B$$ avec un second membre $B$ sous-analytique et continu sur $\sigma$ et analytique \`a l'int\'erieur, a une unique solution sous-analytique et continue donn\'ee par
$$E=\int_0^1B(s^{\frac{1}{r}}u)ds.$$
Elle est analytique \`a l'int\'erieur de $\sigma$.\par\noindent
2) Si de plus la fonction $B$ est arc-analytique, il en est de m\^eme de la fonction $E$.
\end{proposition}

\begin{preuve} 
Posons $$E=\int_0^1B(s^{\frac{1}{r}}u)ds.$$ Nous avons la suite d'\'egalit\'es
$$\begin{array}{lr}
u_i\frac{\partial E}{\partial u_i}=\int_0^1u_i\frac{\partial B}{\partial u_i}(s^{\frac{1}{r}}u)s^{\frac{1}{r}}ds\\
\\
B(u)=\big[sB(s^{\frac{1}{r}}u)\big]_{s=0}^{s=1}=\int_0^1\frac{\partial}{\partial s}(sB(s^{\frac{1}{r}}u))ds=\int_0^1\big[B(s^{\frac{1}{r}}u))+\frac{s}{r}\sum_{i=1}^ku_i\frac{\partial B}{\partial u_i}(s^{\frac{1}{r}}u)s^{\frac{1}{r}-1}\big]ds,\\
\end{array}$$
ce qui montre que la fonction $E$ est une solution de l'\'equation. Elle est continue puisque $B$ l'est, et sous-analytique sur $\sigma$ d'apr\`es le Lemme \ref{SUBAN} appliqu\'e \`a $\phi (s)=s^{\frac{1}{r}}$.\par\noindent
Prouvons maintenant l'unicit\'e de la solution. Il s'agit de prouver que la seule solution sous-analytique continue de l'\'equation 
$$E+\frac{1}{r}\sum_{i=1}^ku_i\frac{\partial E}{\partial u_i}=0$$
est la fonction nulle.
Si $E$ est solution de cette \'equation, nous avons les \'egalit\'es $$\frac{\partial}{\partial s}E(s^{\frac{1}{r}}u)=\frac{s^{\frac{1}{r}-1}}{r}\sum_{i=1}^ku_i\frac{\partial E}{\partial u_i}=-s^{\frac{1}{r}-1}E.$$
Nous allons en d\'eduire que $E(s^{\frac{1}{r}}u)=C(u){\rm exp}(-rs^{\frac{1}{r}})$ o\`u $C(u)$ est une fonction sous-analytique et continue. Si $E\neq 0$ on en d\'eduit en faisant $s=0$ que $C(u)$ est constante et \'egale \`a $E(0)$ et en faisant $s=1$ que  $E(u)=E(0){\rm exp}(-r)$ ce qui montre que la fonction $E$ est constante et doit \^etre nulle. \par
Enfin, calculant
$$\frac{\partial}{\partial u_i}\int_0^1B(s^{\frac{1}{r}}u)ds=\int_0^1s^{\frac{1}{r}}\frac{\partial B}{\partial u_i}(s^{\frac{1}{r}}u)ds,$$ nous voyons que le d\'eveloppement de Taylor de $E$ converge l\`a o\`u $B$ est analytique. Ceci prouve la premi\`ere partie de la
proposition. 

Le 2) est cons\'equence du 2) du Lemme \ref{SUBAN} en prenant $\phi (s)=s^{\frac{1}{r}}$.
\end{preuve}

Posons 
\begin{equation}\label{chache}
C^\phi_\gamma =
\frac{t_j}{n(\phi/\tau)} \int_0^1 A_\phi (s^{\frac{1}{r}} \lambda) ds
\end{equation}
La forme
\begin{equation}\label{hache}
\sum_{\gamma \subset \phi} C_{\gamma}^\phi  \omg(\pi(\gamma)/\tau;\pi(\sg))
\end{equation}
est donc, relativement \`a $\omg(\pi(\phi)/\tau;\pi(\sg))$, une primitive de (\ref{equaf}) 
si nous consid\'erons chaque $C_{\gamma}^\phi $ comme une fonction des coordonn\'ees barycentriques $\lambda$ 
de $\sg$, ind\'e\-pen\-dan\-te de la coordonn\'ee attach\'ee au sommet de $\phi $ non situ\'e dans $ \gamma$. 
\goodbreak

Remarquons que 
\begin{enumerate}
\item toute autre solution du syst\`eme d'\'equations (\ref{equaf}) est une forme qui diff\`ere de la pr\'ec\'edente 
par une forme ferm\'ee. 
\item la forme (\ref{hache}) est d\'efinie localement (en fait, dans l'``\'etoile de $\pi(\phi)$"). 
\end{enumerate}

\noindent Afin de montrer que l'on obtient globalement une forme primitive relative de $\psi^*(t^*\omg)$,  il nous faut montrer~:
\begin{description}
\item{a)} que l'on peut construire \`a l'aide de (\ref{hache}) une forme d\'efinie au dessus de l'int\'erieur de $\tau$, 
et qui soit primitive relative de $\psi^*(t^*\omg)$ (prolongement vertical),
\item{b)} que la forme ainsi d\'efinie se sp\'ecialise ``correctement"  au dessus des faces de $\tau$ 
de fa\c con \`a d\'eterminer 
une forme dans ${\mathcal S}_f$ (prolongement horizontal).
\end{description}

\subsection{Prolongement vertical}

On peut supposer $T$ triangul\'e de telle fa\c con que les fibres singuli\`eres de $f: \Delta \to T$ 
soient situ\'ees au dessus du squelette de codimension $1$ de $T$. 
D'apr\`es l'hypoth\`ese du th\'eor\`eme \ref{THEOPR}, la restriction de $t^*\omg$ \`a toute fibre lisse 
de $f$ est exacte. 
Au dessus de tout point $y$ de l'int\'erieur d'un simplexe $\tau$ de dimension maximum, on peut donc \'ecrire dans la fibre 
lisse $F_y = f^{-1}(y)$, en notant $i_y$ son inclusion dans $\Delta$
$$i^*_yt^*\omg = d \zeta_y$$
o\`u la $(r-1)$-forme $\zeta_y$ est repr\'esentable par une forme 
$$\zeta_y = \sum_{\sg,\gamma_k} b_k \, \omg(\gamma_k / \tau;\sigma),$$
la somme portant sur les simplexes $\sg$ de dimension maximum  au dessus de $\tau$ et sur les 
simplexes $\gamma_k $ de dimension relative $(r-1)$ et d'image $\tau$. 
Avec les notations pr\'ec\'edentes, il vient~:
$$\psi^* \zeta_y = \sum_{\sg,\gamma_k} 
 t_0^{\vert \gamma_{k,0} \vert } \ldots t_s^{\vert \gamma_{k,s} \vert } (b_k \circ \psi) \,
\omg(\pi(\gamma_k) / \tau;\pi(\sigma))$$
qui est une forme d\'efinie dans la fibre de ${\mathcal P}_f$ au dessus de $y$.

D'apr\`es (\ref{equaf}) et le calcul pr\'ec\'edent, sur chaque prisme $\pi(\sg)$, 
la forme $\psi^* (t^*\omg\vert_{F_y}) = \psi^* (d \zeta_y)$ s'\'ecrit  
\begin{equation} 
\begin{array} {lr} 
d \left(\sum_{\sg,\gamma_k}  t_0^{\vert \gamma_{k,0} \vert } \ldots t_s^{\vert \gamma_{k,s} \vert } (b_k \circ \psi) \,
\omg(\pi(\gamma_k) / \tau;\pi(\sigma))\right)= \\
\qquad\qquad d \left(\sum_{j=0}^{s} t_0^{\vert \phi_0\vert}\cdots  t_{j-1}^{\vert\phi_{j-1}\vert} t_j^{\vert\phi_j\vert -1} 
t_{j+1}^{\vert\phi_{j+1}\vert}\cdots t_s^{\vert  \phi_s\vert}
\left(\sum_\gamma C^\phi_\gamma \omg(\pi(\gamma)/\tau;\pi(\sg))\right)\right).
\end{array}
\end{equation}
o\`u chaque face $\gamma_k$ de dimension relative $(r-1)$ et d'image $\tau$ appara\^\i t 
une seule fois  de chaque c\^ot\'e de l'\'equation. Remarquons cependant  que si $\gamma_{k}$ 
(dans la somme de gauche) co\"\i ncide  avec $\gamma$ (dans la somme de droite), alors $\vert \gamma_{k,i} \vert = 
\vert \phi_i\vert$ pour $i\ne j$ et $\vert \gamma_{k,j} \vert = \vert \phi_j\vert -1$, 
autrement dit les 
deux termes $t_0^{\vert \gamma_{k,0} \vert } \ldots t_s^{\vert \gamma_{k,s} \vert } $ et 
$t_0^{\vert \phi_0\vert}\cdots  t_{j-1}^{\vert\phi_{j-1}\vert} t_j^{\vert\phi_j\vert -1} 
t_{j+1}^{\vert\phi_{j+1}\vert}\cdots t_s^{\vert  \phi_s\vert}$ sont \'egaux.

D'apr\`es le lemme de Poincar\'e appliqu\'e au domaine \'etoil\'e $\pi(\sg)\cap f^{-1}(y)$
(fibre non singuli\`ere), et avec des notations \'evidentes, 
il existe une $(r-2)$-forme $\alpha_y$ sur $\pi(\sg) \cap f^{-1}(y)$ 
 telle que~:
$$\sum_{\gamma_k} \left( t_0^{\vert \gamma_{k,0}\vert}\cdots  t_s^{\vert  \gamma_{k,s}\vert} \Big[
b_{k} \circ \psi - C^\phi_{\gamma_k} \Big]\; 
\omg(\pi(\gamma_{k})/\tau;\pi(\sg)) \right) = d\alpha_y.$$
Comme pr\'ec\'edemment ({\it cf.} \ref{equaf}), nous allons chercher des solutions $\alpha_y$ sous la forme 
$$\alpha_y = \sum_{h=0}^s \left(t_0^{\vert \gamma_{k,0}\vert}\cdots  t_{h-1}^{\vert\gamma_{k,h-1}\vert} t_h^{\vert\gamma_{k,h}\vert -1} 
t_{h+1}^{\vert\gamma_{k,h+1}\vert}\cdots t_s^{\vert  \gamma_{k,s}\vert}
\sum_\beta D^{\gamma_k}_{\beta} \omg(\pi(\beta)/\tau;\pi(\sg))\right),
$$
o\`u $\beta$ d\'ecrit l'ensemble des  faces de codimension 1 de $\gamma_k$ ne diff\'erant de $\gamma_k$ qu'au dessus du 
sommet $y_h$ de $\tau$.

Nous sommes donc ramen\'es \`a r\'esoudre, pour chaque $\pi(\gamma_{k})$ l'\'equation
\begin{equation} \label{equazz}
\begin{array} {lr}
 t_0^{\vert \gamma_{k,0}\vert}\cdots  t_s^{\vert  \gamma_{k,s}\vert} \Big[
b_{k} \circ \psi -  C^\phi_{\gamma_k} \Big]\; 
\omg(\pi(\gamma_{k})/\tau;\pi(\sg))) = \\
\qquad d \left( \sum_{h=0}^s t_0^{\vert \gamma_{k,0}\vert}\cdots  t_{h-1}^{\vert\gamma_{k,h-1}\vert} t_h^{\vert\gamma_{k,h}\vert -1} 
t_{h+1}^{\vert\gamma_{k,h+1}\vert}\cdots t_s^{\vert  \gamma_{k,s}\vert}
\left( \sum_\beta D^{\gamma_k}_{\beta} \omg(\pi(\beta)/\tau;\pi(\sg))\right)\right),
\end{array}
\end{equation}
laquelle est du m\^eme type que l'\'equation (\ref{equaf}). 
Les solutions en $D^{\gamma_k}_\beta$ sont donc donn\'ees par la proposition 
\ref{INTEGRATION}. Elles sont sous-analytiques, la forme correspondante
$$\sum_{\beta \subset \gamma_{k}} D^{\gamma_{k}}_\beta \omg(\pi(\beta)/\tau ; \pi(\sigma))$$
est sous-analytique et 
d\'etermin\'ee \`a une forme ferm\'ee pr\`es. 

Posons, comme pr\'ec\'edemment, 
$$C_{\gamma}^\phi = t_j  \widetilde{C_{\gamma}^\phi }\qquad {D}^{\gamma}_\beta = t_h \widetilde{{D}^{\gamma}_\beta}.$$
Dans le faisceau prismal ${\cal P}_f$, et dans chaque fibre lisse au dessus de l'int\'erieur 
d'un simplexe $\tau$ de dimension maximale, on d\'efinit  la forme
\begin{equation} \label{equazu}
\begin{array} {lr}
H = \\ 
\sum_{\gamma\subset \phi\subset \sg} 
 t_0^{\vert \phi_0\vert}\cdots t_s^{\vert  \phi_s\vert}
\widetilde{C}_{\gamma}^\phi  \omg(\pi(\gamma)/\tau;\pi(\sg)) + d \left(
\sum_{\beta \subset \gamma\subset \phi\subset \sg} 
 t_0^{\vert \gamma_0\vert}\cdots t_s^{\vert 
\gamma_s\vert}
\widetilde{D}^{\gamma}_\beta  \omg(\pi(\beta)/\tau;\pi(\sg)) 
\right) \\
\qquad = \sum_{\gamma\subset \phi\subset \sg} 
 t_0^{\vert \phi_0\vert}\cdots t_s^{\vert  \phi_s\vert}
\left(
\widetilde{C}_{\gamma}^\phi  \omg(\pi(\gamma)/\tau;\pi(\sg)) + d \left(
\sum_{\beta \subset \gamma} 
\widetilde{D}^{\gamma}_\beta  \omg(\pi(\beta)/\tau;\pi(\sg)) 
\right) \right),
\end{array}
\end{equation}
o\`u les faces $\phi, \gamma, \beta$ de $\sg$ ont toutes pour image $\tau$ 
et sont de dimension relatives respectives $r, r-1, r-2$.

La forme $H$ est d\'efinie sur le faisceau prismal ${\mathcal P}_f$ au dessus de tout simplexe (ferm\'e) $\tau$. 
D'apr\`es la proposition \ref{IMINVE}, 
$H$ est image r\'eciproque d'une forme diff\'erentielle $H_{\cal S}$ d\'efinie sur le faisceau 
prismal ${\mathcal S}_f$ au dessus de $\tau$~:

$$H_{\cal S} = \sum_{\gamma\subset \phi\subset \sg} \left(
\widetilde{C}_{\gamma}^\phi  \omg(\gamma/\tau;\sg) + d \left(
\sum_{\beta \subset \gamma} 
\widetilde{D}^{\gamma}_\beta  \omg(\beta/\tau;\sg) 
\right) \right) .$$

Il en est de m\^eme de $\alpha_y$ au dessus des points $y$ de l'int\'erieur de $\tau$, laquelle
s'\'ecrit $\alpha_y = \psi^*(\alpha_{{\mathcal S},y})$. 
Au dessus d'un tel point, autrement dit pour toute fibre non singuli\`ere $F_y=g^{-1}(y)$, il vient donc~:
$$\psi^* \zeta_y - \psi^*(\iota_y^*(H_{\cal S})) = d \psi^*(\alpha_{{\mathcal S},y}),$$
 o\`u $\iota_y$ d\'esigne le plongement dans $X$ de la fibre $g^{-1}(y)$.

\subsection{Prolongement horizontal}
Soit alors $\tau'$ une face de $\tau$ et $\phi'$ la face de $\phi$ situ\'ee au dessus de $\tau'$. 
Notons $y_0, \ldots , y_u$ les sommets de $\tau'$ et donc $\phi_0, \ldots , \phi_u$ les faces de $\phi'$ situ\'ees 
au dessus de ces sommets. Supposons dans un premier temps   $\dim_{\rm rel} \phi' < r$. Cela implique
$$\vert \phi_0\vert + \cdots + \vert \phi_u \vert - u < r.$$
Mais comme $\vert \phi_0\vert + \cdots + \vert \phi_s \vert -s = r$, il vient
$$\vert \phi_{u+1} \vert + \cdots + \vert \phi_s \vert - (s-u) > 0$$
ce qui signifie que l'on a $\vert \phi_j\vert \ge 0$ pour l'un  au moins des $j = (u+1), \ldots ,s$. On en conclut que 
le coefficient $ t_0^{\vert \phi_0\vert}\cdots t_s^{\vert  \phi_s\vert}$ tend vers $0$ lorsqu'on s'approche de $\phi'$ 
et donc la forme (\ref{equazu}) s'y annule. 

Supposons maintenant que $\dim_{\rm rel} \phi' \ge r$, cela signifie que 
$$\vert \phi_0\vert + \cdots + \vert \phi_u \vert - u \ge r$$
et donc certains des $\vert \phi_j \vert$ pour $j>u$ peuvent s'annuler. En particulier, pour 
$\vert \phi_{u+1} \vert = \cdots= \vert \phi_s \vert =0$, le prisme au dessus 
d'une fibre de $\tau'$ co\"\i ncide avec le prisme au dessus d'une fibre de $\tau$. 
La d\'efinition (formule \ref{chache}) de $C^\phi_\gamma$ en fonction de $A_\phi$ int\`egre
les signes et coefficients  de la formule \ref{nouvaphi}. Ce sont, d'apr\`es la proposition \ref{IMINVE}, 
les signes et coefficients n\'ecessaires pour que la forme diff\'erentielle (\ref{equazu}) se sp\'ecialise correctement et 
d\'efinisse une forme sur le faisceau prismal $\mathcal S_f$.

Dans tous les cas, la forme que l'on peut d\'efinir sur $\tau'$ par le m\^eme proc\'ed\'e co\"\i ncide donc avec 
la sp\'ecialisation de la forme d\'efinie sur $\tau$.

Comme de plus les exposants de 
$ t_0^{\vert \phi_0\vert}\cdots t_s^{\vert  \phi_s\vert}$ correspondent aux dimensions convenables, 
nous obtenons donc une forme diff\'erentielle $H_{\cal S}$ d\'efinie sur ${\cal S}_f$ et satisfaisant~:
$$\zeta_y - \iota_y^*(H_{\cal S}) = d (\alpha_{{\mathcal S},y})$$
pour toute fibre non singuli\`ere $g^{-1}(y)$. Par l'hom\'eomorphisme sous analytique $t$, on en d\'eduit 
le th\'eor\`eme \ref{THEOPR}.
\end{preuve}

\subsection{Le r\'esultat}\label{princip}

\begin{corollary}\libel{THEODG} Soient $g\colon X\to \R^n$ un
morphisme analytique orient\'e propre et triangulable entre vari\'et\'es analytiques, 
et soit $\omg$ une r-forme diff\'erentielle sur $X$ sous-analytique continue 
telle que la restriction de $\omg$ \`a chaque fibre non-singuli\`ere de
$g$ soit la diff\'erentielle d'une forme sous-analytique. 
Il existe une $(r-1)$-forme \san continue $\Omega$ sur $X$ telle
que l'on ait  $$dg\wedge (\omg -d\Omega)=0$$ 
o\`u $dg$ est l'image inverse de la forme volume sur $\R^n$.
\end{corollary}

\begin{preuve}
Par hypoth\`ese, la restriction de $\omg$ \`a toute fibre non-singuli\`ere de $g$ 
s'\'ecrit: 
$$\omg_y = d\xi_y.$$
D'apr\`es le th\'eor\`eme \ref{THEOPR}, il existe une forme diff\'erentielle $(r-1)$-forme diff\'erentielle 
sous-analytique continue $\Omega$ telle que l'on ait 
$$\xi_y - \iota^*_y\Omega= d\alpha_y$$
pour toute fibre non singuli\`ere de $g$, o\`u $\iota_y$ d\'esigne le plongement dans $X$ de la fibre $g^{-1}(y)$. 
La diff\'erentielle de $\Omega$ au sens des 
distributions admet un repr\'esentant sous-analytique continu et qui v\'erifie (voir le lemme \ref{LECARE})
$$dg \wedge (\omg - d\Omega) =0.$$
Ceci d\'emontre le corollaire. 

\end{preuve}


\end{document}